\documentclass[letterpaper,11pt,oneside,reqno]{amsart}
\usepackage{bm}
\usepackage{mathrsfs}
\usepackage{amsfonts,amsmath, amssymb,amsthm,amscd,stmaryrd}
\usepackage[height=9.6in,width=5.95in]{geometry}
\usepackage[mpexclude,DIV13]{typearea}
\usepackage{verbatim}
\usepackage{hyperref}
\usepackage{graphicx}
\usepackage[latin1]{inputenc}
\usepackage{latexsym}
\usepackage{lscape}
\usepackage{epsfig}
\include{bibtex}

\usepackage{setspace}

\usepackage{parskip}
\begin{document}
\bibliographystyle{alpha}
\newcommand{\cn}[1]{\overline{#1}}
\newcommand{\e}[0]{\epsilon}
\newcommand{\bbf}[0]{\mathbf}

\newcommand{\Pfree}[5]{\ensuremath{\mathbb{P}^{#1,#2,#3,#4,#5}}}
\newcommand{\PfreeShort}{\ensuremath{\mathbb{P}^{BB}}}

\newcommand{\WH}[8]{\ensuremath{\mathbb{W}^{#1,#2,#3,#4,#5,#6,#7}_{#8}}}
\newcommand{\Wfree}[5]{\ensuremath{\mathbb{W}^{#1,#2,#3,#4,#5}}}
\newcommand{\WHShort}[3]{\ensuremath{\mathbb{W}^{#1,#2}_{#3}}}
\newcommand{\WHShortCouple}[2]{\ensuremath{\mathbb{W}^{#1}_{#2}}}

\newcommand{\walk}[3]{\ensuremath{X^{#1,#2}_{#3}}}
\newcommand{\walkupdated}[3]{\ensuremath{\tilde{X}^{#1,#2}_{#3}}}
\newcommand{\walkfull}[2]{\ensuremath{X^{#1,#2}}}
\newcommand{\walkfullupdated}[2]{\ensuremath{\tilde{X}^{#1,#2}}}

\newcommand{\PH}[8]{\ensuremath{\mathbb{Q}^{#1,#2,#3,#4,#5,#6,#7}_{#8}}}
\newcommand{\PHShort}[1]{\ensuremath{\mathbb{Q}_{#1}}}
\newcommand{\PHExp}[8]{\ensuremath{\mathbb{F}^{#1,#2,#3,#4,#5,#6,#7}_{#8}}}

\newcommand{\D}[8]{\ensuremath{D^{#1,#2,#3,#4,#5,#6,#7}_{#8}}}
\newcommand{\DShort}[1]{\ensuremath{D_{#1}}}
\newcommand{\partfunc}[8]{\ensuremath{Z^{#1,#2,#3,#4,#5,#6,#7}_{#8}}}
\newcommand{\partfuncShort}[1]{\ensuremath{Z_{#1}}}
\newcommand{\bolt}[8]{\ensuremath{W^{#1,#2,#3,#4,#5,#6,#7}_{#8}}}
\newcommand{\boltShort}[1]{\ensuremath{W_{#1}}}
\newcommand{\boltNew}{\ensuremath{W}}
\newcommand{\QTLH}{\ensuremath{\mathfrak{H}}}
\newcommand{\QTLHgen}{\ensuremath{\mathfrak{L}}}

\newcommand{\whitenoise}{\ensuremath{\mathscr{\dot{W}}}}
\newcommand{\mf}{\mathfrak}

\newcommand{\EE}{\ensuremath{\mathbb{E}}}
\newcommand{\PP}{\ensuremath{\mathbb{P}}}
\newcommand{\var}{\textrm{var}}
\newcommand{\N}{\ensuremath{\mathbb{N}}}
\newcommand{\R}{\ensuremath{\mathbb{R}}}
\newcommand{\C}{\ensuremath{\mathbb{C}}}
\newcommand{\Z}{\ensuremath{\mathbb{Z}}}
\newcommand{\Q}{\ensuremath{\mathbb{Q}}}
\newcommand{\T}{\ensuremath{\mathbb{T}}}
\newcommand{\E}[0]{\mathbb{E}}
\newcommand{\OO}[0]{\Omega}
\newcommand{\F}[0]{\mathfrak{F}}
\def \Ai {{\rm Ai}}
\newcommand{\G}[0]{\mathfrak{G}}
\newcommand{\ta}[0]{\theta}
\newcommand{\w}[0]{\omega}
\newcommand{\ra}[0]{\rightarrow}
\newcommand{\vectoro}{\overline}
\newcommand{\crairy}{\mathcal{CA}}
\newcommand{\nc}{\mathsf{NoTouch}}
\newcommand{\ncf}{\mathsf{NoTouch}^f}
\newcommand{\wxy}{\mathcal{W}_{k;\bar{x},\bar{y}}}
\newcommand{\AP}{\mathfrak{a}}
\newcommand{\cm}{\mathfrak{c}}
\newtheorem{theorem}{Theorem}[section]
\newtheorem{partialtheorem}{Partial Theorem}[section]
\newtheorem{conj}[theorem]{Conjecture}
\newtheorem{lemma}[theorem]{Lemma}
\newtheorem{proposition}[theorem]{Proposition}
\newtheorem{corollary}[theorem]{Corollary}
\newtheorem{claim}[theorem]{Claim}
\newtheorem{experiment}[theorem]{Experimental Result}

\def\todo#1{\marginpar{\raggedright\footnotesize #1}}
\def\change#1{{\color{green}\todo{change}#1}}
\def\note#1{\textup{\textsf{\color{blue}(#1)}}}

\theoremstyle{definition}
\newtheorem{rem}[theorem]{Remark}

\theoremstyle{definition}
\newtheorem{com}[theorem]{Comment}

\theoremstyle{definition}
\newtheorem{definition}[theorem]{Definition}

\theoremstyle{definition}
\newtheorem{definitions}[theorem]{Definitions}

\theoremstyle{definition}
\newtheorem{conjecture}[theorem]{Conjecture}

\newcommand{\airysh}{\mathcal{A}}
\newcommand{\hfixed}{\mathcal{H}}
\newcommand{\afixed}{\mathcal{A}}
\newcommand{\canopynoarg}{\mathsf{C}}
\newcommand{\canopy}[3]{\ensuremath{\mathsf{C}_{#1,#2}^{#3}}}
\newcommand{\argmax}{x_{{\rm max}}}
\newcommand{\zmax}{z_{{\rm max}}}

\newcommand{\Rkle}{\ensuremath{\mathbb{R}^k_{>}}}
\newcommand{\Ronele}{\ensuremath{\mathbb{R}^k_{>}}}
\newcommand{\ewxy}{\mathcal{E}_{k;\bar{x},\bar{y}}}

\newcommand{\bxyf}{\mathcal{B}_{\bar{x},\bar{y},f}}
\newcommand{\bxyflr}{\mathcal{B}_{\bar{x},\bar{y},f}^{\ell,r}}

\newcommand{\bxyfone}{\mathcal{B}_{x_1,y_1,f}}

\newcommand{\ptac}{p}
\newcommand{\ptact}{v}

\newcommand{\fext}{\mathfrak{F}_{{\rm ext}}}
\newcommand{\gext}{\mathfrak{G}_{{\rm ext}}}
\newcommand{\xext}{{\rm xExt}(\mathfrak{c}_+)}

\newcommand{\dd}{\, {\rm d}}
\newcommand{\signc}{\Sigma}
\newcommand{\wxylr}{\mathcal{W}_{k;\bar{x},\bar{y}}^{\ell,r}}
\newcommand{\wxylrprime}{\mathcal{W}_{k;\bar{x}',\bar{y}'}^{\ell,r}}
\newcommand{\Rklezero}{\ensuremath{\mathbb{R}^k_{>0}}}
\newcommand{\XYfM}{\textrm{XY}^{f}_M}

\newcommand{\upright}{D}
\newcommand{\staircase}{SC}
\newcommand{\energy}{E}
\newcommand{\xmax}{{\rm max}_1}
\newcommand{\ymax}{{\rm max}_2}
\newcommand{\lppls}{\mathcal{L}}
\newcommand{\lpplsre}{\mathcal{L}^{{\rm re}}}
\newcommand{\lpplsarg}[1]{\mathcal{L}_{n}^{\fa \to #1}}
\newcommand{\larg}[3]{\mathcal{L}_{n}^{#1,#2;#3}}
\newcommand{\BP}{M}
\newcommand{\weight}{\mathsf{Wgt}}
\newcommand{\pairweight}{\mathsf{PairWgt}}
\newcommand{\sumweight}{\mathsf{SumWgt}}
\newcommand{\mpgood}{\mathcal{G}}
\newcommand{\mpg}{\mathsf{Fav}}
\newcommand{\mcgone}{\mathsf{Fav}_1}
\newcommand{\radnik}[2]{\mathsf{RN}_{#1,#2}}
\newcommand{\size}[2]{\mathsf{S}_{#1,#2}}
\newcommand{\pdr}{\mathsf{PolyDevReg}}
\newcommand{\pwr}{\mathsf{PolyWgtReg}}
\newcommand{\lwr}{\mathsf{LocWgtReg}}
\newcommand{\hwp}{\mathsf{HighWgtPoly}}
\newcommand{\fbr}{\mathsf{ForBouqReg}}
\newcommand{\bbr}{\mathsf{BackBouqReg}}
\newcommand{\fsc}{\mathsf{FavSurCon}}
\newcommand{\maxpoly}{\mathrm{MaxDisjtPoly}}
\newcommand{\maxswf}{\mathsf{MaxScSumWgtFl}}
\newcommand{\emaxswf}{\e \! - \! \maxswf}
\newcommand{\minswf}{\mathsf{MinScSumWgtFl}}
\newcommand{\eminswf}{\e \! - \! \minswf}
\newcommand{\surreg}{\mathcal{R}}
\newcommand{\scf}{\mathsf{FavSurgCond}}
\newcommand{\disjtpoly}{\mathsf{DisjtPoly}}
\newcommand{\intint}[1]{\llbracket 1,#1 \rrbracket}
\newcommand{\maxsym}{*}
\newcommand{\polynum}{\#\mathsf{Poly}}
\newcommand{\dlp}{\mathsf{DisjtLinePoly}}
\newcommand{\lowb}{\underline{B}}
\newcommand{\highb}{\overline{B}}
\newcommand{\tottt}{t_{1,2}^{2/3}}
\newcommand{\tot}{t_{1,2}}
\newcommand{\btone}{{\bf{t}}_1}
\newcommand{\bttwo}{{\bf{t}}_2}
\newcommand{\formerE}{C}
\newcommand{\rcon}{r_0}
\newcommand{\para}{Q}
\newcommand{\Cstrong}{E}

\newcommand{\mc}{\mathcal}
\newcommand{\vect}{\mathbf}
\newcommand{\bt}{\mathbf{t}}
\newcommand{\scB}{\mathscr{B}}
\newcommand{\scBres}{\mathscr{B}^{\mathrm{re}}}
\newcommand{\rightshadow}{\mathrm{RS}Z}
\newcommand{\dbm}{D}
\newcommand{\edgedbm}{D^{\rm edge}}
\newcommand{\gue}{\mathrm{GUE}}
\newcommand{\edgegue}{\mathrm{GUE}^{\mathrm{edge}}}
\newcommand{\eqdist}{\stackrel{(d)}{=}}
\newcommand{\geqdist}{\stackrel{(d)}{\succeq}}
\newcommand{\leqdist}{\stackrel{(d)}{\preceq}}
\newcommand{\scal}{{\rm sc}}
\newcommand{\fa}{x_0}
\newcommand{\hit}{H}
\newcommand{\scaledle}{\mathsf{Nr}\mc{L}}
\newcommand{\cenleup}{\mathscr{L}^{\uparrow}}
\newcommand{\cenledown}{\mathscr{L}^{\downarrow}}
\newcommand{\eln}{T}
\newcommand{\xmin}{{\rm Corner}^{\mfl,\mc{F}}}
\newcommand{\ymin}{{\rm Corner}^{\mfr,\mc{F}}}
\newcommand{\barxmin}{\overline{\rm Corner}^{\mfl,\mc{F}}}
\newcommand{\barymin}{\overline{\rm Corner}^{\mfr,\mc{F}}}
\newcommand{\qmin}{Q^{\mc{F}^1}}
\newcommand{\barqmin}{\bar{Q}^{\mc{F}^1}}
\newcommand{\test}{T}
\newcommand{\mfl}{\mf{l}}
\newcommand{\mfr}{\mf{r}}
\newcommand{\gfl}{\ell}
\newcommand{\gfr}{r}
\newcommand{\jre}{J}
\newcommand{\highfl}{{\rm HFL}}
\newcommand{\flyleap}{\mathsf{FlyLeap}}
\newcommand{\touch}{\mathsf{Touch}}
\newcommand{\notouch}{\mathsf{NoTouch}}
\newcommand{\close}{\mathsf{Close}}
\newcommand{\abovepar}{\mathsf{High}}
\newcommand{\vecint}{\bar{\iota}}
\newcommand{\cornthree}{{\rm Corner}^\mc{G}_{k,\mfl}}
\newcommand{\cornfour}{{\rm Corner}^\mc{H}_{k,\fa}}
\newcommand{\mpgg}{\mathsf{Fav}_{\mc{G}}}

\newcommand{\lefta}{M_{1,k+1}^{[-2\eln,\gfl]}}
\newcommand{\mida}{M_{1,k+1}^{[\gfl,\gfr]}}
\newcommand{\righta}{M_{1,k+1}^{[\gfr,2\eln]}}

\newcommand{\ipdval}{d}
\newcommand{\ctemp}{d_0}

\newcommand{\wien}{W}
\newcommand{\pole}{P}
\newcommand{\pp}{p}

\newcommand{\const}{D_k}
\newcommand{\numcone}{14}
\newcommand{\numctwo}{13}
\newcommand{\numcthree}{6}
\newcommand{\rsC}{C}
\newcommand{\rsc}{c}
\newcommand{\cone}{c_1}
\newcommand{\Cone}{C_1}
\newcommand{\Ctwo}{C_2}
\newcommand{\smallc}{c_0}
\newcommand{\smallcprime}{c_1}
\newcommand{\smallcanother}{c_2}
\newcommand{\smallcnew}{c_3}
\newcommand{\Cda}{D}
\newcommand{\Kzero}{K_0}
\newcommand{\Rmac}{R}
\newcommand{\rmac}{r}
\newcommand{\conseqmac}{D}
\newcommand{\constn}{C'}
\newcommand{\coninit}{\Psi}
\newcommand{\condee}{\hat{D}}
\newcommand{\conbrac}{\hat{C}}
\newcommand{\Cnew}{\tilde{C}}
\newcommand{\Cbig}{C^*}
\newcommand{\Ctbd}{C_+}
\newcommand{\Ctbs}{C_-}

\newcommand{\imax}{i_{{\rm max}}}

\newcommand{\wlp}{{\rm WLP}}

\newcommand{\canopynumber}{\mathsf{Canopy}{\#}}

\newcommand{\cannum}{{\#}\mathsf{SC}}

\newcommand{\boundgood}{\mathsf{G}}
\newcommand{\lshift}{\mc{L}^{\rm shift}}
\newcommand{\deltapi}{\theta}
\newcommand{\rootneigh}{\mathrm{RNI}}
\newcommand{\rootneighuse}{\mathrm{RNI}}
\newcommand{\manycan}{\mathsf{ManyCanopy}}
\newcommand{\specialpt}{\mathrm{spec}}

\newcommand{\dist}{\vert\vert}
\newcommand{\fik}{\mc{F}_i^{[K,K+1]^c}}
\newcommand{\mcfa}{\mc{H}[\fa]}
\newcommand{\tent}{{\rm Tent}}
\newcommand{\goodk}{\mc{G}_{K,K+1}}
\newcommand{\pairsep}{{\rm PS}}
\newcommand{\mbf}{\mathsf{MBF}}
\newcommand{\nbd}{\mathsf{NoBigDrop}}
\newcommand{\bd}{\mathsf{BigDrop}}
\newcommand{\jleft}{j_{{\rm left}}}
\newcommand{\jright}{j_{{\rm right}}}
\newcommand{\smalljfluc}{\mathsf{SmallJFluc}}
\newcommand{\mfone}{M_{\mc{F}^1}}
\newcommand{\mfthree}{M_{\mc{G}}}
\newcommand{\rhomac}{P}
\newcommand{\phimac}{\varphi}
\newcommand{\chimac}{\chi}
\newcommand{\xnmac}{z_{\mathcal{L}}}
\newcommand{\Cwb}{E_0}
\newcommand{\initcond}{\mathcal{I}}
\newcommand{\neargeod}{\mathsf{NearGeod}}
\newcommand{\polyunique}{\mathrm{PolyUnique}}
\newcommand{\latecoal}{\mathsf{LateCoal}}
\newcommand{\nolatecoal}{\mathsf{NoLateCoal}}
\newcommand{\normalcoal}{\mathsf{NormalCoal}}
\newcommand{\regfluc}{\mathsf{RegFluc}}
\newcommand{\mdeltaweight}{\mathsf{Max}\Delta\mathsf{Wgt}}
\newcommand{\ovbar}[1]{\mkern 1.5mu\overline{\mkern-1.5mu#1\mkern-1.5mu}\mkern 1.5mu}

\newcommand{\maxmin}{\pwr}
\newcommand{\nmac}{N}

\newcommand{\high}{{\rm High}}
\newcommand{\notlow}{{\rm NotLow}}

\newcommand{\rmreg}{{\rm Reg}}

\newcommand{\down}{\mathsf{Fall}}
\newcommand{\up}{\mathsf{Rise}}

\newcommand{\paradelta}{\Delta^{\cup}\,}

\title[Modulus of continuity of polymer weight profiles]{Modulus of continuity of polymer weight profiles \\ in Brownian last passage percolation}

\author[A. Hammond]{Alan Hammond}
\address{A. Hammond\\
  Department of Mathematics and Statistics\\
 U.C. Berkeley \\
 899 Evans Hall \\
  Berkeley, CA, 94720-3840 \\
  U.S.A.}
  \email{alanmh@berkeley.edu}
  \thanks{The author is supported by NSF grant DMS-$1512908$.}

\begin{abstract} 
In last passage percolation models lying in the KPZ universality class, the energy of long energy-maximizing paths may be studied as a function of the paths' pair of endpoint locations.
Scaled coordinates may be introduced, so that these maximizing paths, or polymers, now cross unit distances with unit-order fluctuations, and have scaled energy, or weight, of unit order. In this article, we consider Brownian last passage percolation in these scaled coordinates.
In the narrow wedge case, one endpoint of such polymers is fixed, say at $(0,0) \in \R^2$, and the other is varied horizontally, over $(z,1)$, $z \in \R$, so that the polymer weight profile is a function of $z \in \R$. This profile  is known to manifest a one-half power law, having $1/2-$-H\"older continuity. The polymer weight profile may be defined beginning from a much more general initial condition.
In this article, we present a more general assertion of this one-half power law, as well as a bound on the poly-logarithmic correction. The polymer weight profile admits a modulus of continuity of order $x^{1/2} \big( \log x^{-1} \big)^{2/3}$, with 
a high degree of uniformity in the scaling parameter and over a very broad class of initial data.  

\end{abstract}

\maketitle


\tableofcontents

\section{Introduction}

The $1 + 1$ dimensional Kardar-Parisi-Zhang (KPZ) universality class includes a wide range of interface models suspended over a one-dimensional domain, in which growth in a direction normal to the surface competes with a smoothening surface tension in the presence of a local randomizing force that roughens the surface.
Such surfaces typically grow linearly, with fluctuations after that linear growth is subtracted being described by scaling exponents:
if linear growth has order $n$, then interface height above a given point has typical deviation from the mean of order $n^{1/3}$, while non-trivial correlations in this height as the spatial coordinate is varied are encountered on scale~$n^{2/3}$. Moreover, an exponent of one-half dictates the interface's regularity, with the interface height being expected to vary between a pair of locations at distance of order at most $n^{2/3}$
on the order of the square root of the distance between these locations.

Such growth models may be initiated at time zero with a given interface profile. In the narrow wedge case, when growth is initiated from a unique point, 
a limiting description of the late time interface, suitably scaled in light of the one-third and two-thirds powers and up to the subtraction of a parabola, is offered by the {\rm Airy}$_2$ process,
 which is a random  function $\mc{A}: \R \to \R$, whose finite dimensional distributions are specified by Fredholm determinants, that was introduced by~\cite{PrahoferSpohn}. 
 Another well-known initial condition is the flat case, when growth begins from a zero initial condition. Here, the ${\rm Airy}_1$ process describes the interface at late time. 
  The one-half power law for interface regularity is expressed by the H{\"o}lder-$1/2-$-continuity of the processes~${\rm Airy}_1$ and~${\rm Airy}_2$, which was proved in~\cite{QR13}.

Growth may be initiated from a much more general initial condition than in these narrow wedge or flat cases.
For initial conditions that grow at most linearly, it has been anticipated that a limiting description of the suitably scaled late-time interface should exist in these cases also. Indeed, in a recent preprint~\cite{MQR17}, Matetski, Quastel and Remenik
have utilized a biorthogonal ensemble representation found by~\cite{Sas05,BFPS07} associated to the totally asymmetric exclusion process in order to  find Fredholm determinant formulas for the multi-point distribution of the height function of this growth process begun from an arbitrary initial condition. Using these formulas to take the KPZ scaling limit, the authors construct a scale invariant Markov process that lies at the heart of the KPZ universality class. The time-one evolution of this Markov process may be applied to very general initial data, and the result is the scaled profile begun from such data, which generalizes the ${\rm Airy}_1$ and~${\rm Airy}_2$ processes seen in the flat and narrow wedge cases. 
These more general limiting processes also enjoy  H{\"o}lder-$1/2-$-continuity:  see~\cite[Theorem 4.4]{MQR17}.

The broad range of interface models that are rigorously known or expected to lie in the KPZ universality class includes many last passage percolation models.
Such an LPP model comes equipped with a planar random environment, which is independent in disjoint regions. 
Directed paths, that are permitted say to move only in a direction in the first quadrant, are then assigned energy via this randomness, by say integrating the environment's value along the path.
For a given pair of planar points, the path attaining the maximum energy over directed paths with such endpoints is called a geodesic. 
The random interface model that we alluded to at the outset is then specified as the maximum geodesic energy when one geodesic endpoint is varied and the other held fixed, in the narrow wedge case, or when the other is free to vary and is rewarded according to the initial condition, in the more general case.   
The one-third and two-thirds power laws for typical deviation of maximum energy and for lateral correlation have been rigorously demonstrated for only a few LPP models, each of which enjoys an integrable structure: the seminal work of Baik, Deift and Johansson~\cite{BDJ1999} 
 rigorously established the one-third exponent, and moreover obtained the GUE Tracy-Widom distributional limit, for the case of Poissonian last passage percolation, while the two-thirds power law for transversal fluctuation was derived for this model by Johansson~\cite{Johansson2000}. 
For models in which these two exponents have been rigorously identified,
the exponent pair dictates a system of scaled coordinates in which the concerned maximizing paths and their weights are unit-order, random, quantities: the scaled geodesics may be called polymers, and their scaled energies, weights.

 Brownian last passage percolation  is an LPP model with attractive integrable and probabilistic features. 
In this article, we study the scaled interface profile (that is, the polymer weight profile) in Brownian LPP begun from a very general initial condition. We present results proving a more precise version of the one-half power law for interface regularity than has been established hitherto. Here are two of the  main conclusions:
\begin{itemize}
\item In Theorem~\ref{t.differenceweight}, we prove that the maximum difference in the weight of two point-to-point polymers whose endpoints differ by at most a small scaled quantity~$\e$ exceeds $\e^{1/2}R$
with probability at most $\exp \big\{  - O(1) R^{3/2}  \big\}$ for a very broad range of values of $R$, uniformly in the scaling parameter for Brownian LPP.
\item In Theorem~\ref{t.wlp.one}, we prove that any weak limit point of the scaled interface profiles, as the scaling parameter tends to infinity, has sample paths that admit a modulus of continuity of the order of $x^{1/2} \big( \log x^{-1} \big)^{2/3}$. This assertion, alongside a finite-$n$ counterpart, Theorem~\ref{t.nmodcon}, is proved uniformly over a large class of the data that initiates the random growth. 
\end{itemize}

For a given choice of initial condition, the weak limit point in Theorem~\ref{t.wlp.one} may be expected to be unique and to coincide with the interface profile obtained from this initial data by evolving for a given duration the Markov operator constructed in~\cite{MQR17}. However, this Markov operator has been constructed as a limit of totally asymmetric exclusion, so at present this assertion is not proved.
Were the techniques of~\cite{MQR17} to be adapted to hold for Brownian last passage percolation, it would then presumably be possible to assert the upper bound of order $x^{1/2} \big( \log x^{-1} \big)^{2/3}$ on  modulus of continuity for general initial condition interface profiles under the KPZ fixed point.

The strongly on-scale assertion of the one-half power law for profile regularity in Theorem~\ref{t.differenceweight} plays a significant role in two companion papers. In~\cite{NonIntPoly}, it is harnessed to prove that, in Brownian last passage percolation, it is a superpolynomial rarity that a large number of disjoint polymers coexist in a unit-order scaled region. In~\cite{Patch}, this assertion is exploited
to make a strong unit-order Brownian comparison for polymer weight profiles (about which more momentarily).

Beyond the conclusions just discussed, the present article also presents a useful tool, Proposition~\ref{p.maxminweight}. Although the weight of a polymer is random, this weight is dictated in the large by parabolic curvature, with the randomness playing a unit-order role once this curvature is accounted for. The proposition shows that the discrepancy between polymer weight and parabola is controlled uniformly 
as the polymer's endpoints are varied over compact intervals lying in a very broad region. This tool is needed in the present article and in~\cite{NonIntPoly}. For exponential or Poissonian last passage percolation, a similar tool has been developed, in~\cite[Propositions~$10.1$ and~$10.5$]{SlowBondSol}.

We also mention that an alternative expression of the one-half power law for interface regularity is the assertion that Airy processes such as   {\rm Airy}$_1$ and {\rm Airy}$_2$, or scaled interface models in the last passage percolation setting, locally resemble Brownian motion.
Such statements may be understood in a local limit, when Gaussianity of a process $\mc{A}$ is proved for the low~$\e$ limit for the random variable $\e^{-1/2} \big( \mc{A}(x+\e) - \mc{A}(x) \big)$
associated to any given $x \in \R$. 
Finite dimensional distributional convergence to Brownian motion (of diffusion rate two) in this limit  has been proved for the Airy$_2$ process in~\cite{Hagg}, for the Airy$_1$ process in~\cite{QR13},
and for the more general versions of these Airy processes constructed in \cite{MQR17} in Theorem~$4.4$ of that paper; in~\cite{Pimentel17}, similar 
 local limit results for general initial condition profiles have been obtained for geometric last passage percolation models. 
Comparison to Brownian motion may also be made without taking such a local limit. In~\cite{AiryLE}, 
the ${\rm Airy}_2$ process was understood to be absolutely continuous with respect to Brownian motion on a unit-order interval, by a technique in which this process is embedded as the uppermost curve in a random ensemble of, in effect, mutually avoiding Brownian motions. (This {\em Brownian Gibbs} technique will play a fundamental role in the present article, and we will return to it.) This inference was improved in~\cite{BrownianReg}, where the implied Radon-Nikodym derivative is shown to lie in all $L^p$-spaces for $p \in (1,\infty)$, albeit after an affine shift is applied to the ${\rm Airy}_2$ process, so that comparison is made not to Brownian motion but to Brownian bridge. In a companion paper to the present article~\cite{Patch}, 
the problem of unit-order scale Brownian comparison is made for the class of Brownian LPP polymer weight profiles, begun from general initial data,  that are the subject of the present article.
It is in essence shown there that a given unit-order interval may be split into a random but controlled number of intervals in such a way that the profile when restricted to the smaller intervals has, after affine adjustment, a Radon-Nikodym derivative with respect to Brownian bridge that lies in $L^p$ for $p \in (1,3)$.

\subsection{Brownian last passage percolation [LPP]}\label{s.brlpp}
We now define this model.
On a probability space carrying a law labelled~$\PP$,
let $B:\Z \times \R \to \R$ denote an ensemble of independent  two-sided standard Brownian motions $B(k,\cdot):\R\to \R$, $k \in \Z$.

Let $i,j \in \Z$ with $i \leq j$.
We denote the integer interval $\{i,\cdots,j\}$ by $\llbracket i,j \rrbracket$.
Further let $x,y \in \R$ with $x \leq y$.
Consider the collection of  non-decreasing lists 
 $\big\{ z_k: k \in \llbracket i+1,j \rrbracket \big\}$ of values $z_k \in [x,y]$. 
With the convention that $z_i = x$ and $z_{j+1} = y$,
we associate an energy $\sum_{k=i}^j \big( B ( k,z_{k+1} ) - B( k,z_k ) \big)$ to any such list.
We then define  the maximum energy
$$
M^1_{(x,i) \to (y,j)} \, = \, \sup \, \bigg\{ \, \sum_{k=i}^j \Big( B ( k,z_{k+1} ) - B( k,z_k ) \Big) \, \bigg\} \, , 
$$
where the supremum is taken over all such lists. The random process $M^1_{(0,1) \to (\cdot,n)}: [0,\infty) \to \R$ was introduced by~\cite{GlynnWhitt} and further studied in~\cite{O'ConnellYor}.


The one-third and two-thirds KPZ scaling considerations that we outlined earlier in the introduction are manifest in Brownian LPP. When the ending height $j$ exceeds the starting height $i$ by a large quantity $n \in \N$, and the location $y$ exceeds $x$ also by $n$, then the maximum energy grows linearly, at rate $2n$,
and has a fluctuation about this mean of order $n^{1/3}$. Moreover, if $y$ is permitted to vary from this location, then it is changes of $n^{2/3}$ in its value that result in a non-trivial correlation of the maximum energy with its original value.

These facts prompt us to introduce scaled coordinates to describe the two endpoint locations, and a notion of scaled maximum energy, which we will refer to as weight. 
Let  $n \in \N$, and 
suppose that $x,y \in \R$ satisfy 
 $y \geq x - 2^{-1} n^{1/3}$.
Define 
\begin{equation}\label{e.weightmzeroone} 
  \weight_{n;(x,0)}^{(y,1)} \,     =  \,   2^{-1/2} n^{-1/3} \Big(  M^1_{(2n^{2/3}x,0) \to (n  + 2n^{2/3}y,n)} - 2n  -  2n^{2/3}(y-x) \Big) \, .
\end{equation}
(Clearly, $n$ must be positive. In fact, $\N$ will denote $\{1,2,\cdots \}$ throughout.)

Consistently with the facts just mentioned, the quantity  $\weight_{n;(x,0)}^{(y,1)}$ may be expected to be, for given real choices of $x$ and $y$, a unit-order random quantity, whose law is tight in the scaling parameter $n \in \N$. The quantity describes, in units chosen to achieve this tightness, the maximum possible energy associated to journeys which in the original coordinates occur between  $(2n^{2/3}x,0)$ and $(n  + 2n^{2/3}y,n)$.
In scaled coordinates, this is a journey between $(x,0)$ and $(y,1)$.
We view the first coordinate as space and the second as time, so this journey is between $x$ and $y$ over the unit time interval $[0,1]$.

 Underlying this definition is a geometric picture of scaled maximizing paths, or polymers, that achieve these weight values. We will defer explicitly defining these polymers, but it is useful to bear in mind that 
  $\weight_{n;(x,0)}^{(y,1)}$ equals the weight of a polymer that travels between locations that in scaled coordinates equal $(x,0)$ and $(y,1)$.

\subsection{Main results}

In four subsections, we present the principal conclusions: on  polymer weight difference under horizontal perturbation of endpoints; our finite-$n$ 
assertion concerning the modulus of continuity of polymer weight profiles from general initial condition; the inference made about weak limit points of such profiles in the high~$n$ limit; 
and a general tool, on the rarity of deviation from parabolic curvature by polymer weights.


  

\subsubsection{Polymer weight change under horizontal perturbation of endpoints}\label{s.pwc}
Set $Q:\R \to \R$, $Q(z) = 2^{-1/2} z^2$.
The polymer weight   $\weight_{n;(x,0)}^{(y,1)}$  has a globally parabolic profile, hewing to the shape $-Q(y-x)$. When this parabolic term is added to the polymer weight, the result is a random process in $(x,y)$
which typically suffers changes of order $\e^{1/2}$
when $x$ or $y$ are varied on a small scale $\e > 0$.
Our first main result gives rigorous expression to this statement, uniformly in $(n,x,y) \in \N \times \R \times \R$ for which the difference $\vert y - x \vert$ is permitted to inhabit an expanding region about the origin, of scale~$n^{1/18}$.

\begin{theorem}\label{t.differenceweight}
Let $\e \in (0,2^{-4}]$. 
Let $n \in \N$ satisfy
$n \geq 10^{32} c^{-18}$ and let $x,y \in \R$ satisfy  $\big\vert x - y  \big\vert \leq 2^{-2} 3^{-1} \rsc  n^{1/18}$.
Let 
 $R \in \big[10^4 \, , \,   10^3 n^{1/18} \big]$.
Then
\begin{equation}\label{e.differenceweight}
\PP \left( \sup_{\begin{subarray}{c} u_1,u_2 \in [x,x+\e] \\
    v_1,v_2 \in [y,y+\e]  \end{subarray}} \Big\vert \weight_{n;(u_2,0)}^{(v_2,1)} + Q(v_2 - u_2) - \weight_{n;(u_1,0)}^{(v_1,1)} - Q(v_1 - u_1) \Big\vert  \, \geq \, \e^{1/2}
  R  \right)
\end{equation}
  is at most  $10032 \, C  \exp \big\{ - c_1 2^{-21}   R^{3/2}   \big\}$.
\end{theorem}
Here, we set $c_1 = 2^{-5/2} c \wedge 1/8$, where $\wedge$ denotes minimum. Bounds in Theorem~\ref{t.differenceweight}, and many later results, have been expressed explicitly up to two positive constants $c$ and $C$. See Subsection~\ref{s.regular} for a discussion of the role of this pair of constants.

The imposition in Theorem~\ref{t.differenceweight} that 
 $R \in \big[10^4 \, , \,   10^3 n^{1/18} \big]$ is rather weak, with the case where $R$ is fixed being of interest; and indeed, the decay rate asserted by the theorem 
 is already very fast when 
$R$ is of order $n^{1/18}$.

\subsubsection{Maximum local variation of polymer weight profiles from general initial data}

What do we mean by such polymer weight profiles?  The random function $y \to   \weight_{n;(0,0)}^{(y,1)}$
may be viewed as the weight profile obtained by scaled maximizing paths that travel from the origin at time zero to the variable location $y$ at time one. This insistence that the paths must begin at the origin, called the narrow wedge by physicists, is of course rather special. 
We now make a more general definition, of the  $f$-rewarded line-to-point polymer weight  $\weight_{n;(*:f,0)}^{(y,1)}$. Here, $f$ is an initial condition, defined on the real line. 
Paths may begin anywhere on the real line at time zero; they travel to $y \in \R$ at time one. (Because they are free at the beginning and fixed at the end, we refer to these paths as `line-to-point'.) They begin with a reward given by evaluating $f$ at the starting location, and then gain the weight associated to the journey they make.  The value $\weight_{n;(*:f,0)}^{(y,1)}$, which we will define momentarily, denotes the maximum $f$-rewarded weight of all such paths.   In the notation $\weight_{n;(*:f,0)}^{(y,1)}$, we again use subscript and superscript expressions to refer to space-time pairs of starting and ending locations. The starting spatial location is being denoted $*:f$. The star is intended to refer to the free time-zero endpoint, which may be varied, and the $:f$ to the reward offered according to where this endpoint is placed.

The next definition specifies essentially the broadest class of $f$ suitable for a study of the weight profiles $y \to \weight_{n;(*:f,0)}^{(y,1)}$ for all sufficiently high $n \in \N$.

\begin{definition}\label{d.if}
Writing  $\ovbar\coninit = \big( \coninit_1,\coninit_2,\coninit_3 \big) \in (0,\infty)^3$ for a triple of positive reals, we let $\initcond_{\ovbar\coninit}$
denote the set of measurable functions $f:\R \to \R \cup \{ - \infty \}$ such that
$f(x) \leq \coninit_1 \big( 1 + \vert x \vert \big)$
and $\sup_{x \in [-\coninit_2,\coninit_2]} f(x) > - \coninit_3$.
\end{definition}

For $f$ lying in one of the function spaces 
$\initcond_{\ovbar\coninit}$, we now formally define the $f$-rewarded line-to-point polymer weight  $\weight_{n;(*:f,0)}^{(y,1)}$ to be 
$$
    \sup_{x \in (-\infty,2^{-1}n^{1/3} + y]} \big(  \weight_{n;(x,0)}^{(y,1)}    + f(x) \big)  \, .
$$

Our second main result asserts that the maximal variation in $f$-rewarded $n$-polymer weight over length~$\e > 0$ intervals in $[-1,1]$ is a controlled random multiple of $\e^{1/2} \big( \log \e^{-1}\big)^{2/3}$. 
The bound on probability, above scale 
$e^{-O(1)n^{1/12}}$, is asserted uniformly in initial 
data, and in $(n,\e)$, except for very small $\e \leq e^{-O(1)n^{1/12}}$. 
\begin{theorem}\label{t.nmodcon}
For $\ovbar\coninit \in (0,\infty)^3$, some  $c',r_0 = c'(\ovbar\coninit),r_0(\ovbar\coninit) > 0$ and all 
$f \in \initcond_{\ovbar\coninit}$,~$n~\in~\N$~and~$r~\geq~r_0$,
\begin{equation}\label{e.nmodcon}
\PP \left( \sup_{\begin{subarray}{c} y,z \in [-1,1], \\
    2\exp \{-c' n^{1/12} \} < z - y < e^{-1}  \end{subarray}} \frac{\Big\vert \, \weight_{n;(*:f,0)}^{(z,1)} - \weight_{n;(*:f,0)}^{(y,1)} \, \Big\vert}{(z-y)^{1/2} \big( \log (z-y)^{-1} \big)^{2/3}} \, \geq \, r  \right) \leq  2^{47} c^{-4/3} r^{-2} \big( \log r \big)^{4/3} \vee 4e^{-c' n^{1/12}} \, .
\end{equation}
\end{theorem}

\subsubsection{Modulus of continuity of weak limits of weight profiles from general initial data}

Let $n \in \N$,  $\ovbar\coninit \in (0,\infty)^3$ and  $f \in \initcond_{\ovbar\coninit}$.
Let $\nu_{n;(*:f,0)}^{([-1,1],1)}$ denote the law of the random function
$$
[-1,1] \to \R : y \to \weight_{n;(*:f,0)}^{(y,1)} \, . 
$$

The control offered by Theorem~\ref{t.nmodcon} is certainly sufficient to show that the curves of any weak limit point of $\nu_{n;(*:f,0)}^{([-1,1],1)}$ as $n \to \infty$ admit modulus of continuity $z^{1/2} \big( \log z^{-1} \big)^{2/3}$, up to a random factor that is controlled uniformly in the choice of limit point.

To formulate a theorem in this regard,
let $\mathcal{A}$ be an arbitrary index set, and let  $\big\{ \nu_{n,\alpha}: n \in \N \big\}$, $\alpha \in \mathcal{A}$, be an $\mathcal{A}$-indexed collection of sequences of probability measures on the Borel $\sigma$-algebra of a  given Hausdorff topological space. The collection is here called $\mc{A}$-uniformly tight  if, for each $\e > 0$, there exist $n_0 \in \N$ and a compact set $K$ such that 
 $\nu_{n,\alpha}(K) \geq 1 - \e$ whenever $n \geq n_0$ and $\alpha \in \mathcal{A}$.

\begin{theorem}\label{t.wlp.one}
Let $\ovbar\coninit \in (0,\infty)^3$ denote a triple of positive reals. 
\begin{enumerate}
\item  Suppose that $n \in \N$ satisfies
$n > 2^{-3/2} \coninit_1^3 \vee 8 (\coninit_2  + 1)^3$.
 Let    $f \in \initcond_{\ovbar\coninit}$.  Then  the  measure  $\nu_{n;(*:f,0)}^{([-1,1],1)}$ is supported on
the space $\mc{C}$ of continuous real-valued functions on $[-1,1]$.
\item The collection of sequences of probability measures
$\big\{ \nu_{n;(*:f,0)}^{([-1,1],1)}: n \in \N  \big\}$ indexed by $f \in \initcond_{\ovbar\coninit}$ is $\initcond_{\ovbar\coninit}$-uniformly tight. Here, the space $\mc{C}$
is endowed with the topology of uniform convergence.
\item A law on $\mc{C}$ is said to belong to the {\em weak limit point} set 
$\wlp_{\ovbar{\coninit}}$ if, for some sequence $f_n \in  \initcond_{\ovbar\coninit}$, $n \in \N$, it equals the weak limit of the laws    $\nu_{n;(*:f_n,0)}^{([-1,1],1)}$ along some subsequence of $n \in \N$. By~$(2)$, $\wlp_{\ovbar{\coninit}} \not= \emptyset$. 
For any $\nu \in \wlp_{\ovbar{\coninit}}$, let $X$ be $\nu$-distributed. Then, for~$r~\geq~r_0$,
\begin{equation}\label{e.nur}
\nu \left( \, \sup_{\begin{subarray}{c} x,y \in [-1,1], \\
    x < y < x + e^{-1}  \end{subarray}} \frac{\big\vert X(y) - X(x) \big\vert}{(y-x)^{1/2} \big( \log (y-x)^{-1} \big)^{2/3}} \, \geq \, r \,  \right) \, \leq \,   2^{47} c^{-4/3} r^{-2} \big( \log r \big)^{4/3} \, ,
\end{equation}
where Theorem~\ref{t.nmodcon} provides the constant $r_0 = r_0(\ovbar\coninit)$.
\end{enumerate}
\end{theorem}

Brownian motion on a unit interval has modulus of continuity of order $x^{1/2} \big( \log x^{-1} \big)^{1/2}$,
and it may be expected that some version of Theorem~\ref{t.wlp.one}(3) is valid with the logarithmic power of two-thirds replaced by one-half. Indeed, in the special case of narrow wedge initial data, such a result has been proved: see~\cite[Theorem~$2.13$]{BrownianReg}, or \cite[Theorem~$1.11(1)$]{BrownianReg}
for a result concerning the {\rm Airy}$_2$ process.

\subsubsection{Tail behaviour of polymer weight suprema and infima}

Theorem~\ref{t.differenceweight} quantifies polymer weight changes in response to horizontal endpoint perturbation after the weight has been adjusted by the addition of the parabola~$Q(z) = 2^{-1/2} z^2$. 
In our fourth result, we verify that the point-to-point polymer weight profile indeed strongly hews to this given parabola. 
The regime where this is verified is that in which the polymer endpoints differ by at most an order of $n^{1/18}$.
Within this zone, the inference is made uniformly as the endpoints vary over any given unit-order region.
\begin{proposition}\label{p.maxminweight}
Let $n \in \N$ satisfy 
$n \geq 10^{29} \vee 2(c/3)^{-18}$. Let $x,y \in \R$ 
satisfy
$\big\vert x - y  \big\vert \leq   6^{-1}  \rsc  n^{1/18}$.  Let $t \in \big[  34 \, , \, 4 n^{1/18} \big]$. 
Then
\begin{equation}\label{e.weightuvsupremum} 
 \PP \bigg( \sup_{u,v \in [0,1]} \Big( \weight_{n;(x+u,0)}^{(y+v,1)} + Q( y+v - x - u ) \Big) \geq t \bigg)  \leq    139 C  \exp \big\{ - c_1 2^{-10} t^{3/2} \big\}  
\end{equation}
and
\begin{equation}\label{e.weightuvinfimum}
  \PP \bigg(  \inf_{u,v \in [0,1]} \Big( \weight_{n;(x+u,0)}^{(y+v,1)} + Q(y+v-x - u) \Big) \leq - t  \bigg) 
   \leq  261 C \exp \big\{ - c_1 2^{-3} t^{3/2} \big\} \, . 
\end{equation}
\end{proposition}
In~\cite[Propositions~$10.1$ and~$10.5$]{SlowBondSol}, comparable bounds are proved for exponential and Poissonian LPP, with bounds  of the form
$\exp \big\{ - O(1) t \big\}$. These propositions have the flexibility of treating extremal weights of polymers whose endpoints are permitted to vary over compact regions in space as well as time, rather than merely time, as it is the case for Proposition~\ref{p.maxminweight}.

\subsection{The road map}

Theorem~\ref{t.differenceweight}, which is a key result underlying Theorems~\ref{t.nmodcon} and~\ref{t.wlp.one}, is proved using ideas similar to the proof of the Kolmogorov continuity criterion. The authors of~\cite{QR13} note in Section~$1.1$ that the task of checking the Kolmogorov criterion on the basis of suitable two-point information for such processes as {\rm Airy}$_1$
has turned out to be surprisingly difficult. 
Similar subtleties arise in our context:
two-point information has to be presented in a way that is valid on arbitrarily small scales, without the index $n$ needing to rise. The crucial tool that will enable the derivation of Theorem~\ref{t.differenceweight}
is a powerful two-point estimate with the necessary attributes, Proposition~\ref{p.locreg}.

Section~\ref{s.basicnot} introduces notation for the use of scaled coordinates and presents some basic results about polymer weight.

 In Section~\ref{s.lineensembles},  the engine for our main results,  Proposition~\ref{p.locreg}, namely the two-point estimate for narrow wedge weight profiles, is stated and proved. In this section, we will explain how the narrow wedge profile may be embedded as the uppermost curve in a certain system of ordered random continuous curves called a line ensemble. A suitably normalized version of any such line ensemble has the {\em Brownian Gibbs} property, which in essence means it is a system of mutually avoiding Brownian bridges. Its curves moreover have a globally parabolic shape, and a definition of {\em regular} ensemble is made to capture these attributes. The short proof of Proposition~\ref{p.locreg} harnesses the Brownian Gibbs property in an essential way. Certain further properties of regular  ensembles are needed, and these also appear in Section~\ref{s.lineensembles} as Proposition~\ref{p.mega}, quoted from~\cite{BrownianReg}.

In three further sections  are then respectively proved 
Proposition~\ref{p.maxminweight}; Theorem~\ref{t.differenceweight}; and Theorems~\ref{t.nmodcon} and~\ref{t.wlp.one}.

\subsubsection{Comment on the companion papers~\cite{BrownianReg},~\cite{NonIntPoly} and~\cite{Patch}.} 
Via the upcoming Proposition~\ref{p.mega} and Lemma~\ref{l.pardom}, this article draws on Brownian Gibbs results developed in~\cite{BrownianReg}. That article is long and it is worth pointing out that the concerned results in~\cite{BrownianReg} are simple and have short proofs.
The present article's main conclusions about scaled Brownian LPP are applied in the later two companion papers. The article  may be read on its own,
or viewed as part of this four-paper study, an overview of which appears in~\cite[Section~$1.2$]{BrownianReg}.

\subsubsection{Acknowledgments.}
The author thanks Riddhipratim Basu, Ivan Corwin,  Shirshendu Ganguly and Jeremy Quastel for valuable conversations, and three referees for helpful comments.

\section{The basics: notation, scaling, polymers and their weight}\label{s.basicnot}
In consecutive subsections, we introduce notation; describe Brownian LPP in scaled coordinates; offer a principle that aids in working with these coordinates, and an application; discuss basics about polymers; and provide a simple result about them.
\subsection{General notation and structure}
\subsubsection{Notation}
Let $i,j \in \Z$ with $i \leq j$. Recall that $\llbracket i,j \rrbracket$ denotes the integer interval $\{ i,\cdots  , j \}$.

For $k \geq 1$, we write $\R^k_\leq$
for the subset of $\R^k$ whose elements $(z_1,\cdots,z_k)$
are non-decreasing sequences. When the sequences are increasing, we instead write $\R^k_<$. We also use the notation $A^k_\leq$ and $A^k_<$.
Here, $A \subset \R$ and the sequence elements are supposed to belong to $A$.
We will typically use this notation when $k=2$.

A bar over a symbol indicates a vector, as in the usage
$\ovbar\coninit = \big( \coninit_1,\coninit_2,\coninit_3 \big)$
in Theorems~\ref{t.nmodcon} and~\ref{t.wlp.one}.

\subsubsection{The role of hypotheses invoked during proofs}\label{s.calcderexplain}

Our proofs invoke several inputs, notably  three $\rmreg$ conditions that specify the notion of a regular ensemble, and Propositions~\ref{p.mega} and~\ref{p.locreg}. 
Whenever such results are invoked, certain conditions on the concerned hypotheses will be needed. We will always note explicitly what these conditions are, whenever such an application is made.
Clearly, it is necessary that the hypotheses of the result that is being proved collectively imply all the conditions that are invoked during its proof. 
The work needed to do this for a given result may be called the {\em calculational derivation} of that result. 
These derivations have almost no conceptual content, reach conclusions that in their overall form are plausible, consist of largely trivial steps, and will be of interest to only the most committed of readers (perhaps only those who are actually applying the results). 
In some cases, however, the derivations occupy a fair amount of space. We have chosen to separate the principal calculational derivations from the body of the proofs in this article.
The concerned results are  Theorem~\ref{t.differenceweight}, 
Proposition~\ref{p.maxminweight}, Proposition~\ref{p.dyadic}, and Lemmas~\ref{l.regfluc} and~\ref{l.equicty}.  
Their calculational derivations are presented in 
Appendix~$A$, which is a supplement
that appears in the version~\url{math.berkeley.edu/~alanmh/papers/ModCon.pdf} of this paper on the author's webpage. The latex source code for this version is an ancillary file to the present arXiv submission.

\subsection{Scaling: staircases to zigzags,  energy to weight, and geodesics to polymers}

\subsubsection{Staircases.}

Taking  $i,j \in \N$ with $i \leq j$, and 
$x,y \in \R^2_\leq$,
we have ascribed in Section~\ref{s.brlpp} an energy to  any non-decreasing list $\big\{ z_k: k \in \llbracket i+1,j \rrbracket \big\}$ of values $z_k \in [x,y]$.
In order to emphasise the geometric aspects of this definition, and in the hope that it may aid the visualization of the concerned concepts,
 we associate to each list a  subset of $[x,y] \times [i,j] \subset \R^2$, which will be the range of a piecewise affine path,  
 that we call a staircase.
 
 To define the staircase associated to $\big\{ z_k: k \in \llbracket i+1,j \rrbracket \big\}$, 
we again adopt the convention that $z_i = x$ and  $z_{j+1} = y$. 
The staircase is specified as the union of certain horizontal planar line segments, and certain vertical ones.
The horizontal segments take the form $[ z_k,z_{k+1} ] \times \{ k \}$ for $k \in \llbracket i , j \rrbracket$.
The right and left endpoints of each consecutive pair of horizontal segments are interpolated by a vertical planar line segment of unit length. It is this collection of vertical line segments that form
the vertical segments of the staircase.

The resulting staircase may be depicted as the range of an alternately rightward and upward moving path from starting point $(x,i)$ to ending point $(y,j)$. 
The set of staircases with these starting and ending points will be denoted by $\staircase_{(x,i) \to (y,j)}$.
Such staircases are in bijection with the collection of non-decreasing lists already considered. Thus, any staircase $\phi \in \staircase_{(x,i) \to (y,j)}$
is assigned an energy $E(\phi) = \sum_{k=i}^j \big( B ( k,z_{k+1} ) - B( k,z_k ) \big)$ via the associated $z$-list.

\subsubsection{Energy maximizing staircases are called geodesics.}
A staircase  $\phi \in \staircase_{(x,i) \to (y,j)}$ whose energy  attains the maximum value $M^1_{(x,i) \to (y,j)}$ is called a geodesic from $(x,i)$ to~$(y,j)$.
It is a simple consequence of the continuity of the constituent Brownian paths $B(k,\cdot)$
that this geodesic exists for all choices of $(x,y) \in \R^2_\leq$.
It is also true, and is proved in~\cite[Lemma~$A.1$]{Patch},
 that, for any given such choice of the pair $(x,y)$, there is an almost surely unique geodesic from  $(x,i)$ to~$(y,j)$.
However, this uniqueness will not be needed in the present article.

\subsubsection{The scaling map.}
For $n \in \N$, consider the $n$-indexed {\em scaling} map $R_n:\R^2 \to \R^2$ given by
$$
 R_n \big(v_1,v_2 \big) = \Big( 2^{-1} n^{-2/3}( v_1 - v_2) \, , \,   v_2/n \Big) \, .
$$ 
 
The scaling map acts on subsets $C$ of $\R^2$ by setting
$R_n(C) = \big\{ R_n(x): x \in C \big\}$.

\subsubsection{Scaling transforms staircases to zigzags.}
The image of any staircase under $R_n$
will be called an $n$-zigzag. The starting and ending points of an $n$-zigzag $Z$ are defined to be the image under $R_n$
of such points for the staircase $S$ for which $Z = R_n(S)$.
 
Note that the set of horizontal lines is invariant under $R_n$, while vertical lines are mapped to lines of gradient  $- 2 n^{-1/3}$.
As such, an $n$-zigzag is the range of a piecewise affine path from the starting point to the ending point which alternately moves rightwards  along horizontal line segments  and northwesterly along sloping line segments, where each sloping line segment has gradient  $- 2 n^{-1/3}$.

 \subsubsection{Compatible triples.}
 Let $(n,t_1,t_2) \in \N \times \R^2_<$, which is to say that $n \in \N$ and $t_1,t_2 \in \R$ with $t_1 < t_2$.
 Taking $x,y \in \R$, does there exist an $n$-zigzag from $(x,t_1)$ and $(y,t_2)$?
 As far as the data $(n,t_1,t_2)$ is concerned, such an $n$-zigzag may exist only if 
 \begin{equation}\label{e.ctprop}
     \textrm{$t_1$ and $t_2$ are integer multiplies of $n^{-1}$} \, .
\end{equation}
We say that data $(n,t_1,t_2)  \in \N \times \R^2_<$ is a {\em compatible triple} if it verifies the last condition. 
We will consistently impose this condition, whenever we seek to study $n$-zigzags whose lifetime is $[t_1,t_2]$.
The use of compatible triples should be considered to be a fairly minor, microscopic, detail. As the index $n$ increases, the $n^{-1}$-mesh becomes finer, so that the space of $n$-zigzags better approximates a field of functions, defined on arbitrary finite intervals of the vertical coordinate, and taking values in the horizontal coordinate.     

An important piece of notation associated to a compatible triple is $\tot$, which will denote the difference $t_2 - t_1$. The law of the underlying Brownian ensemble $B: \Z \times \R \to \R$ is invariant under integer shifts in the first, curve indexing, coordinate. This translates to an invariance in law of scaled objects under vertical shifts by multiples of $n^{-1}$, something that makes the parameter $\tot$
of far greater relevance than $t_1$ or $t_2$.

\begin{figure}[ht]
\begin{center}
\includegraphics[height=7cm]{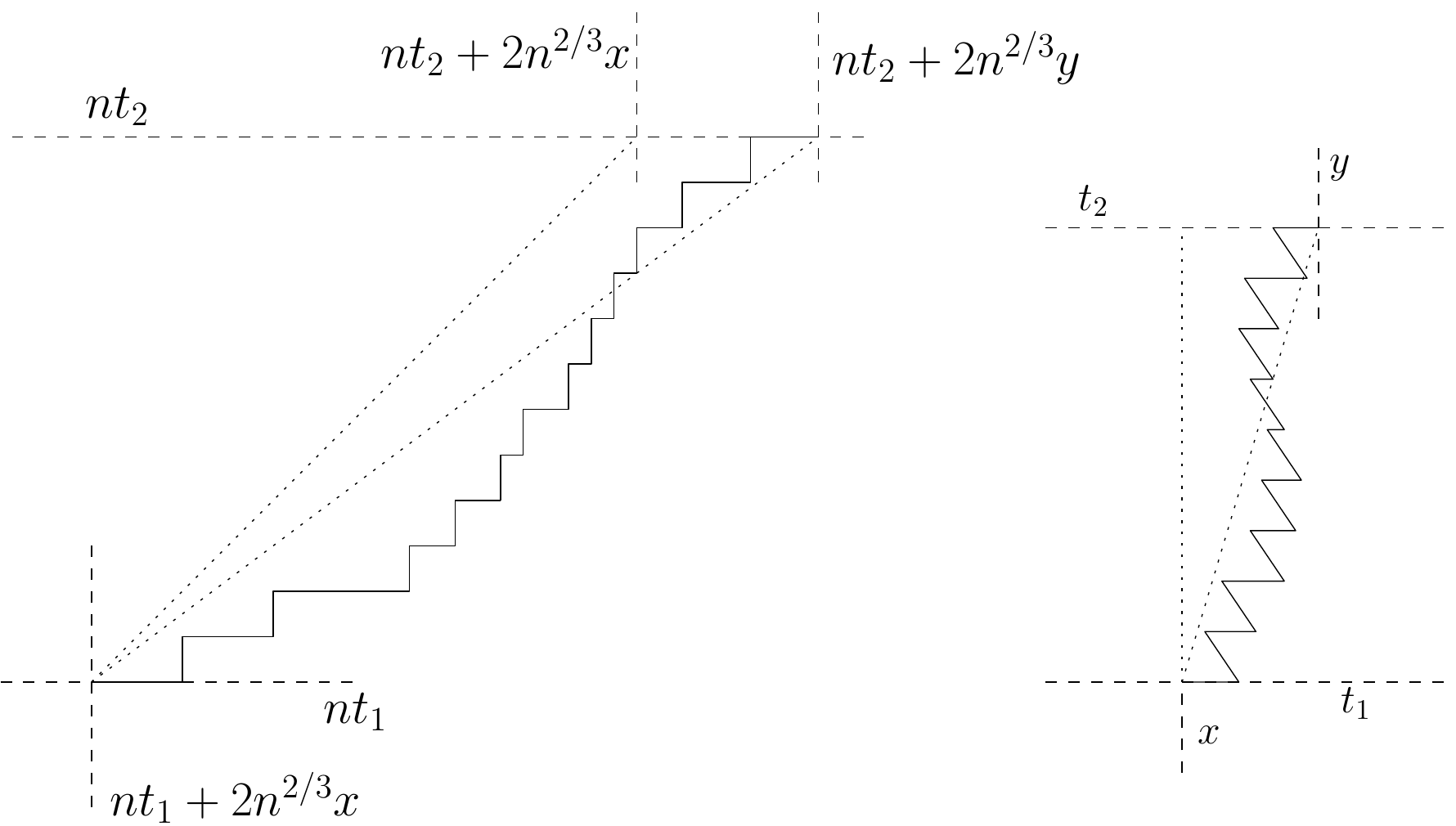}
\caption{Let $(n,t_1,t_2)$ be a compatible triple and let $x, y \in \R$. The endpoints of the geodesic in the left sketch are such that, when the scaling map~$R_n$ is applied to produce the right sketch,  the result is an $n$-polymer from $(x,t_1)$ to $(y,t_2)$.
  }
\label{f.scaling}
\end{center}
\end{figure}

\subsubsection{Staircase energy scales to zigzag weight.}
Let $n \in \N$ and $(i,j) \in \N^2_<$.
Any $n$-zigzag $Z$ from $(x,i/n)$ to $(y,j/n)$  is ascribed a scaled energy, which we will refer to as its weight, 
$\weight(Z) = \weight_n(Z)$, given by 
\begin{equation}\label{e.weightzigzag}
 \weight(Z) =  2^{-1/2} n^{-1/3} \Big( E(S) - 2(j - i)  - 2n^{2/3}(y-x) \Big) 
\end{equation}
where $Z$ is the image under $R_n$ of the staircase $S$.

\subsubsection{Maximum weight.} Let $n \in \N$. The quantity
 $\weight_{n;(x,0)}^{(y,1)}$ specified in~(\ref{e.weightmzeroone}) is nothing other than 
 the maximum weight ascribed to any $n$-zigzag from $(x,0)$ to $(y,1)$.
 
 Let $(n,t_1,t_2) \in \N \times \R^2_<$ be a compatible triple. 
Suppose that $x,y \in \R$ satisfy 
 $y \geq x - 2^{-1} n^{1/3} \tot$. We now offer a definition of  $\weight_{n;(x,t_1)}^{(y,t_2)}$ such that this quantity equals maximum weight of any $n$-zigzag from $(x,t_1)$ to $(y,t_2)$.
 We must set 
\begin{equation}\label{e.weightm}
  \weight_{n;(x,t_1)}^{(y,t_2)} \,     =  \,   2^{-1/2} n^{-1/3} \Big(  M^1_{(n t_1 + 2n^{2/3}x,n t_1) \to (n t_2 + 2n^{2/3}y,n t_2)} - 2n \tot -  2n^{2/3}(y-x) \Big) \, .
\end{equation}
 \subsubsection{Highest weight zigzags are called polymers.}\label{s.polymer}
 An $n$-zigzag that attains this maximum will be called an $n$-polymer, or usually, simply a polymer.
 Thus, geodesics map to polymers under the scaling map. 
We will write~$\rho_{n;(x,t_1)}^{(y,t_2)}$ for {\em any} $n$-polymer  $(x,t_1)$ to $(y,t_2)$ (see Figure~\ref{f.scaling}), since we do not invoke polymer uniqueness results in this article. 

 This usage of the term `polymer' for `scaled geodesic' is apt for our study, due to the central role played by these objects. The usage is not, however, standard: the term `polymer' is often used to refer to typical realizations of the path measure in LPP models at positive temperature.
\subsubsection{Parameter settings in applications of results will be indicated in boldface.}\label{s.boldface}
Often we apply results involving several parameters. Typically these include the index $n$, times $t_1$ and $t_2$, and spatial locations such as $x$ and $y$.
Whenever we apply results, we will always state what these parameter settings are. When we do so, we will use boldface to indicate the variables in the result being applied, expressing these in terms of non-boldface variables, which assume their values from the current context. This device permits notational conflicts to be deescalated (so that notational choices need not proliferate, as they would were these conflicts to be   eliminated by other means).
\subsection{The scaling principle}
Let $(n,t_1,t_2) \in \N \times \R^2_<$ be a compatible triple. The quantity $n \tot$ is a positive integer, in view of the defining property~(\ref{e.ctprop}).
The scaling map $R_k: \R^2 \to \R^2$ has been defined whenever $k \in \N$, and thus we may speak of $R_n$  and $R_{n \tot}$.
The map $R_n$ is the composition of $R_{n \tot}$ and the transform $S_{\tot^{-1}}$ given by $\R^2 \to \R^2: (a,b) \to \big(a\tot^{-2/3},b\tot^{-1}\big)$. That is, the system of $n\tot$-zigzags is transformed into the system of $n$-zigzags
by an application of  $S_{\tot^{-1}}$.  Note from~(\ref{e.weightm}) that $\weight_{n;(x,t_1)}^{(y,t_2)} =  \tot^{1/3} \weight_{n \tot ;(x \tot^{-2/3},\kappa)}^{(y \tot^{-2/3},\kappa + 1)}$, where $\kappa = t_1 \tot^{-1}$; indeed this weight transformation law is valid for all zigzags, rather than just polymers, in view of~(\ref{e.weightzigzag}).
 
We may summarise these inferences by saying that the system of $n\tot$-zigzags, including their weight data, is transformed into the $n$-zigzag system, and its accompanying weight data, by the transformation $\big( a,b,c \big) \to \big( a \tot^{-1/3}  ,b \tot^{-2/3} , c \tot^{-1} \big)$, where the coordinates refer to the changes suffered in weight, 
 horizontal and vertical coordinates.
 This fact leads us to what we call the {\em scaling principle}. 
 
 \noindent{\em The scaling principle.} 
Let $(n,t_1,t_2) \in \N \times \R^2_<$ be a compatible triple.
 Any statement concerning the system of $n$-zigzags, including weight information, is equivalent to the corresponding statement concerning the system of $n\tot$-zigzags, provided that the following changes are made:
 \begin{itemize}
 \item the index $n$ is replaced by $n\tot$;
 \item any time is multiplied by $\tot^{-1}$;
 \item any weight is multiplied by $\tot^{1/3}$;
 \item and any horizontal distance is multiplied by $\tot^{-2/3}$.
 \end{itemize}
 
\subsubsection{The scaling principle applied:  uniform control on polymer weight for a general time-pair.}\label{s.spa}
Proposition~\ref{p.maxminweight}
provides a uniform control on polymer weights whose starting and ending points lie in $[x,x+1] \times \{0\}$ and  $[y,y+1] \times \{1\}$.
We now illustrate the scaling principle by using it to extend the proposition to treat the case where these intervals are replaced by   $[x,x+ \tot^{2/3}] \times \{ t_1 \}$ and  $[y,y + \tot^{2/3}] \times \{t_2\}$ for a general time pair $(t_1,t_2)$.


We phrase this more general result as an upper bound on the probability of a {\em polymer weight regularity} event. To define the new event, 
we again consider a compatible triple $(n,t_1,t_2) \in \N \times \R^2_<$. 
For $x,y \in \R$, $w_1,w_2 \geq 0$ and  $r > 0$, let $\maxmin_{n;([x,x+w_1],t_1)}^{([y,y+w_2],t_2)}( r )$ denote the event that, for all $(u,v) \in [0,w_1] \times [0,w_2]$,
$$
\Big\vert \, \tot^{-1/3} \weight_{n;(x+u,t_1)}^{(y+v,t_2)} +  2^{-1/2} \tot^{-4/3}  \big(y+v-x - u \big)^2 \, \Big\vert \, \leq \, r \, .
$$
We write $\neg \, A$ for the complement of the event $A$.
\begin{corollary}\label{c.maxminweight}
Let $(n,t_1,t_2) \in \N \times \R^2_<$ be a compatible triple for which 
$n\tot \in \N$ is at least $10^{29} \vee 2(c/3)^{-18}$.
 Let $x,y \in \R$ and  $a,b \in \N$ be   such that
$\big\vert x - y  \big\vert \tot^{-2/3} + \max\{ a,b\} - 1 \leq   6^{-1}  \rsc  (n\tot)^{1/18}$.
 Let $r \in \big[  34 \, , \, 4 (n \tot)^{1/18} \big]$.
Then
$$
 \PP \Big( \neg \, \maxmin_{n;([x,x+a \tot^{2/3}],t_1)}^{([y,y+b\tot^{2/3}],t_2)}(r) \Big) \leq  
 ab \cdot  400 C \exp \big\{ - c_1 2^{-10} r^{3/2} \big\} \, .
$$
\end{corollary}
\noindent{\bf Proof.} 
It is immediate from Proposition~\ref{p.maxminweight} that  
\begin{equation}\label{e.negmaxmin}
 \PP \Big( \neg \, \maxmin_{n;([x,x+1],0)}^{([y,y+1],1)}(r) \Big) \leq  
  400 C \exp \big\{ - c_1 2^{-10} r^{3/2} \big\} \, .
\end{equation}
when
$n \geq 10^{29} \vee 2(c/3)^{-18}$,
$\big\vert x - y  \big\vert \leq   6^{-1}  \rsc  n^{1/18}$ and $r \in \big[  34 \, , \, 4 n^{1/18} \big]$.

By the scaling principle and invariance under vertical shift, we know that
$$
 \PP \Big( \maxmin_{n;([x,x+ \tot^{2/3}],t_1)}^{([y,y+\tot^{2/3}],t_2)}(r)
  \Big) \, = \,  
 \PP \bigg( \maxmin_{n \tot;([0,1],0)}^{\big( \big[ (y-x)\tot^{-2/3} ,  (y-x)\tot^{-2/3} + 1 \big] , 1 \big)}(r)
  \bigg) \, .
$$
Consider (\ref{e.negmaxmin}) with parameter settings ${\bf n} = n \tot$, ${\bf x} = 0$, ${\bf y} =  (y-x)\tot^{-2/3}$ and ${\bf r} = r$. 
(This is the first use of boldface notation as specified in Subsection~\ref{s.boldface}. Its use here deprives ${\bf n} = n \tot$ of ambiguity of meaning, for example.)
We find then that  
\begin{equation}\label{e.probmaxmin}
 \PP \Big( \neg \, \maxmin_{n;([x,x+ \tot^{2/3}],t_1)}^{([y,y+\tot^{2/3}],t_2)}(r)
  \Big) \leq 
  400 C \exp \big\{ - c_1 2^{-10} r^{3/2} \big\}
\end{equation}
provided that 
$n\tot \geq 10^{29} \vee 2(c/3)^{-18}$,
$\big\vert x - y  \big\vert \tot^{-2/3} \leq   6^{-1}  \rsc  (n\tot)^{1/18}$ and $r \in \big[  34 \, , \, 4 (n\tot)^{1/18} \big]$.

This last bound is then summed over the $ab$ choices of pairs 
$$
\big({\bf x},{\bf y} \big) \in \big\{ x, x + \tot^{2/3}, \cdots,x+(a-1)\tot^{2/3}\big\} \times \big\{ y, y + \tot^{2/3},  \cdots, y+(b -1) \tot^{2/3} \big\}
$$
in order to obtain the corollary. Note that we hypothesise that 
$\big\vert x - y  \big\vert \tot^{-2/3} + \max\{ a,b\} - 1$ be at most   $6^{-1} \rsc  (n\tot)^{1/18}$ in order that the bound~(\ref{e.probmaxmin}) be valid for these parameter choices. \qed

Actually, when we apply this corollary in the present article, it will be in the case that $t_1 = 0$ and $t_2 = 1$, so in fact the use of the scaling principle is unnecessary in this regard.
It is useful, however, to have a general form for the corollary: for example, it is used in~\cite{NonIntPoly}.

\subsection{Polymers}

\subsubsection{Polymer concatenation}

Let $n \in \N$ and $(t_1,t_2,t_3) \in \R^3_\leq$ be such that $(n,t_1,t_2)$ and $(n,t_2,t_3)$ are compatible triples.
Let $x,y,z \in \R$.
In accordance with the convention stated in Subsection~\ref{s.polymer}, we may consider $n$-polymers $\rho_{n;(x,t_1)}^{(y,t_2)}$ and $\rho_{n;(y,t_2)}^{(z,t_3)}$.
The union of these two subsets of $\R^2$ is clearly an $n$-zigzag from $(x,t_1)$ and $(z,t_3)$.
In the union, the journey over the latter polymer follows that over the former. For this reason, we regard the union polymer as the concatenation of the two given polymers, and denote it by 
$\rho_{n;(x,t_1)}^{(y,t_2)} \circ \rho_{n;(y,t_2)}^{(z,t_3)}$. The new polymer's weight equals $\weight_{n;(x,t_1)}^{(y,t_2)} + \weight_{n;(y,t_2)}^{(z,t_3)}$.

\subsubsection{Polymer splitting}

Opposite to the operation of polymer concatenation is the splitting of a given polymer into two pieces. 
Let  $(n,t_1,t_2) \in \N \times \R^2_\leq$ be a compatible triple, and let $(x,y) \in \R^2$
satisfy  $y \geq x - 2^{-1} n^{1/3} \tot$. Let $t \in (t_1,t_2)$
be 
such that $(n,t_1,t)$ and $(n,t,t_2)$ are also compatible triples.
For any polymer $\rho_{n;(x,t_1)}^{(y,t_2)}$, we may select an element $(z,t) \in \rho_{n;(x,t_1)}^{(y,t_2)}$.
The removal of $(z,t)$ from $\rho_{n;(x,t_1)}^{(y,t_2)}$ creates two connected components. 
Taking the closure of each of these amounts to adding the point $(z,t)$ to each of them. The resulting sets are $n$-zigzags from $(x,t_1)$ to $(z,t)$, and from $(z,t)$ to $(y,t_2)$; in fact,
it is straightforward to see that they are $n$-polymers.

\subsubsection{Polymer crossing and rewiring}\label{s.polymercross}

We now make some comments about the implications of the event that two polymers cross. 
We do so for line-to-point polymers. 
For $\ovbar\coninit \in (0,\infty)^3$, let $f \in \initcond_{\ovbar\coninit}$.
Let $(n,t_1,t_2) \in \N \times \R^2_\leq$ be a compatible triple, and let $y \in \R$.
An $n$-zigzag $\phi$ from $(x,t_1)$, where $x \in \R$, to $(y,t_2)$, whose $f$-rewarded weight $\weight(\phi) + f(x)$ attains the maximum value 
$\weight_{n;(*:f,t_1)}^{(y,t_2)}$, is called an $f$-rewarded line-to-point polymer.
Such polymers are born free, but not equal: an endowment of $f$ is bestowed according to the place of birth. Pursuing a similar convention to that used in the point-to-point case, 
{\em any} such polymer will be denoted by $\rho_{n;(*:f,t_1)}^{(y,t_2)}$.

Suppose that two $f$-rewarded line-to-point polymers cross. That is, suppose that $(x_1,x_2) \in \R^2_<$
and $(y_1,y_2) \in \R_<^2$ are such that there exist such polymers, labelled  $\rho_{n;(*:f,t_1)}^{(y_1,t_2)}$ and  $\rho_{n;(*:f,t_1)}^{(y_2,t_2)}$ by our convention,
whose journeys are $(x_2,t_1) \to (y_1,t_2)$ and $(x_1,t_1) \to (y_2,t_2)$. The two polymers necessarily meet, and indeed the union of the horizontal segments of the two polymers also meet.
If $(z,t)$ is such a point of intersection, the operation of polymer splitting at $(z,t)$ may be applied to the two polymers, resulting in decompositions that may be respectively labelled $\rho_1 \circ \rho_2$
and $\rho_3 \circ \rho_4$. 
The zigzags $\rho_1$ and $\rho_3$
share their $f$-rewarded weights, $\weight(\rho_1) + f(x_2)$ and  $\weight(\rho_3) + f(x_1)$, because the weight maximality of  $\rho_1 \circ \rho_2$
and $\rho_3 \circ \rho_4$ each enforce one of the two inequalities between these quantities. 
Thus,  $\rho_1 \circ \rho_2$ and  $\rho_3 \circ \rho_2$
share their $f$-rewarded weight, and so do $\rho_3 \circ \rho_4$ and $\rho_1 \circ \rho_4$. The new, rewired, zigzags $\rho_3 \circ \rho_2$ and $\rho_1 \circ \rho_4$
are thus seen to be $f$-rewarded line-to-point polymers.

In summary, when two $f$-rewarded line-to-point polymers cross, the rewiring just undertaken results in an alternative pair of such polymers so that the old pair and the new share their set of starting and ending points.


\subsection{A simple lemma concerning polymer weight}

\begin{lemma}\label{l.basic}
Let $(n,t_1,t_2) \in \N \times \R^2_<$
be a compatible triple. 
\begin{enumerate}
\item The random function $(x,y) \to \weight_{n;(x,t_1)}^{(y,t_2)}$,
which is defined on the set of $(x,y) \in \R^2$ 
satisfying $y \geq x - 2^{-1} n^{1/3} \tot$, is continuous almost surely.
\item
Further consider an intermediate time $t \in (t_1,t_2)$
such that $(n,t_1,t)$ and $(n,t,t_2)$ are compatible triples.
Let $x,y,z \in \R$. Then 
$$
 \weight_{n;(x,t_1)}^{(y,t_2)} \geq  \weight_{n;(x,t_1)}^{(z,t)} +  \weight_{n;(z,t)}^{(y,t_2)} \, ,
$$
provided that these three weights are well defined. (The explicit conditions that ensure that the definitions make sense are 
$y \geq x - 2^{-1} n^{1/3} \tot$, $z \geq x - 2^{-1} n^{1/3} (t - t_1)$ and  $y \geq z - 2^{-1} n^{1/3} (t_2 - t)$.)
 \item Let  $\ovbar\coninit \in (0,\infty)^3$ and  $f \in \initcond_{\ovbar\coninit}$.  Suppose that $n \in \N$ satisfies
$n > 2^{-3/2} \coninit_1^3 \vee 8 (\coninit_2  + 1)^3$. Then  $[-1,1] \to \R: y \to \weight_{n;(*:f,0)}^{(y,1)}$ is almost surely finite and continuous.
\end{enumerate}
\end{lemma}
\noindent{\bf Proof. (1):}
By~(\ref{e.weightm}), it is enough to prove, for each $i,j \in \N$, $i \leq j$, that $M^1_{(x,i) \to (y,j)}$ is a continuous function of $(x,y) \in \R_\leq^2$.
Let $x_1$ and $x_2$ satisfy $x_1 \leq x_2 \leq y$.  It is a simple matter to verify that 
$$
 B(i,x_2) - B(i,x_1) \leq  M^1_{(x_1,i) \to (y,j)} -   M^1_{(x_2,i) \to (y,j)} \leq \sup_{k \in \llbracket i,j \rrbracket}     M^1_{(x_1,i) \to (x_2,k)} \, ; 
$$
$$
\textrm{and} \, \, \, \, \, \, \, \, \,  \sup_{k \in \llbracket i,j \rrbracket}     M^1_{(x_1,i) \to (x_2,k)} \leq \sum_{k=i}^j \Big(  \sup_{z \in [x_1,x_2]} B(k,z)  \, - \,   \inf_{z \in [x_1,x_2]} B(k,z)  \Big) \, .
$$
Continuity of $M^1_{(x,i) \to (y,j)}$ in the $x$-variable thus follows from continuity of the two-sided Brownian motions $B(k,\cdot)$.
This continuity in the $y$-variable follows similarly.

\noindent{\bf (2):} The weight $\weight_{n;(x,t_1)}^{(z,t)} + \weight_{n;(z,t)}^{(y,t_3)}$ of the concatenation of  two polymers $\rho_{n;(x,t_1)}^{(z,t)}$ and $\rho_{n;(z,t)}^{(y,t_3)}$
   offers a lower bound on  $\weight_{n;(x,t_1)}^{(y,t_2)}$.

\noindent{\bf (3):} In the proof of~\cite[Lemma~$4.6(2)$]{Patch}, which appears in~\cite[Appendix~$A$]{Patch}, it is noted that   $\mathsf{W}_{n;(*:f,0)}^{(y,1)}$ equals
$$
  2^{-1/2} n^{-1/3}  \sup_{u \in (-\infty,n + 2 n^{2/3}y]} \Big(    M^1_{(u , 0) \to (n+2 n^{2/3}y,n)} - n -   2 n^{2/3}y + u    +  h(u) \Big)    \, .
$$
Here,  
$h:\R \to \R \cup \{ - \infty \}$ is given by
$h(x) = 2^{1/2} n^{1/3}  f\big(n^{-2/3}x/2 \big)$.
Using a notation for unscaled line-to-point energy, namely
 $$
  M^1_{(*:g,0) \to (y,n)}  : = \sup \Big\{ E(\phi) + g(x): \phi \in  \upright^1_{(x,0) \to (y,n)} \, , \, x \leq y \Big\} \, ,
$$
the quantity  $\mathsf{W}_{n;(*:f,0)}^{(y,1)}$ is seen to equal $2^{-1/2} n^{-1/3}  M^1_{(*:g,0) \to (n+2 n^{2/3}y,n)}$ where the function $g$ is given by
 $g(u) =  - n -   2 n^{2/3}y + u    +  h(u)$. 
It is further noted in the same proof  that  $\limsup_{u \to -\infty} g(u)/\vert u \vert < 0$
is satisfied when $n > 2^{-3/2} \coninit_1^3$; and that, since $f \in \initcond_{\ovbar\coninit}$, the condition that $g(u) > -\infty$ for some $u \leq n+2n^{2/3}y$ (where $y \in [-1,1]$ is given) is verified when
$n \geq 8(\coninit_2  - y)^3$. We may thus apply~\cite[Lemma~$A.2$]{Patch} to learn that $\mathsf{W}_{n;(*:f,0)}^{(y,1)}$ almost surely assumes finite real values whenever $y \in [-1,1]$
under our present hypotheses. Regarding the continuity of this function of $y \in [-1,1]$, note first that, for $n$ given satisfying these hypotheses, the location of the maximizer  `$*:f$' is tight as $f$ varies over $\initcond_{\ovbar\coninit}$: this is a consequence of the square-root growth of $M^1_{(x,0) \to (y,n)}$ in the variable $y -x$, which is explained in and after equation~(28) of~\cite{Patch},
 which growth cannot compete with the linear decrease in the function~$g$.
The question of the continuity of $y \to \mathsf{W}_{n;(*:f,0)}^{(y,1)}$
has in essence  been reduced to the first part of the present lemma; we omit the details of this reduction.
\qed
\section{Line ensembles, their Brownian Gibbs property, and a key two-point estimate}\label{s.lineensembles}
In four subsections: we discuss the embedding of the narrow wedge polymer weight profile in a certain line ensemble; explain the Brownian Gibbs property enjoyed by a normalized version of that line ensemble, and an associated notion of {\em regular} ensemble; 
we gather the needed inputs from elsewhere, which assert  that the narrow wedge profile is embedded in a regular ensemble (Proposition~\ref{p.scaledreg}), 
and certain useful properties of regular ensembles  (Proposition~\ref{p.mega} and Lemma~\ref{l.pardom});  and we give the short, Brownian Gibbs, proof of the vital two-point estimate, Proposition~\ref{p.locreg}.
\subsection{Polymer weight profiles as the uppermost curves in line ensembles}\label{s.polymerweight}
Let $(n,t_1,t_2) \in \N \times \R^2_\leq$ be a compatible triple and let $x \in \R$.
Define the  forward polymer weight profile 
$$
\mc{L}_{n;(x,t_1)}^{\uparrow;t_2} (1,\cdot) :  \big[ x - 2^{-1} n^{1/3} \tot  ,\infty\big) \to \R
$$
 with base-point $(x,t_1)$ and end height $t_2$ by setting 
 $\mc{L}_{n;(x,t_1)}^{\uparrow;t_2}(1,y) = \weight_{n;(x,t_1)}^{(y,t_2)}$ for $y \geq x - 2^{-1} n^{1/3} \tot$. We call this weight profile `forward', and adorn the notation with the symbol $\uparrow$, to reflect that it is the spatial location~$y$
 associated to the more advanced time $t_2$ that is treated as the variable: we stand at $(x,t_1)$ and look forward in time to witness the weight profile as a function of $(\cdot,t_2)$.

Retaining the triple  $(n,t_1,t_2) \in \N \times \R^2_\leq$  but now fixing  $y \in \R$ (and treating $x \in \R$ as a variable),
we also introduce the {\em backward} polymer weight profile
$$
\mc{L}_{n;t_1}^{\downarrow;(y,t_2)}(1,\cdot)  :  \big( - \infty ,  y + 2^{-1} n^{1/3} \tot \big] \to \R
$$
 with base-point $(y,t_2)$ and end height $t_1$
by setting 
$\mc{L}_{n;t_1}^{\downarrow;(x,t_2)}(1,x) =  \weight_{n;(x,t_1)}^{(y,t_2)}$ for each $x \leq y +  2^{-1} n^{1/3} \tot$.
In other words, we now stand at $(y,t_2)$ and 
 look backwards in time at those polymers,  ending at our location, which begin at time (and height) $t_1$;  
 it is the weight profile of these polymers that is being recorded.

This new and elaborate looking notation may seem  merely to describe the already denoted weight profile  $\weight_{n;(x,t_1)}^{(y,t_2)}$, viewed as a function either of $x$ or $y$.
Its conceptual significance is suggested by our calling the argument of either profile `$(1,\cdot)$' rather than simply `$(\cdot)$'.
Indeed,  we will view~$\mc{L}_{n;(x,t_1)}^{\uparrow;t_2}$ as an ensemble of $n \tot + 1$ curves of which the lowest indexed curve, just defined, is the uppermost; and the backward object is just the same.
Either ensemble collectively has the Brownian Gibbs property, a fundamental tool for the analysis of the weight profile that is its uppermost curve. 

It is useful to retain a vivid picture of both of the processes $\mc{L}_{n;(x,t_1)}^{\uparrow;t_2} (1,\cdot)$ and $\mc{L}_{n;t_1}^{\downarrow;(x,t_2)}(1,\cdot)$
  as random curves that locally resemble Brownian motion but that globally follow the shape of a parabola.
 The parabola in question is  $- 2^{-1/2} (y-x)^2 \tot^{-4/3}$. The forward process adopts values of order $\tot^{1/3}$ for argument values $y$ that differ from $x$ by order $\tot^{2/3}$, and is forced downwards rapidly by parabolic curvature outside this region. When $\tot$ is small, for example, the weight profile is sharply peaked, and it broadens out as $\tot$ rises.   
(This description neglects the role of the index $n$, but roughly it develops accuracy as $n$ rises.) 

It is valuable to bring the forward and backward weight profiles for differing values of $\tot$ on to the same footing, by using a parabolic change of coordinates that, for example, flattens out the sharp peak witnessed when $\tot$ is small. The coordinate change will also bring the peak centre to the origin (from $x$ or $y$, according to the forward or backward case). The above weight profiles are already scaled objects, and so we introduce the term {\em normalized} to refer to the profiles viewed after this new, $\tot$-determined, change of coordinates.

Indeed, we  define the normalized forward polymer weight profile   
$$
 \scaledle_{n;(x,t_1)}^{\uparrow;t_2}(1,\cdot) :  \big[-  2^{-1} (n \tot)^{1/3}  , \infty \big) \to \R \, ,
$$
by setting 
$$
\scaledle_{n;(x,t_1)}^{\uparrow;t_2}\big( 1, z \big)
 = \tot^{-1/3} \mc{L}_{n;(x,t_1)}^{\uparrow;t_2}\big(  x + \tot
^{2/3} z \big) \, .
$$
Its backwards counterpart 
$\scaledle_{n;t_1}^{\downarrow;(y,t_2)}(1,\cdot)  :  \big( - \infty ,  2^{-1} (n \tot)^{1/3}  \big] \to \R$
is obtained by setting $\scaledle_{n;t_1}^{\downarrow;(y,t_2)}(1,z)$ equal to  $\tot^{-1/3} 
\mc{L}_{n;t_1}^{\downarrow;(y,t_2)}(1,   y + \tot
^{2/3} z \big)$.


Brownian motion is invariant under the parabolic change of coordinates, while the parabola  $x \to - 2^{-1/2} (y-x)^2 \tot^{-4/3}$ maps to $x \to -2^{1/2}x^2$.
Thus, our normalized processes should be pictured as locally Brownian as before, but with curvature dictated by the curve $-2^{-1/2} x^2$.
This picture in fact expands its domain of validity as the index increases, encompassing an expanding region about the origin, where the relevant indexing variable is now $n \tot$, rather than $n$. 
These heuristic comments find rigorously expressed counterparts in the next section.

\subsection{Brownian Gibbs ensembles and a regularity property}\label{s.bg}

Our weight profiles 
$\mc{L}_{n;(x,t_1)}^{\uparrow;t_2} (1,\cdot)$, 
 $\scaledle_{n;(x,t_1)}^{\uparrow;t_2}(1,\cdot)$, and their backward counterparts, may be embedded as uppermost curves in systems (or `ensembles') of random curves.
 (In  \cite[Figure~$4$]{NonIntPoly}, such a scaled and a normalized forward ensemble are illustrated.) 
  The normalized forward and backward ensembles satisfy the Brownian Gibbs property; moreover, these objects adhere well enough to the informal description of being locally Brownian and globally parabolic that we describe them as {\em regular} ensembles. 
 
In the next paragraphs: 
\begin{itemize}
\item we offer in outline the definition of the curves, indexed by higher values $k > 1$, in ensembles such as 
$\mc{L}_{n;(x,t_1)}^{\uparrow;t_2} (k,\cdot)$ and
 $\scaledle_{n;(x,t_1)}^{\uparrow;t_2}(k,\cdot)$ (our treatment is informal because only the uppermost ensemble curves, indexed by $k=1$, concern us in this article);
 \item we explain the implication for the uppermost curve of an ensemble being Brownian Gibbs;
\item and we specify in Definition~\ref{d.regularsequence} what it means for a Brownian Gibbs ensemble to be regular. 
\end{itemize} 
 
\subsubsection{Embedding weight profiles into ensembles as the uppermost curve.}\label{s.embedding}  For $(i,j) \in \N^2_\leq$ and $(x,y) \in \R^2_\leq$, recall that $M^1_{(x,i) \to (y,j)}$ is the energy maximum over staircases from $(x,i)$ to $(y,j)$. It has a counterpart $M^k_{(x,i) \to (y,j)}$, the maximum collective energy of a $k$-tuple of staircases with these endpoints satisfying a natural disjointness condition. The scaling relation~(\ref{e.weightm}) is extended to define the weight associated to the maximizing $k$-tuple after scaling.
This weight is then defined to equal the sum (over $i$) of the $k$ lowest indexed ensemble curves  $\mc{L}_{n;(x,t_1)}^{\uparrow;t_2} (i,\cdot)$ evaluated at $\cdot = y$.
 
We attempt a heuristic explanation of the upcoming appearance of the Brownian Gibbs property. The ensemble $\intint{n + 1} \times [0,\infty) \to \R: (k,y) \to  M^k_{(0,0) \to (y,n)} - M^{k-1}_{(0,0) \to (y,n)}$  (where $M^0 = 0$ is taken) is a microscopic counterpart to the scaled system  
$\intint{n+1} \times [0,\infty) \to \R: (k,y) \to  \mc{L}_{n;(0,0)}^{\uparrow;1}(k,y)$, because it is
the inverse image under the scaling map $R_n$ of the latter ensemble. 
In~\cite{O'ConnellYor}, this microscopic counterpart was identified in law as Dyson's Brownian motion with $n+1$ particles: a system of $n+1$ one-dimensional Brownian motions, all begun at the origin at time zero, and conditioned on mutual avoidance at all positive times. The linear scaling map $R_n$ preserves the diffusion rate of locally Brownian processes, so that the scaled ensemble $\mc{L}_{n;(0,0)}^{\uparrow;1}$ maintains the basic character of Dyson's Brownian motion: it is an ordered system of standard Brownian motions, conditioned on mutual avoidance, with boundary conditions that ensure that the first few uppermost curves are at unit distance, with the system's curves  globally following the parabola $-2^{-1/2} x^2$.


\subsubsection{Line ensembles and their Brownian Gibbs property.}  
Let $n \in \N$ and 
let $I \subseteq \R$ be closed.
A $\intint{n}$-indexed line ensemble defined on $I$ is a random collection of continuous curves  $\mc{L}:\intint{n} \times I \to \R$ specified under a probability measure $\PP$. The $i\textsuperscript{th}$ curve is thus $\mc{L}(i,\cdot): I \to \R$. (The adjective `line' has been applied to these systems perhaps because of their origin in such models as Poissonian LPP, where the counterpart object has piecewise constant curves. We will omit it henceforth.)
An ensemble is called {\em ordered} if $\mc{L}(i,x) > \mc{L}(i+1,x)$ whenever $i \in \intint{n-1}$ and $x$ lies in the interior of $I$.
The curves may thus assume a common value at any finite endpoint of $I$.
We will consider ordered ensembles that satisfy a key  condition called the Brownian Gibbs property.
We specify this property only in the special case that we need, in regard to the uppermost curve $\mc{L}(1,\cdot)$.

For $a,b \in \R$, $a < b$, and $y,z \in \R$,
let $\mc{B}^{[a,b]}_{y,z}$ denote the law of Brownian bridge $B:[a,b] \to \R$ with $B(a) = y$ and $B(b) =z$, given by conditioning standard Brownian motion to have these endpoints.

For any $[a,b] \subseteq I$, and $y,z \in \R$, 
consider
the conditional distribution of the marginal  process $\mc{L}(1,\cdot):[a,b] \to \R$ given that $\mc{L}(1,a) = x$ and $\mc{L}_N(1,b) =y$ and further given the form $g:[a,b] \to \R$ of $\mc{L}_N(2,\cdot)$ on $[a,b]$. 
The Brownian Gibbs property of $\mc{L}$ asserts that
this law
equals Brownian bridge $B$ under $\mc{B}^{[a,b]}_{y,z}$ conditioned by $B(u) > g(u)$ for all $u \in [a,b]$.


\subsubsection{Defining $(c,C)$-regular ensembles} 
The next definition specifies a  $(\bar\phi,\rsc,\rsC)$-regular ensemble from~\cite[Definition~$2.4$]{BrownianReg},
in the special case where the vector $\bar\phi$  equals  $(1/3,1/9,1/3)$.

\begin{definition}\label{d.regularsequence} 
%
Consider a Brownian Gibbs ensemble of  the form 
$$
\mc{L}: \intint{\nmac} \times \big[ - \xnmac , \infty \big) \to \R  \,  ,
$$
defined on a probability space under the law~$\PP$.
The number $\nmac = \nmac(\mathcal{L})$ of ensemble curves and the absolute value $\xnmac$ of the finite endpoint may take any values in $\N$ and $[0,\infty)$.
 
Let $\rsC$ and $\rsc$ be two positive constants. The ensemble $\mc{L}$
is said to be $(\rsc,\rsC)$-regular if the following conditions are satisfied.
\begin{enumerate}
\item {\bf Endpoint escape.} $\xnmac \geq  \rsc N^{1/3}$.
\item {\bf One-point lower tail.} If $z \in [ -\xnmac, \infty)$ satisfies $\vert z \vert \leq \rsc \nmac^{1/9}$, then
$$
\PP \Big( \mc{L} \big( 1,z\big) + 2^{-1/2}  z^2 \leq - s \Big) \leq \rsC \exp \big\{ - \rsc s^{3/2} \big\}
$$
for all $s \in \big[1, \nmac^{1/3} \big]$.
\item {\bf One-point upper tail.}  If $z \in [ -\xnmac, \infty)$ satisfies $\vert z \vert \leq \rsc \nmac^{1/9}$, then
$$
\PP \Big( \mc{L} \big( 1,z\big) +  2^{-1/2} z^2 \geq  s \Big) \leq \rsC \exp \big\{ - \rsc s^{3/2} \big\}
$$
for all $s \in [1, \infty)$.
\end{enumerate}
A Brownian Gibbs ensemble of the form 
$$
\mc{L}: \intint{\nmac} \times \big( -\infty , \xnmac  \big] \to \R
$$
is also said to be $(\rsc,\rsC)$-regular if the reflected ensemble $\mc{L}( \cdot, - \cdot)$ is. This is equivalent to the above conditions when instances of $[ - \xnmac, \infty)$
are replaced by $(-\infty, \xnmac]$.
\end{definition}

We will refer to these three regular ensemble conditions as $\rmreg(1)$,  $\rmreg(2)$ and $\rmreg(3)$.

\subsection{Inputs concerning regular ensembles}

\subsubsection{The normalized forward and backward ensembles are $(c,C)$-regular}\label{s.regular}

In Subsection~\ref{s.embedding}, we informally described how $\mc{L}_{n;(0,0)}^{\uparrow;1}$
is a globally parabolic object whose curves are mutually avoiding Brownian motions with a boundary condition suitable to ensuring that for example $\mc{L}_n(1,0)$ and $\mc{L}_n(1,0) - \mc{L}_n(2,0)$
are random but unit-order quantities. We capture this notion by asserting that this ensemble is {\em regular}. When the time-pair $(t_1,t_2) \in \R_\leq^2$  (as well as $x \in \R$) is general, it is the normalized ensemble  $\scaledle_{n;(x,t_1)}^{\uparrow;t_2}$ which is regular. 
Our assertion to this effect is~\cite[Proposition~$4.2$]{NonIntPoly}.  
\begin{proposition}\label{p.scaledreg}
Let $(n,t_1,t_2) \in \N \times \R^2_<$ be a compatible triple, and let $x \in \R$. 
The normalized forward weight profile   $\scaledle_{n;(x,t_1)}^{\uparrow;t_2}(1,\cdot)$, defined on  $\big[- 2^{-1} (n \tot)^{1/3}  , \infty \big)$, may be represented as the lowest indexed curve in an ensemble 
$$
 \scaledle_{n;(x,t_1)}^{\uparrow;t_2}: \intint{n \tot + 1} \times \big[- 2^{-1} (n \tot)^{1/3}  , \infty \big) \to \R 
$$
that enjoys the  Brownian Gibbs property.
Denoting this ensemble by $\mc{L}$, we naturally have  $\nmac(\mc{L}) = n \tot  + 1$ and $\xnmac = 2^{-1} (n \tot)^{1/3}$.

 There exist positive constants $\rsC$ and $\rsc$, which may be chosen independently of all such choices of the parameters $t_1$, $t_2$, $x$ and $n$,
such that the ensemble~$\mc{L}$ is $(\rsc,\rsC)$-regular. 

Similarly, the backward weight profile 
$\scaledle_{n;t_1}^{\downarrow;(y,t_2)}(1,\cdot)  :  \big( - \infty ,  2^{-1} (n \tot)^{1/3}  \big] \to \R$ may be embedded in an ensemble
$$
\scaledle_{n;t_1}^{\downarrow;(y,t_2)}(1,\cdot)  : \intint{n \tot + 1} \times  \big( - \infty ,   2^{-1} (n \tot)^{1/3}  \big] \to \R \, .
$$
This new ensemble also enjoys the properties just described for its forward counterpart, uniformly in the concerned parameters.
\end{proposition}

Hypothesis bounds in most of our results have been  
stated explicitly up to the appearance of 
two positive constants $c$ and $C$. The value is this pair is  fixed  by Proposition~\ref{p.scaledreg}.
Since bounding the constants would render hypotheses to be explicit, we mention that they are determined in~\cite[Appendix~$A.1$]{BrownianReg}
via Ledoux~\cite[(5.16)]{Ledoux} and Aubrun's~\cite[Proposition~$1$]{Aubrun} bounds on the lower and upper tail of the maximum eigenvalue of a matrix in the Gaussian unitary ensemble. In applications, we will harmessly suppose that $C \geq 1$ and $c \leq 1/2$.

\subsubsection{Basic properties of $(c,C)$-regular ensembles}\label{s.nontrivial}


Recall from Subsection~\ref{s.pwc} that $c_1 = 2^{-5/2}c \wedge \tfrac{1}{8}$.

\begin{proposition}\label{p.mega}
Suppose that  $\mc{L} = \mc{L}_N$, mapping either $\intint{N} \times [-\xnmac,\infty)$ or  $\intint{N} \times (-\infty,\xnmac]$, to $\R$,
is a $(\rsc,\rsC)$-regular ensemble, where $N \in \N$ and $\xnmac \geq 0$.
\begin{enumerate}
\item (Uniform curve lower bound) 
 Whenever   $(t,r,y) \in \R$ satisfy $N \geq (c/3)^{-18} \vee  6^{36}$, 
  $t \in \big[ 0 ,N^{1/18} \big]$, $r \in \big[ 2^{3/2} \, , \, 2N^{1/18} \big]$
  and $\vert y \vert \leq 2^{-1} \rsc N^{1/18}$,
$$
\PP \Big( \inf_{x \in [y-t,y+t]} \big( \mc{L}_N(1,x) + 2^{-1/2} x^2 \big) \leq - r \Big) \, \leq \, \Big( t \vee   (3 - 2^{3/2})^{-1} 5^{1/2}   \Big) \cdot 10 C \exp \big\{ - c_1 r^{3/2} \big\} \, .
$$
\item (No Big Max) 
For  $\vert y \vert \leq 2^{-1} c  N^{1/9}$, $r \in \big[0,4^{-1} \rsc  N^{1/9}\big]$, $t \in \big[ 2^{7/2} , 2 N^{1/3} \big]$ and $N \geq c^{-18}$,
$$
\PP \Big( \sup_{x \in [y-r,y+r]} \big( \mc{L}_N ( 1,x ) + 2^{-1/2}x^2 \big) \geq t \Big) \leq  (r + 1) \cdot  6  \rsC \exp \big\{ - 2^{-11/2} \rsc  t^{3/2} \big\} \, . 
$$
\item (Collapse near infinity) 
For $\eta \in (0,\rsc]$, let
$\ell = \ell_\eta:\R \to \R$ denote the even function which is affine on $[0,\infty)$ and has gradient $ - 5 \cdot 2^{-3/2} \eta \nmac^{1/9}$ on this interval, and which satisfies
$\ell(\eta \nmac^{1/9}) = \big( - 2^{-1/2} + 2^{-5/2} \big) \eta^2 \nmac^{2/9}$. If  $\nmac \geq 2^{45/4} \rsc^{-9}$, then  
\begin{eqnarray*}
 & & \PP \Big( \mc{L}_N \big(1,z\big) > \ell(z) \, \, \textrm{for some} \, \,  z \in  D \setminus \big[ - \eta  \nmac^{1/9} , \eta \nmac^{1/9} \big] \Big) \\
  & \leq &  
6C \exp \Big\{ - c \eta^3  2^{-15/4}   \nmac^{1/3} \Big\} \, .
\end{eqnarray*}
The set $D$ is the spatial domain of $\mc{L}$, either $[ - \xnmac , \infty )$  or  $(-\infty,\xnmac]$.  
\end{enumerate}
\end{proposition}

These three assertions are proved in~\cite{BrownianReg}. Respectively, they appear as, or are special cases of, the following results in that article:
Proposition~$A.2$,
Proposition~$2.28$,
and
Proposition~$2.30$.

A few words about the assertions' meaning. 
The first is a bound on the  lower tail of  the minimum value of  the lowest indexed curve on a compact interval.
The result is a strengthening of the defining property $\rmreg(2)$, which treats the one-point case. 
The second is a  similar strengthening of the one-point upper tail $\rmreg(3)$. In regard to the third, note that $\rmreg(2)$ and $\rmreg(3)$ do not assert that the lowest indexed curve hews to the parabola $-2^{-1/2}z^2$
globally, but only in an expanding region about the origin, of width $2cN^{1/9}$ centred at the origin, where $N$ is the ensemble curve cardinality. Proposition~\ref{p.mega}(3) offers a substitute control on curves far from the origin, showing them to decay at a rapid but nonetheless linear rate in the region beyond scale $N^{1/9}$.
  
Two further basic properties of regular ensembles are needed. One concerns a {\em parabolic symmetry} for whose explanation a little notation is helpful. 
Write $\para:\R \to \R$ for the parabola $\para(u) = 2^{-1/2} u^2$, and
let $l:\R^2 \to \R$ be given by $l(u,v) = - 2^{-1/2}v^2 - 2^{1/2}v(u-v)$. Note that $u \to l(u,v)$
is the tangent line of the parabola $u \to - \para(u)$ at the point $\big(v,-\para(v)\big)$. Note also that, for any $u,v \in \R$,
\begin{equation}\label{e.plp}
\para(u) = - l(u,v) + \para(u-v) \, .
\end{equation}

For $\xnmac \geq 0$, consider a $(c,C)$-regular ensemble $\mc{L}_N:\intint{N} \times [-\xnmac,\infty) \to \R$. 
For any $y_N > - \xnmac$, define $\lshift_{N,y_N}:\intint{N} \times [-\xnmac - y_N,\infty) \to \R$ to be the shifted ensemble given by 
$$
\lshift_{N,y_N}(i,u) = \mc{L}_{N}(i,u + y_N) - l(u+y_N,y_N)  \, .
$$

\begin{lemma}\label{l.pardom}
Let  $\rsc, \rsC > 0$ and $N \in \N$.
Suppose that  $\mc{L}_N:\intint{N} \times [-\xnmac,\infty) \to \R$
is  a    $\big(\rsc,\rsC\big)$-regular ensemble.
\begin{enumerate}
\item Whenever $y_N  \in \R$
satisfies $\vert y_N \vert \leq \rsc/2 \cdot N^{1/9}$, the ensemble $\lshift_{N,y_N}$ is    $\big(\rsc/2,\rsC\big)$-regular.
\item 
For any $[a,b] \subseteq [-\xnmac,\infty)$, and $y,z \in \R$,
the conditional distribution of the marginal process $\mc{L}_N(1,\cdot):[a,b] \to \R$ given that $\mc{L}_N(1,a) = y$ and $\mc{L}_N(1,b) =z$ stochastically dominates the Brownian bridge law  $\mc{B}^{[a,b]}_{y,z}$. 
\end{enumerate}
\end{lemma} 
\noindent{\bf Proof (1).} This is \cite[Lemma~$2.26$]{BrownianReg}.

\noindent{\bf Proof (2).}
By the Brownian Gibbs property of $\mc{L}_N$,
the conditional distribution of the marginal  process $\mc{L}_N(1,\cdot):[a,b] \to \R$ given that $\mc{L}_N(1,a) = x$ and $\mc{L}_N(1,b) =y$ and further given the form $g:[a,b] \to \R$ of $\mc{L}_N(2,\cdot)$ on $[a,b]$ equals Brownian bridge $B$ under $\mc{B}^{[a,b]}_{y,z}$ conditioned by $B(u) > g(u)$ for all $u \in [a,b]$. A stochastic monotonicity result~\cite[Lemma 2.6]{AiryLE} implies that this conditional distribution stochastically dominates the law specified by $g \equiv - \infty$. The latter law is~$\mc{B}^{[a,b]}_{y,z}$.  \qed



\subsection{Two-point estimate for the uppermost curve in a regular ensemble}

Here is the critical input for the proof of  Theorem~\ref{t.differenceweight}.

\begin{proposition}\label{p.locreg}
Suppose that  $\mc{L} = \mc{L}_N$, mapping either $\intint{N} \times [-\xnmac,\infty)$ or  $\intint{N} \times (-\infty,\xnmac]$, to $\R$,
is a $(\rsc,\rsC)$-regular ensemble, where $N \in \N$ and $\xnmac \geq 0$.
For $x \geq -\xnmac + 2$ and $t > 0$, define 
$$
\boundgood_t(x) = \bigcap_{x-2 \leq y \leq x+2} \Big\{ -t \leq  \mc{L}_N(1,y) + 2^{-1/2} y^2 \leq t  \, \Big\} \, .
$$
If $\vert x \vert \leq  2^{-1 }\rsc  N^{1/9}$, $\e \in (0,1]$, $K \geq 9$ and  
$N \geq 6^3 c^{-3}$, then
$$
 \PP \, \bigg( \, \Big\vert \mc{L}_N(1,x+\e) + 2^{-1/2} (x+\e)^2 - \mc{L}_N(1,x) + 2^{-1/2} x^2 \Big\vert \, \geq \,  K\e^{1/2} \, , \, \boundgood_{K/4}(x) \, \bigg) 
 $$
 is at most $2^{3/2} \pi^{-1/2} K^{-1} \exp \big\{ - 2^{-3} K^2 \big\}$.
\end{proposition}
 Theorem~\ref{t.differenceweight} (and Theorem~\ref{t.wlp.one}), 
 and Proposition~\ref{p.locreg} applied via Proposition~\ref{p.scaledreg},  all give expression to the one-half power law that governs polymer weight: when the endpoints of polymers are varied by short horizontal displacements of order~$\e$, the change in polymer weight has an order of~$\e^{1/2}$.   Proposition~\ref{p.locreg} is notably flexible, in that the parameters for horizontal scale, $\e$, and scaled fluctuation, $K$, may be selected without imposing any dependence on the lower bound demanded on the ensemble curve cardinality~$N$. This favourable feature comes at the price that the result gauges the small probability of high two-point difference only when we impose a global boundedness event $\boundgood_{K/4}(x)$
on the ensemble $\mc{L}_N$. We will have more to say about the role of  Proposition~\ref{p.locreg}  and the implications of its strengths and its drawback early in Section~\ref{s.polyweightreg}, when 
Theorem~\ref{t.differenceweight} is proved.

\noindent{\bf Proof of Proposition~\ref{p.locreg}.}
We first argue that an application of Lemma~\ref{l.pardom}(1) reduces to the case that $x = 0$.
To see this, note that~(\ref{e.plp}) implies that  $\lshift_{N,y_N}(i,u) + \para(u)$ equals
$\mc{L}_N(i,u + y_N) + \para(u+y_N)$, whenever $y_N > - \xnmac$. Selecting $y_N$ in Lemma~\ref{l.pardom}(1) to be $x$ in Proposition~\ref{p.locreg}, the proposition's conclusion is seen to be unchanged when the adjustments $x \to 0$ and $\mc{L}_N \to \lshift_{N,x}$ are made.
For this reason, we may, and will, consider only $x =0$. We must also work now under the assumption that $\mc{L}_N$ is $(c/2,C)$-regular. This information is used alongside the hypothesis  
$N \geq 6^3 c^{-3}$ to ensure that 
$c/2 \cdot  N^{1/3} \geq 3$, so that the interval $[-3,3]$ lies in the spatial domain $[-\xnmac,\infty)$ or $(-\infty,\xnmac]$~of~$\mc{L}_N$.

For a stochastic process $X$ whose domain of definition includes $[0,\e]$, define the events
$$
\down_{\e,K}[X] = \Big\{  X(\e) + 2^{-1/2} \e^2  \leq X(0) - K \e^{1/2} \Big\} 
$$
and  
$$
\up_{\e,K}[X] = \Big\{ X(\e) + 2^{-1/2} \e^2  \geq X(0) + K \e^{1/2} \Big\} \, . 
$$

For brevity, we write $W:[-\xnmac,\infty) \to \R$, $W(u) = \mc{L}_N(1,u)$, so that $W$ maps $[-\xnmac,\infty)$ or $(-\infty,\xnmac]$ to $\R$; and for the same reason, we also denote $t = K/4$. 
To obtain Proposition~\ref{p.locreg} for $x=0$, it is enough to verify two bounds:
\begin{equation}\label{e.downbound}
\PP \Big(  \down_{\e,K}[W] \, , \, \boundgood_t(0) \, \Big) 
\leq 2^{1/2} \pi^{-1/2} K^{-1} \exp \big\{ - 8^{-1} K^2 \big\} 
\end{equation}
and
\begin{equation}\label{e.upbound}
\PP \Big(  \up_{\e,K}[W] \, , \, \boundgood_t(0) \, \Big) 
\leq 2^{1/2} \pi^{-1/2} K^{-1} \exp \big\{ - 8^{-1} K^2 \big\} \, .
\end{equation}
\noindent{\em Deriving~(\ref{e.downbound}).}
Note that
\begin{eqnarray*}
 & & \PP \Big(  \down_{\e,K}[W] \, , \, \boundgood_t(0) \, \Big)  \leq  
 \PP \Big( \down_{\e,K}[W]  \, , \, 
  W(0) \leq t   \, , \, W(2) + 2^{3/2} \geq - t \, \Big) \\
   & \leq & \sup \bigg\{ \, \PP \Big( \, \down_{\e,K}[W] \, \Big\vert \, W(0) = y \, , \, W(2) = z \Big) : y \leq t \, , \, z \geq - t - 2^{3/2} \bigg\} \\
    & {\bm \leq} & \sup \bigg\{ \, \mc{B}_{y,z}^{[0,2]} \Big( \, \down_{\e,K}[B] \, \Big) : y \leq t \, , \, z \geq - t - 2^{3/2} \bigg\}   
     \,  =  \, \mc{B}_{t,-t - 2^{3/2}}^{[0,2]} \Big( \, \down_{\e,K}[B]  \, \Big) \, . 
\end{eqnarray*}
The third bound (which is highlighted ${\bm \leq}$) follows from Lemma~\ref{l.pardom}(2), while the final equality is a consequence of the coupling of Brownian bridge laws via affine shift.
Brownian bridge $B:[0,2] \to \R$ subject to $B(0) = t$
and $B(2) = -t - 2^{3/2}$ may be mapped via a further affine shift to Brownian bridge $\mc{B}_{0,0}^{0,2}$ with vanishing endpoint values. We thus see that 
\begin{eqnarray}
 & & \mc{B}_{t,-t - 2^{3/2}}^{[0,2]} \Big( \, \down_{\e,K}[B]  \, \Big) \, \leq \, \mc{B}_{t,-t - 2^{3/2}}^{[0,2]} \Big( \,  B(\e)  \leq t - K \e^{1/2}  \, \Big) \nonumber \\
 & = &  \mc{B}_{0,0}^{[0,2]} \Big( \,  B(\e)  \leq  - K \e^{1/2} \, + \, \e/2 \cdot \big(2t + 2^{3/2}\big) \, \Big) 
   \leq  \mc{B}_{0,0}^{[0,2]} \Big( \,   B(\e)  \leq  - 2^{-1} K \e^{1/2} \, \Big) \, , \label{e.infquant}
\end{eqnarray}
where the bound $\e \leq (2t + 2^{3/2})^{-2} K^2$ that permits the final inequality is due to $K = 4t$, $t \geq 2^{1/2}$ (in the guise that $2t + 2^{3/2} \leq 4t$) and $\e \leq 1$. Since $B(\e)$ under $\mc{B}_{0,0}^{[0,2]}$
is normally distributed with mean zero and variance $2^{-1}\e(2-\e) \leq \e$, the right-hand side of~(\ref{e.infquant}) is at most $\nu ( x ,\infty) \leq (2\pi)^{-1/2} x^{-1} e^{-x^2/2}$ with $x=K/2$, where $\nu$ is the standard normal distribution, and the inequality~\cite[Section~$14.8$]{Williams} is classical.  We have  obtained~(\ref{e.downbound}).


\noindent{\em Deriving~(\ref{e.upbound}).}
In summary of the last argument, $W$ may fall no more suddenly between zero and $\e$ than does a suitable Brownian bridge that begins at time zero. But neither may it rise suddenly between this same pair of times, because such a rise is a fall when viewed from right to left, and this fall entails a fall that is at least as great on the part of a suitable Brownian bridge whose rightmost time is~$\e$. That is, the argument for~(\ref{e.upbound}) is  almost symmetrical to that for~(\ref{e.downbound}). It takes the form
\begin{eqnarray*}
 & & \PP \Big(  \up_{\e,K}[W] \, , \, \boundgood_t(0) \, \Big)  \leq  
 \PP \Big( \up_{\e,K}[W]  \, , \, 
  W(-2) + 2^{3/2} \geq - t   \, , \, W(\e) + 2^{-1/2}\e^2 \leq  t \, \Big) \\ 
   & \leq & \sup \bigg\{ \, \PP \Big( \, \up_{\e,K}[W] \, \Big\vert \, W(-2) = y \, , \, W(\e) = z \Big) : y \geq - t - 2^{3/2} \, , \, z \leq  t -  2^{-1/2}\e^2  \bigg\} \\
    & {\bm \leq} & \sup \bigg\{ \, \mc{B}_{y,z}^{[-2,\e]} \Big( \, \up_{\e,K}[B] \, \Big) : y \geq -t - 2^{3/2} \, , \, z \leq  t -  2^{-1/2}\e^2  \bigg\}   
     \,  =  \, \mc{B}_{-t - 2^{3/2},t -  2^{-1/2}\e^2 }^{[-2,\e]} \Big( \, \up_{\e,K}[B]  \, \Big) \\
      & {\bm =}  &
 \mc{B}_{0,0}^{[0,2 + \e]} \Big( \,  B(\e) \leq  - K \e^{1/2} + 2^{-1/2}\e^2 + \tfrac{\e}{2 + \e} ( 2t + 2^{3/2} -  2^{-1/2}\e^2 ) \, \Big) \leq  \mc{B}_{0,0}^{[0,2]} \Big( \,   B(\e)  \leq  - 2^{-1} K \e^{1/2} \, \Big)  \, ,
\end{eqnarray*}
where Lemma~\ref{l.pardom}(2) again implies the third bound, which is again highlighted ${\bm \leq}$. To obtain the latter equality (which is highlighted~${\bm =}$), consider the map that sends $B:[-2,\e] \to \R$ to $W:[0,2+\e] \to \R$, where $W$ is the affine shift of $[0,2+\e] \to \R: y \to B(\e - y)$ 
for which $W(0) = W(2 + \e) = 0$, and note that this map sends the law  $\mc{B}_{-t - 2^{3/2},t -  2^{-1/2}\e^2 }^{[-2,\e]}$ to $\mc{B}_{0,0}^{[0,2 + \e]}$ and the event $\up_{\e,K}[B]$ to the event that   $W(\e) \leq - K \e^{1/2} + 2^{-1/2}\e^2 + \tfrac{\e}{2 + \e} ( 2t + 2^{3/2} -  2^{-1/2}\e^2 )$. 
It is now  $\e \leq 1 \wedge  (t + 2^{1/2} +  2^{-1/2})^{-2} 2^{-2}K^2$ that permits the final displayed  inequality,
 with the assumptions that $t \geq 2^{1/2} + 2^{-1/2}$ (in the guise that $t + 2^{1/2} + 2^{-1/2}$ is at most~$2t$), $K = 4t$ and $\e \leq 1$ implying this condition. Since the last displayed quantity is the right-hand side of~(\ref{e.infquant}), the rest of the derivation of~(\ref{e.downbound}) applies, and~(\ref{e.upbound}) is obtained.

In summary,~(\ref{e.downbound}) and~(\ref{e.upbound})
imply Proposition~\ref{p.locreg} with $x=0$, to whose derivation the proof of this proposition has been reduced. \qed

\section{Collective control on polymer weights: the proof of Proposition~\ref{p.maxminweight}}\label{s.collective}

Proposition~\ref{p.maxminweight} will be proved first in the case that the parameter $n \in \N$ is even; the case of odd $n$ will then be obtained by reducing to the preceding one. 
In this section, we prove the proposition for even $n$, and then explain in outline the argument for $n$ odd. The detailed argument for odd $n$ appears in 
the online Appendix~$A$. 

In order to permit the reduction argument, we will in fact prove the result when $n \in 2\N$
under slightly weaker hypotheses and with a slightly stronger conclusion.

That is, we will now prove the proposition in the case that $n \in \N$ is even, when $n$ is supposed to satisfy 
$n \geq 10^{29} \vee 2(c/3)^{-18}$; with $x,y \in \R$ 
satisfying
$\big\vert x - y  \big\vert \leq   3^{-1} 2^{-2/3} \rsc  n^{1/18}$; and with $t \in \big[  33 \, , \, 4 n^{1/18} \big]$. We will derive the conclusions~(\ref{e.weightuvsupremum}) and~(\ref{e.weightuvinfimum}) where these bounds are strengthened by the respective replacements $2^{-10} \to 2^{-9}$ and $2^{-3} \to 2^{-5/2}$ on their right-hand sides.

The planar line segment with endpoints $(x,0)$ and $(y,1)$
crosses the horizontal line whose  height is one-half
at the location $(x+y)/2$. We record this location in the form $z =  (x+y)/2$. 
We will prove the lower-tail bound~(\ref{e.weightuvinfimum}) by bounding the polymer weight 
 $\weight_{n;(x+u,0)}^{(y+v,1)}$ for any given $(u,v) \in [0,1]^2$ below by considering routes from $(x+u,0)$ to $(y+v,1)$ that pass via $(z,1/2)$. The two polymer weights in the lower bound
 concern journeys from $(z,1/2)$ to $(x+u,0)$ (after time reversal) and from $(z,1/2)$ to $(y+v,1)$. Rooting in this way at $(z,1/2)$, we may gauge the probability of low weight values for polymers emanating from $(z,1/2)$ and ending in a compact interval at time zero or one by applying Proposition~\ref{p.mega}(1) to the duration one-half normalized ensembles rooted at $(z,1/2)$, of forward or backward type according to whether a time one or time zero endpoint is being considered. 

Thus, we let $u,v \in [0,1]$. Note that
\begin{eqnarray*}
 & & 2^{-1/2} \big(2^{2/3}(z -x -u) \big)^2 + 2^{-1/2} \big(2^{2/3}(y+ v - z) \big)^2 \\
  & = & 2^{-1/6}  (y+v-x - u)^2 +2^{-1/6} ( u + v  )^2 \, .
\end{eqnarray*}
Note further that
\begin{eqnarray*}
 \weight_{n;(x+u,0)}^{(y+v,1)} & \geq & \weight_{n;(x+u,0)}^{(z,1/2)} + \weight_{n;(z,1/2)}^{(y+v,1)}  
  =  \mc{L}_{n;0}^{\downarrow;(z,1/2)}(1,x+u) + 
 \mc{L}_{n;(z,1/2)}^{\uparrow;1}(1,y+v) \\
  & = & 2^{-1/3} \scaledle_{n;0}^{\downarrow;(z,1/2)}\big(1, 2^{2/3} (x+u - z) \big) + 
  2^{-1/3} \scaledle_{n;(z,1/2)}^{\uparrow;1}\big(1,  2^{2/3} (y+v -z) \big) \, .
\end{eqnarray*}
The  inequality here invokes Lemma~\ref{l.basic}(2) with ${\bf t_1} = 0$, ${\bf t_2} = 1$, ${\bf t} = 2^{-1}$, ${\bf x} = x+u$, ${\bf y} = y + v$ and ${\bf z} = z$. The two equalities invoke the definitions of the four ensembles whose top curves are being evaluated.
Regarding the inequality, we may note that the use of  Lemma~\ref{l.basic}(2) entails that certain bounds on ${\bf z}$ be satisfied. We omit reference to these bounds now because they are anyway implicated later in the argument, when we come to analyse the above two right-hand terms. Finally, note that, in order to enable our use of  Lemma~\ref{l.basic}(2) with the choice ${\bf t} = 2^{-1}$,
we have imposed in the present case of the proof of Proposition~\ref{p.maxminweight} that $n \in \N$ be even, in order that the triples $(n,0,2^{-1})$ and  $(n,2^{-1},1)$ be compatible.

Adding to the above inequality the $2^{-1/3}$-rd multiple of the display that preceded it, we find that
\begin{eqnarray}
 & & \weight_{n;(x+u,0)}^{(y+v,1)} + 2^{-1/2}  (y+v-x - u)^2 + 2^{-1/2} ( u + v  )^2 \nonumber \\
  & \geq & 2^{-1/3} \Big( \scaledle_{n;0}^{\downarrow;(z,1/2)}\big(1, 2^{2/3} (x+u - z) \big) +  2^{-1/2} \big(2^{2/3}(z -x -u) \big)^2 \Big) \label{e.weightbound} \\
   & &  \quad + \,   2^{-1/3} \Big( \scaledle_{n;(z,1/2)}^{\uparrow;1}\big(1,  2^{2/3} (y+v -z) \big) + 2^{-1/2} \big(2^{2/3}(y+ v - z) \big)^2 \Big) \, . \nonumber
\end{eqnarray}

Note that
\begin{eqnarray}
 & &  \PP \bigg(  \inf_{u \in [0,1]} \Big( \scaledle_{n;0}^{\downarrow;(z,1/2)}\big(1, 2^{2/3} (x+u - z) \big) +  2^{-1/2} \big(2^{2/3}(z -x -u) \big)^2 \Big) \leq - s \bigg) \label{e.pinf} \\
 & =  &  \PP \bigg(  \inf_{u \in [0,2^{2/3}]} \Big( \scaledle_{n;0}^{\downarrow;(z,1/2)}\big(1, 2^{2/3} (x - z) + u \big) +  2^{-1/2} \big(2^{2/3}(x - z) + u \big)^2 \Big) \leq - s \bigg) \, . \nonumber
\end{eqnarray}
The latter term equals
$$
\PP \Big( \inf_{x' \in [y-t',y+t']} \big( \mc{L}_N(1,x') + 2^{-1/2}(x')^2 \big) \leq - r \Big)
$$
when $t' = 2^{-1/3}$; when $t' \geq 2^{-1/3}$, the new expression is an upper bound. 
Here, we take $\mc{L}_N =  \scaledle_{n;0}^{\downarrow;(z,1/2)}$, $y =  2^{2/3}(x - z) + 2^{-1/3}$ and $r = s$.

We seek then to apply 
Proposition~\ref{p.mega}(1) 
to the $(c,C)$-regular ensemble $\mc{L}_N =  \scaledle_{n;0}^{\downarrow;(z,1/2)}$, 
doing so 
with the choice of $({\bf r},{\bf y},{\bf t}) = (r,y,t')$.
It is 
Proposition~\ref{p.scaledreg}
that permits this choice of ensemble. If we set $t_1 = 0$ and $t_2 = 1/2$, so that $\tot = 1/2$, the number of curves in the ensemble   $\mc{L}_N =  \scaledle_{n;t_1}^{\downarrow;(z,t_2)}$ equals $n \tot  + 1$ and thus is at least $n \tot$.
For the application to be valid, our parameters must thus satisfy 
\begin{equation}\label{e.tothyp}
\tot n \geq 1 
  \vee  (c/3)^{-18} \vee  6^{36} \, , \, 
   2^{-1/3} \leq (n\tot)^{1/18} \, , \, s \in \big[ 2^{3/2} \, , \, 2(n \tot)^{1/18} \big] \, ,
\end{equation}
and  $\big\vert 2^{2/3}(x - z) + 2^{-1/3} \big\vert = \big\vert 2^{2/3}(x - y)/2 + 2^{-1/3} \big\vert \leq \rsc/2 \cdot (n \tot)^{1/18}$.
 From this application of Proposition~\ref{p.mega}(1), we find that the probability in~(\ref{e.pinf}) is at most 
\begin{equation}\label{e.twoub}
  \Big( 2^{-1/3} \vee (3 - 2^{3/2})^{-1} 5^{1/2} \Big) \cdot 10 C  \exp \big\{ - c_1 s^{3/2} \big\} \, . 
\end{equation}

By applying   Proposition~\ref{p.mega}(1) to the ensemble $\mc{L}_N =  \scaledle_{n;(z,1/2)}^{\uparrow;1}$, 
with ${\bf y}$ now chosen equal to  $2^{2/3}(y - z) + 2^{-1/3}$, and with $({\bf r},{\bf t})$ again set to be $\big(s,2^{-1/3}\big)$, we find that the quantity
$$
 \PP \bigg(  \inf_{v \in [0,1]} \Big( \scaledle_{n;(z,1/2)}^{\uparrow;1}\big(1,  2^{2/3} (y+v -z) \big) + 2^{-1/2} \big(2^{2/3}(y+ v - z) \big)^2 \Big) \leq - s \bigg) 
$$
is also bounded above by~(\ref{e.twoub}). This application of the proposition requires in addition to the bounds in~(\ref{e.tothyp}) that $\big\vert 2^{2/3}(y - z) + 2^{-1/3} \big\vert = \big\vert 2^{2/3}(y - x)/2 + 2^{-1/3} \big\vert \leq \rsc/2 \cdot (n \tot)^{1/18}$.

Since $2^{-1/2}(u + v)^2 \leq 2^{3/2}$ whenever $u,v \in [0,1]$, we find from the inequality~(\ref{e.weightbound}) and the upper bound by~(\ref{e.twoub}) on the two probabilities that
\begin{eqnarray*}
 & & \PP \bigg(  \inf_{u,v \in [0,1]} \Big( \weight_{n;(x+u,0)}^{(y+v,1)} + 2^{-1/2}  (y+v-x - u)^2 \Big) \leq - 2 \cdot 2^{-1/3}s \, - \, 2^{3/2} \bigg) \\
  & \leq & 2  \Big(  2^{-1/3} \vee (3 - 2^{3/2})^{-1} 5^{1/2} \Big)  \cdot 10 C  \exp \big\{ - c_1 s^{3/2} \big\} \, .
\end{eqnarray*}
Setting $t = 2 \cdot 2^{2/3}s$ and using $s \geq 2^{5/6}$ (so that $t \geq 2^{2/3}s + 2^{3/2}$), 
and noting  $20  \Big( 2^{-1/3} \vee  5^{1/2} (3 - 2^{3/2})^{-1}    \Big) \leq 261$, 
we obtain~(\ref{e.weightuvinfimum}).

Let $\high_{n;([x,x+1],0)}^{([y,y+1],1)}(t)$ denote the event that  $\sup_{u,v \in [0,1]} \Big( \weight_{n;(x+u,0)}^{(y+v,1)} + 2^{-1/2} ( y+v - x - u )^2 \Big) \geq t$.
Clearly it is this event whose probability we must bound above as we turn to derive~(\ref{e.weightuvsupremum}). The event entails the presence of a high weight polymer that crosses a square, but both of its endpoints may have exceptional locations. The derivation of~(\ref{e.weightuvsupremum}) will proceed by noting that, typically, one of the endpoints can be made typical. Indeed, when the $\high$ event occurs,
so that a high weight polymer runs between say $(x+U,0)$ and $(y+V,1)$, where the pair $(U,V) \in [0,1]^2$ is random, a fairly high weight polymer will also typically exist between the deterministic location $(2x-y,-1)$
and $(y+V,1)$. This is because a rather high lower bound on the weight of such a polymer is obtained by considering the pair of polymers, from $(2x-y,-1)$ to $(x+U,0)$, and from $(x+U,0)$ to $(y+V,1)$, whose weights are typically not too low, and high. The probability of the presence of this fairly high weight polymer may then be gauged by the No Big Max Proposition~\ref{p.mega}(2), because this event entails that the duration-two normalized forward ensemble rooted at $(2x-y,-1)$ assumes a high value within a compact interval.

To begin implementing this approach, we consider the event
$\high_{n;([x,x+1],0)}^{([y,y+1],1)}(t)$, and let $(U,V) \in [0,1]^2$ be the lexicographically minimal pair of $(u,v) \in [0,1]^2$  that realize this event (a definition which makes sense because, by Lemma~\ref{l.basic}(1), the set of such pairs is closed).

Now {\em reset} the value of $z$ to be $2x -y$, so that 
the planar line segment interpolating $(z,-1)$ and $(y,1)$ passes through $(x,0)$.
Note that
\begin{eqnarray}
 & &  2^{-1/2} (z -x -U)^2 + 2^{-1/2} (x+ U - y - V)^2 \label{e.eqfive} \\ 
  & = & 2^{-3/2} (z - y - V)^2    + 2^{-1/2} \big( 2U^2 + V^2/2 - 2UV  \big) \, . \nonumber
  \end{eqnarray}

Let $\notlow_{n;(z,-1)}^{([x,x+1],0)}(t/2)$ denote the event that  $\inf_{u \in [0,1]} \Big( \weight_{n;(z,-1)}^{(x+u,0)} + 2^{-1/2} ( x +u - z )^2 \Big) \geq - t/2$.
Note that
$$
 \neg \, \, \notlow_{n;(z,-1)}^{([x,x+1],0)}(t/2) \,  = \, \bigg\{ \inf_{u \in [0,1]} \Big( \scaledle_{n;(z,-1)}^{\uparrow;0}\big(1,x+u - z\big) + 2^{-1/2} ( x +u - z )^2 \Big) < - t/2   \bigg\} \, .
$$
We apply   Proposition~\ref{p.mega}(1)
to the ensemble $\scaledle_{n;(z,-1)}^{\uparrow;0}$,  the choice admissible by  Proposition~\ref{p.scaledreg},
in order to find an upper bound on the probability of the displayed event. The application is made with 
${\bf y} = x - z + 1/2$, ${\bf t} = 1/2$ and ${\bf r} = t/2$.
Since the ensemble $\scaledle_{n;(z,-1)}^{\uparrow;0}$ has $n+1$, and therefore at least $n$, curves, we see that the next bounds suffice for the application to be made:
 $n  \geq (c/3)^{-18} \vee  6^{36}$, 
  $1/2 \leq n^{1/18}$, $t/2 \in \big[ 2^{3/2} \, , \, 2n^{1/18} \big]$
  and $\big\vert x - z + 1/2 \big\vert =  \big\vert y - x + 1/2 \big\vert
 \leq 2^{-1} \rsc   n^{1/18}$. 
 (In regard to the latter, recall that $z$ is now $2x-y$.) We learn from the application that
$$
 \PP \Big(  \neg \, \, \notlow_{n;(z,-1)}^{([x,x+1],0)}(t/2) \Big) \, \leq \,   
  (3 - 2^{3/2})^{-1} 5^{1/2}  \cdot 10 C \exp \big\{ - c_1 2^{-3/2} t^{3/2} \big\} \, .
$$

When $\notlow_{n;(z,-1)}^{([x,x+1],0)}(t/2) \cap \high_{n;([x,x+1],0)}^{([y,y+1],1)}(t)$ occurs, consider the concatenation $\rho_{n;(z,-1)}^{(x+U,0)} \circ \rho_{n;(x+U,0)}^{(y+V,1)}$ 
of any pair of polymers with the endpoints implied by our convention governing this notation. The concatenation 
has weight at least 
\begin{eqnarray*}
 & &  \Big( -t/2  - 2^{-1/2} ( x + U - z )^2 \Big) \, + \, \Big( t -  2^{-1/2} ( y+V - x - U )^2  \Big) \\
 & \geq & t/2 -  2^{-3/2} (z - y - V)^2  -   2^{-1/2} \big( 2U^2 + V^2/2 - 2UV  \big)  \\
  & \geq & t/2 -  2^{-3/2} (z - y - V)^2    - 5 \cdot 2^{-3/2} 
 \, ,
\end{eqnarray*}
the displayed inequalities due to~(\ref{e.eqfive}) and $U ,V \in [0,1]$.
Thus,
$$
  \high_{n;([x,x+1],0)}^{([y,y+1],1)}(t) \, \cap \, \notlow_{n;(z,-1)}^{([x,x+1],0)}(t/2) \subseteq \bigg\{ \sup_{v \in [0,1]} \Big( \weight_{n;(z,-1)}^{(y+v,1)} + 2^{-3/2} (z - y - v)^2 \Big)  \geq t/2  - 5 \cdot 2^{-3/2}  \bigg\} \, .
$$
The right-hand event equals 
\begin{eqnarray}
 & &  \bigg\{ \sup_{v \in [0,1]} \Big( \mc{L}_{n;(z,-1)}^{\uparrow;1}(1,y+v) + 2^{-3/2} (z - y - v)^2 \Big)  \geq t/2 - 5 \cdot 2^{-3/2} \bigg\} \label{e.lnrhs} \\ 
  & = &   \bigg\{ \sup_{v \in [0,1]} \Big(  \scaledle_{n;(z,-1)}^{\uparrow;1}\big(1, 2^{-2/3} (y+v - z) \big) + 2^{-1/2} \big( 2^{-2/3}(z - y - v) \big)^2 \Big)  \geq 2^{-4/3}t - 5 \cdot 2^{-11/6} \bigg\} \nonumber \\
  & = &   \bigg\{ \sup_{v \in [0,2^{-2/3}]} \Big(  \scaledle_{n;(z,-1)}^{\uparrow;1}\big(1, 2^{-2/3} (y - z) + v \big) + 2^{-1/2} \big( 2^{-2/3}(z - y) - v \big)^2 \Big)  \geq 2^{-4/3}t - 5 \cdot 2^{-11/6} \bigg\} \, . \nonumber 
\end{eqnarray}

We now apply the No Big Max Proposition~\ref{p.mega}(2) 
to the ensemble $\mc{L}_N =  \scaledle_{n;(z,-1)}^{\uparrow;1}$,  the application permitted by Proposition~\ref{p.scaledreg}.
In the application, we set ${\bf y} =  2^{-2/3} (y - z) + 2^{-5/3}$, ${\bf r} = 2^{-5/3}$ and ${\bf t} =  2^{-4/3}t - 5 \cdot 2^{-11/6}$. The curve cardinality of the ensemble in queston is $2n+1 \geq 2n$.
As such, it is sufficient for the application to be valid that
$$
\big\vert 2^{-2/3} (y - z) + 2^{-5/3} \big\vert =  
\big\vert 2^{1/3} (y - x) + 2^{-5/3} \big\vert  \leq c/2 \cdot (n\tot)^{1/9} \, , \, 2^{-5/3} \leq \rsc/4 \cdot (n\tot)^{1/9} 
$$
and $2^{-4/3}t - 5 \cdot 2^{-11/6} \in \big[ 2^{7/2} , 2 (n\tot)^{1/3} \big]$
as well as $\tot n \geq c^{-18}$ where here $\tot$ equals $2$ (in accordance with the time-pair $t_1 = -1$ and $t_2 = 1$ being considered).
This application tells us that the $\PP$-probability of the event~(\ref{e.lnrhs}) is at most
\begin{equation}\label{e.fivethirds}
   \big( 2^{-5/3} + 1 \big) \cdot  6  \rsC \exp \big\{ - 2^{-11/2} \rsc  \big(  2^{-4/3}t - 5 \cdot 2^{-11/6} \big)^{3/2} \big\} \leq  8  \rsC \exp \big\{ -     2^{-9} \rsc t^{3/2} \big\}  \, ,
\end{equation}
where we used $2^{-4/3}t - 5 \cdot 2^{-11/6} \geq 2^{-7/3} t$ when $t \geq 5 \cdot 2^{1/2}$.

We find then that
$$
 \PP \Big( \high_{n;([x,x+1],0)}^{([y,y+1],1)}(t) \Big) \,  \leq \, \PP \Big( \high_{n;([x,x+1],0)}^{([y,y+1],1)}(t) \cap \notlow_{n;(z,-1)}^{([x,x+1],0)}(t/2) \Big) \, + \, \PP \Big( \neg \, \, \notlow_{n;(z,-1)}^{([x,x+1],0)}(t/2)  \Big)
$$
is bounded above by
$$
  8  \rsC \exp \big\{ -     2^{-9} \rsc t^{3/2} \big\} 
 \, + \,
  10 (3 - 2^{3/2})^{-1} 5^{1/2}    C \exp \big\{ - c_1 2^{-3/2} t^{3/2} \big\} \, .
$$
We now  use $c_1 \leq \rsc$ and $8 + 10 (3 - 2^{3/2})^{-1}  5^{1/2}     \leq 139$ to obtain~(\ref{e.weightuvsupremum}). This completes the proof of the slightly strengthened form of Proposition~\ref{p.maxminweight} specified at the start of this section in the case that $n \in \N$ is even. 

It remains to derive the proposition in the case that $n \in \N$ is odd. The argument is simple enough and here we outline it; but its details occupy a few pages and are deferred to 
the online Appendix~$A$.

The new case is handled by reducing to that where $n$ is even; this involves retreating to unscaled coordinates. The argument harnesses a trivial bound on unscaled energy, namely that $M^1_{(z_1,0) \to (z_2,m)} \leq M^1_{(z_1,0) \to (z_2,m+1)}$ whenever $m \in \N$ and $z_1,z_2 \in \R$ satisfy $z_1 \leq z_2$. Presently given are $n \in \N$ odd; $x,y \in \R$; and $u,v \in [0,1]$. First we seek to obtain~(\ref{e.weightuvsupremum}). Applying the unscaled inequality with $m=n$, $z_1 = 2 n^{2/3}(x+u)$ and $z_2 = n + 2 n^{2/3}(y+v)$, and writing the resulting bound in scaled coordinates via~(\ref{e.weightmzeroone}),
we find that
\begin{equation}\label{e.twoweights}
 \weight_{n;(x + u,0)}^{(y + v,1)} \, \leq \,  \big( 1 + n^{-1} \big)^{1/3} \, \weight_{n+1;(x' + u',0)}^{(y' + v',1)} \, + \, 2^{-1/2} n^{-1/3} \, ,
\end{equation}
where the new variables $x', y',u',v' \in \R$ satisfy $x' + u' = \big( 1 + n^{-1} \big)^{-2/3}(x+u)$
and $y' + v' =  \big( 1 + n^{-1} \big)^{-2/3}(y+v) - 2^{-1}(n+1)^{-2/3}$.
Thus the weight term in~(\ref{e.weightuvsupremum}) may be replaced by  $\weight_{n+1;(x' + u',0)}^{(y' + v',1)}$ at the expense of adjustments that are negligible when $n$ is high, with the new term exceeding the old. Similarly may the parabolic term $Q(y+v - x -u)$ be replaced by $Q(y'+v' - x' -u')$.  
Thus we are able to derive~(\ref{e.weightuvsupremum}) in the present case of odd $n$, by using~(\ref{e.twoweights}) and the already derived~(\ref{e.weightuvsupremum}) with ${\bf n} = n + 1$, ${\bm x} = \big( 1 + n^{-1} \big)^{-2/3} x$ and ${\bf y} =  \big( 1 + n^{-1} \big)^{-2/3}y - 2^{-1}(n+1)^{-2/3}$; and with ${\bf t}$ set to be $t$ multiplied by a factor that is slightly less than one. 
Proposition~\ref{p.maxminweight} has slightly stronger hypotheses on parameters, and slightly weaker conclusions, in its general case $n \in \N$ than in the form obtained when $n$ is even.
It is these changes that permit us to apply the even version of~(\ref{e.weightuvsupremum}) 
in order to obtain~(\ref{e.weightuvsupremum}) when $n$ is odd. 
The deriviation of~(\ref{e.weightuvinfimum}) with $n$ odd similarly proceeds from the unscaled inequality noted above and reduces to~(\ref{e.weightuvinfimum}) with ${\bf n} = n-1$. The details are similar enough that we defer further discussion to 
the online Appendix~$A$. \qed


The reader may have noticed that every application of a result concerning a $(c,C)$-regular ensemble in the preceding proof invoked Proposition~\ref{p.scaledreg} 
in order to justify that the ensemble in question indeed enjoys this property. Every subsequent such application is no different, and henceforth we omit mention of Proposition~\ref{p.scaledreg}'s role.

\section{Polymer weight regularity: proving Theorem~\ref{t.differenceweight}}\label{s.polyweightreg}

We begin this section by introducing a notation $\paradelta\weight$ to denote the difference in {\em parabolically adjusted} weight of two polymers 
crossing between the opposite endpoints of two intervals. (The superscript $\cup$ is intended to suggest a parabola.)
When $(x_1,x_2)$ and $(y_1,y_2)$ belong to $\R^2_\leq$, 
we set  
\begin{equation}\label{e.paradelta}
\paradelta\weight_{n;([x_1,x_2],0)}^{([y_1,y_2],1)} \, := \, \Big( \weight_{n;(x_2,0)}^{(y_2,1)} + 2^{-1/2}(y_2 - x_2)^2 \Big) - \Big( \weight_{n;(x_1,0)}^{(y_1,1)}  + 2^{-1/2}(y_1 - x_1)^2 \Big) \, .
\end{equation}
We will abuse this notation when one of the concerned intervals collapses a point, writing for example
 $\paradelta\weight_{n;(x_1,0)}^{([y_1,y_2],1)}  = \weight_{n;(x_1,0)}^{(y_2,1)} - \weight_{n;(x_1,0)}^{(y_1,1)} +  2^{-1/2}(y_2 - x_1)^2 -  2^{-1/2}(y_1 - x_1)^2$. 
 
 We also write $y + U = \{ y+u:u \in U \}$ when $y \in \R$ and $U \subset \R$.

Any integer is called 
a dyadic rational of scale zero. A dyadic rational of scale $i \in \N$, $i \geq 1$, has for the form $p 2^{-i}$ where $p \in \Z$ is odd. 
A dyadic interval of scale $i \in \N$ is a closed interval of length~$2^{-i}$ that has an endpoint which is a dyadic rational of scale~$i$.

Recall from Subsection~\ref{s.spa} the polymer weight regularity events $\pwr$.

\begin{proposition}\label{p.dyadic}
Let $n \in \N$, $n \geq (4/c)^9$, let $x,y \in \R$ satisfy  $\big\vert x - y  \big\vert \leq \rsc/4 \cdot n^{1/9}$,  
 and let $K_0 \geq 9$. 
 Let $i,k_0 \in \N$ with $i \geq k_0$ (which is at least one). 
Consider the quantities
\begin{equation}\label{e.deltaone}
\PP \Big( \sup \Big\vert  \paradelta\weight_{n;(x+ z,0)}^{(y+ U,1)} \Big\vert \geq   K_0 2^{-i/2}  \, , \,  \maxmin_{n;([x,x+1],0)}^{([y-2,y+3],1)}(K_0/4)   \Big)  
\end{equation}
 and 
\begin{equation}\label{e.deltatwo}
\PP \Big( \sup \Big\vert  \paradelta\weight_{n;(x+ U,0)}^{(y+ z,1)} \Big\vert \geq  K_0 2^{-i/2}   \, , \,  \maxmin_{n;([x-2,x+3],0)}^{([y,y+1],1)}(K_0/4)  \Big) \, , 
\end{equation}
where in each case the supremum is taken over all 
dyadic rationals $z \in [0,2^{-k_0}]$ of scale $i$ and dyadic intervals $U \subseteq [0,2^{-k_0}]$ of scale at least $i$. 
Each of these quantities is at most 
$2^{2(i-k_0)} \exp \big\{ - 2^{-3}  \Kzero^2 \big\}$.
\end{proposition}
This result is naturally a key technical ingredient in the proof of Theorem~\ref{t.differenceweight}.
The theorem concerns polymer weight differences when small changes are made in endpoint locations. The proposition is similar, but restricts to endpoint locations that are dyadic rationals.
The theorem will follow from the proposition by noting that one may skip between the two nearby endpoint locations $x$ and $x + \e$ (and similarly for $y$ and $y + \e$)
by jumping through a possibly infinite sequence of intermediate dyadic rational  locations where no dyadic scale need be visited more than twice and the minimal dyadic scale is of the order of the difference~$\e$.
This property of this `stepping stone' sequence makes the union bounds over the estimate in Proposition~\ref{p.dyadic} manageable. 
This inference of the theorem from the proposition is similar to the derivation of the Kolmogorov continuity criterion, in which moment bounds on the difference of a stochastic process between a generic pair of times imply H\"older continuity of the process: see~\cite[Thoerem~$8.13$]{Durrett10}.


One further aspect of the plan for proving Theorem~\ref{t.differenceweight} deserves mention before we proceed. 
Note that in the theorem the parameter $R$, which measures the degree of polymer weight fluctuation,
must verify an $n$-dependent upper bound. Although this bound in a sense is insignificant for the purpose of analysing high $n$ behaviour,
Proposition~\ref{p.dyadic} has been stated so that there is no comparable hypothesis:
the quantities $K_0$ and $2^{-i}$ are counterparts to $R$ and $\e$, and the proposition holds for all high~$n$, where the lower bound on $n$ deteriorates neither as the dyadic scale~$2^{-i}$ decreases, nor as the parameter $K_0$ increases.
Now this is a valuable property, but it comes at a certain price, about which more in a moment. 
The reason that the property is valuable is that, for a given high value of $n$, Proposition~\ref{p.dyadic}
may be applied as $i$ increases to infinity, while in the meantime, $K_0$ also increases; indeed, this is how we will derive 
 Theorem~\ref{t.differenceweight}, with a union bound over the infinite number of applications of the proposition indexed by~$i$ being controllable due to the ongoing increase in $K_0$.
  As for the price to be paid, we mention that,
in order that the property obtains, it has been necessary in the events whose probabilities are gauged in~(\ref{e.deltaone}) and~(\ref{e.deltatwo})
to include the global polymer weight regularity events $\pwr$; this in turn is because the proposition's proof  will invoke the two-point Proposition~\ref{p.locreg}, in which the spatial scale $\e > 0$ is permitted to be arbitrarily small without forcing the ensemble curve cardinality $N$ to be higher in an $\e$-dependent way, a favourable circumstance which is only possible at the expense of introducing the global regularity event $\boundgood_{K/4}$, counterpart to the above $\pwr$ events, into the probability upper bound in that result.
In any case, it would seem that we can derive Theorem~\ref{t.differenceweight} from 
Proposition~\ref{p.dyadic}  only by invoking the $\pwr$ event. In fact, we do impose this event, but, at the very end of the derivation, we gauge the 
 probability of the complementary event $\neg \, \pwr$ via Corollary~\ref{c.maxminweight}; thus, this probability cost is paid only once, rather than with each of the infinitely many applications of Proposition~\ref{p.dyadic}.

\noindent{\bf Proof of Proposition~\ref{p.dyadic}.}
Consider a given dyadic interval $U \subseteq [0,2^{-k_0}]$ of a scale $j$ such that $j  \geq i$ and write $U = [u_1,u_2]$.
For a given dyadic rational $z \in [0,2^{-k_0}]$, note that 
$$
 \paradelta\weight_{n;(x+ z,0)}^{(y+ U,1)}  =  \mc{L}^{\uparrow;1}_{n;(x+z,0)}\big(1,y+ u_2 \big) + 2^{-1/2}(y+u_2 - x+z)^2 -  \mc{L}^{\uparrow;1}_{n;(x+z,0)}\big(1, y + u_1 \big)  - 2^{-1/2}(y+u_1 - x+z)^2  
$$
and that
\begin{equation}\label{e.twoomegas}
  \paradelta\weight_{n;(x+ U,0)}^{(y+ z,1)}   = \mc{L}^{\downarrow;(y+z,1)}_{n;0}\big(1,x+u_2 \big)
+ 2^{-1/2}(x+u_2 - y+z)^2  
   -  \mc{L}^{\downarrow;(y + z,1)}_{n;0}\big(1,x + u_1 \big) - 2^{-1/2}(x+u_1 - y+z)^2   \, . 
\end{equation}
Let $K > 0$ be a parameter that we will later specify. We now seek to apply 
Proposition~\ref{p.locreg}
in order to bound the $\PP$-probability that the last two displayed quantities
exceed $K \vert U \vert^{1/2} = K 2^{-j/2}$.
In order to do so for the first quantity, the proposition will be applied to the ensemble
$\mc{L}_N =  \scaledle^{\uparrow;1}_{n;(x+z,0)}$.
 Proposition~\ref{p.locreg}'s 
 parameters in this case are set 
 ${\bf x} = y+u_1 - x - z$, ${\bm \e} = 2^{-j}$ and ${\bf K} = K$.
 Note that the event whose probability is bounded above by the proposition is a subset of $\boundgood_{K/4}({\bf x})$. With the present choice of ensemble $\mc{L}_n$,   the event  $\boundgood_{K/4}({\bf x})$ equals $\maxmin_{n;(\{x+z\},0)}^{([y+u_1 - 2,y+u_1 + 2],1)}(K/4)$. 
 (Note that here $\{ x+z \}$ is a singleton set, so that a space of polymers emanating from the point $(x+z,0)$ is at stake.)
 Thus,  $\maxmin_{n;([x,x+1],0)}^{([y-2,y+3],1)}(K/4) \subseteq \boundgood_{K/4}({\bf x})$.
 As such, 
 Proposition~\ref{p.locreg} implies that
 \begin{eqnarray}
  & &\PP \bigg(   \Big\vert  \paradelta\weight_{n;(x+ z,0)}^{(y+ U,1)} 
  \Big\vert   \geq  K 2^{-j/2}  \, , \,   \maxmin_{n;([x,x+1],0)}^{([y-2,y+3],1)}(K/4) \bigg) \nonumber \\
  & \leq & \PP \bigg(   \Big\vert  \paradelta\weight_{n;(x+ z,0)}^{(y+ U,1)} 
  \Big\vert   \geq  K 2^{-j/2}  \, , \,  \boundgood_{K/4}({\bf x})  \bigg)   \leq  
 2^{3/2} \pi^{-1/2}   K^{-1}  \exp \big\{ - 2^{-3}  K^2 \big\}  
  \, . \label{e.concone}
 \end{eqnarray}
The hypotheses of Proposition~\ref{p.locreg} that are invoked to obtain this bound are 
$$
n + 1 \geq 6^3 c^{-3} \, , \, \vert  y + u_1 - x - z \vert \leq 2^{-1} \rsc n^{1/9} \, , \,  2^{-j} \in (0,1] 
  \, \, \, \textrm{and} \, \, \, K \geq  9 \, ;
$$
note that $j \geq 0$ 
is used to validate the third of these.

Another application of  Proposition~\ref{p.locreg}
 is made in regard to the quantity~(\ref{e.twoomegas}).
On this occasion, the ensemble $\mc{L}_n$ is set equal to $\scaledle^{\downarrow;(y+z,1)}_{n;0}$, and the parameters are set:
 ${\bf x} = x + u_1  - y - z$, ${\bm \e} = 2^{-j}$ and ${\bf K} = K$. 
 In this instance, the event  $\boundgood_{K/4}({\bf x})$ equals  
$\maxmin_{n;([x+u_1 - 2,x + u_1 + 2],0)}^{(y+z,1)}(K/4)$, so that  $\maxmin_{n;([x-2,x+3],0)}^{([y,y+1],1)}(K/4) \subseteq \boundgood_{K/4}({\bf x})$.
The outcome of the application in this case is the conclusion that
\begin{eqnarray}
 & & \PP \bigg(  \Big\vert \paradelta\weight_{n;(x+ U,0)}^{(y+ z,1)} \Big\vert \geq K 2^{-j/2} \, , \,    \maxmin_{n;([x-2,x+3],0)}^{([y,y+1],1)}(K/4)  \bigg)  \nonumber  \\
  & \leq &
 2^{3/2} \pi^{-1/2}    K^{-1}  \exp \big\{ - 2^{-3}  K^2 \big\} 
  \, , \label{e.conctwo}
\end{eqnarray}
while the hypothesis $\vert x + u_1  - y - z \vert \leq 2^{-1} \rsc  n^{1/9}$ is used to make the application.

The dyadic scale $j$ is at least the scale $i$ by assumption, and we now denote $\ell = j - i \geq 0$. 
We recall from the statement of Proposition~\ref{p.dyadic} that we consider a parameter $\Kzero \geq 9$. 
We now set the value of our parameter~$K$, choosing it to be $\Kzero 2^{(j-i)/2}$. Our conclusions~(\ref{e.concone}) and~(\ref{e.conctwo}) tell us that
$$
 \PP  \Big( \Big\vert  \paradelta\weight_{n;(x+ z,0)}^{(y+ U,1)} 
  \Big\vert   \geq  \Kzero 2^{-i/2}  \, , \,  \maxmin_{n;([x,x+1],0)}^{([y-2,y+3],1)}(K_0/4) \Big) 
$$
and 
\begin{equation}\label{e.backfluc}
 \PP  \Big(  \Big\vert \paradelta\weight_{n;(x+ U,0)}^{(y+ z,1)} \Big\vert \geq  \Kzero 2^{-i/2}  \, , \,  \maxmin_{n;([x-2,x+3],0)}^{([y,y+1],1)}(K_0/4) \Big)
\end{equation}
are both at most
$2^{3/2} \pi^{-1}   \Kzero^{-1} \exp \big\{ - 2^{-3} \cdot 2^\ell \Kzero^2 \big\}$, 
where we used $K/4 \geq \Kzero/4$.

We will now sum the stated bound on quantities of the form~(\ref{e.backfluc})
in order to find an upper bound on the expression~(\ref{e.deltatwo})
where recall that, in this expression, the supremum is taken over choices of scale~$i$ dyadic rational $z \in [0,2^{-k_0}]$ and dyadic interval $U = [u_1,u_2] \subset [0,2^{-k_0}]$ of scale at least~$i$. 
Indeed, 
since the number of dyadic rationals of scale $i$ in $[0,2^{-k_0}]$
is at most $2^{i-k_0} + 1$, while the number of scale~$j$ dyadic intervals contained in $[0,2^{-k_0}]$ equals $2^{j-k_0}$,
we may sum over $j \geq i$, also using that $\Kzero \geq 2^{3/2}$, to find that~(\ref{e.deltatwo}) is at most 
$$
 (2^{i - k_0} +1) \cdot 2^{i-k_0} \cdot \tfrac{1}{e-2} \cdot
 2^{3/2} \pi^{-1} K_0^{-1}     \exp \big\{ - 2^{-3}  \Kzero^2 \big\} \, .
$$
Using $i \geq k_0$ and $K_0 \geq 9$,
 this quantity is found to be bounded above by
$$
 2^{2(i-k_0)+ 5/2} 3^{-2} (e-2)^{-1}  \pi^{-1}  \exp \big\{ - 2^{-3}  \Kzero^2 \big\}  \leq 2^{2(i-k_0)}  \exp \big\{ - 2^{-3}  \Kzero^2 \big\}    \, .
$$
 We have obtained the upper bound claimed in Proposition~\ref{p.dyadic} on the probability~(\ref{e.deltatwo}) and, since
 this upper bound on~(\ref{e.deltaone}) is similarly obtained,  we have completed the proof of 
this proposition. \qed

\medskip


We now state an estimate also needed for the proof of Theorem~\ref{t.differenceweight}. Suppose that $n \in \N$  satisfies
\begin{equation}\label{e.twentyfiveapp}
n \geq 10^{29} \vee 2(c/3)^{-18} \, , \,  \big\vert x - y  \big\vert + 4 \leq   6^{-1}  \rsc  n^{1/18} \, \, \textrm{and} \, \, r_0 \in \big[  34 \, , \,  4 n^{1/18} \big] \, .
\end{equation}
Corollary~\ref{c.maxminweight} with ${\bf n} = n$, ${\bf t_1} = 0$,  ${\bf t_2} = 1$, ${\bf x} = x-2$, ${\bf y} = y-2$, ${\bf a} = {\bf b} = 5$ and ${\bf r} = r_0$ implies that 
\begin{equation}\label{e.maxminbound}
 \PP \Big( \neg \,  \maxmin_{n;([x-2,x+3],0)}^{([y-2,y+3],1)}(r_0) \Big) \leq 5^2 \cdot  400 \, C \exp \big\{ - c_1 2^{-10} r_0^{3/2} \big\} \, .
\end{equation}
Note that $\big\vert x - y  \big\vert + 4 \leq   6^{-1}  \rsc  n^{1/18}$
is implied by  $\big\vert x - y  \big\vert \leq  2^{-1} 6^{-1}  \rsc  n^{1/18}$
and $4 \leq 2^{-1} 6^{-1} \rsc  n^{1/18}$ and the latter is implied by 
$n \geq \big( 48 c^{-1} \big)^{18}$.

\medskip

\noindent{\bf Proof of Theorem~\ref{t.differenceweight}.}
Write $I = [x,x+\e]$ and $J = [y,y+\e]$, and let $u \in I$ and $v \in J$ be arbitrary. Recalling that $\e \leq 2^{-4}$ is less than one, we
consider the  binary expansion $u - x = \sum_{j=1}^\infty s_j 2^{-j}$. If the expansion is not unique, we choose its terminating version for definiteness.
Let the increasing sequence $u_0,u_1,\cdots$ 
enumerate the set 
$$
 \bigg\{ \,  x +  \sum_{j=1}^k s_j 2^{-j} : k \in \N \, \bigg\} \, .
$$
This sequence may be finite or infinite. It begins
$u_0 = x$, and, when it is infinite, each term $u_n$ equals the sum of $x$ and the quantity given by the truncation of the binary expansion of $u-x$
that contains $n$ instances of the digit one.
Let $n_1 \in \N \cup \{ \infty \}$ denote the maximal index of a term in the $u$-sequence. If $n_1 < \infty$, then $u_{n_1} = u$, and if $n_1 = \infty$, then $u_n \nearrow u$, and we set $u_\infty = u$.

Similarly, we specify an increasing sequence  $v_0,v_1,\cdots$ by replacing $(x,u)$ by $(y,v)$ above, and let~$n_2$ denote the maximal index of a term in the $v$-sequence. If $n_2 = \infty$, set $v_\infty = v$.

Call the planar points $\big\{ (u_i,0): 0 \leq i \leq n_1 \big\}$ lower pegs, and the points  $\big\{ (v_i,1): 0 \leq i \leq n_2 \big\}$ upper pegs.
Think of a  planar line segment {\em cord}  that runs in the first instance between $(u_0,0) = (x,0)$
and $(v_0,1) = (y,1)$.  Permitting the cord to be pegged at its lower and upper end to any of the pegs, we see that the cord begins in its leftmost possible location.
The rightmost available location is given by lower peg $(u_\infty,0) = (u,0)$ and upper peg $(v_\infty,1) = (v,1)$.
We now specify a possibly infinite sequence of cord moves by which  the cord will achieve, or at least converge towards, this rightmost location.
Let $(L_i,U_i) \in \{0,\cdots,n_1\} \times \{0,\cdots,n_2 \}$ denote the indices of lower and upper peg locations at step $i \in \N$, where the original location is indexed by $i = 0$,
so that $(L_0,U_0) = (0,0)$.
Let $k \in \N$ and consider the value of $(L_k,U_k)$. If this value is $(n_1,n_2)$, then the cord movement is complete and the value of  $(L_{k+1},U_{k+1})$
is not recorded. In the other case, there are two possible moves for the cord at the next step: a lower peg advance, in which  
$L_{k+1} = L_k + 1$ and $U_{k+1} = U_k$, or an upper peg advance, in which $U_{k+1} = U_k + 1$ and $L_{k+1} = L_k$. It may be that one of these moves is inadmissible, because $L_k = n_1$, which renders the lower peg advance unavailable, or $U_k = n_2$, which does likewise for the upper peg advance.
If this is so, then $(L_{k+1},U_{k+1})$ is set equal to the value given by the only available advance.
In the case where both advances are possible, note that each move entails displacing a peg to the right by a distance of the form $2^{-i}$ for some $i \in \N$.
The decision of which advance to make is taken so that this distance is the larger for the available two advances,
with say the upper advance being made if the distances are equal.
In this way, we specify the value of $(L_{k+1},U_{k+1})$; we may also record the dyadic scale, $D_{k+1} \in \N$, of the advance associated to this index increase $k \to k+1$:
this scale is the value of $i \in \N$ such that the peg displacement made in the peg advance resulting in the new peg locations $(L_{k+1},U_{k+1})$ equals $2^{-i}$.
  
The sequence of location pairs  $\big( u_{L_k} , v_{U_k} \big)$ either reaches its terminal state $(u,v)$ after finitely many moves, or it converges to this state as $k$ increases.
Note also that the dyadic scale sequence $D_1,D_2,\cdots$ is a non-decreasing $\N$-valued sequence that assumes any given value at most twice. This sequence depends on the pair $(u,v)$, and we may indicate this dependence by writing $D_k = D_k(u,v)$.

Moreover, we define a dyadic scale $i_0 \in \N$ by setting $i_0 \geq 0$ to be minimal such that $2^{i_0} \e \geq 1$. Then, since the first peg is displaced by a distance $2^{-D_1}$ which is at most $\e < 2^{1 - i_0}$, we see that $D_1 \geq i_0$.
We see then that the dyadic scales   $\big\{ D_j: j \geq 1 \big\}$ of the distances of peg moves for the advancing cord grow linearly from around $i_0$: 
\begin{equation}\label{e.djlb}
D_j \geq i_0 - 1 + j/2 \, \, \,  \textrm{for} \, \, \,  j \geq 1 \, . 
\end{equation}
To any cord location we may associate the parabolically adjusted weight $\weight_{n;(z,0)}^{(z',1)} + 2^{-1/2} (z' - z)^2$ of the polymer whose endpoints $(z,0)$ and $(z',1)$ are the pegs to which the cord is pinned.
Using the notation~(\ref{e.paradelta}), we may further set, for $k \geq 0$, 
\begin{equation}\label{e.wkplusone}
W_{k+1} = W_{k+1}(u,v) = \paradelta\weight_{n;([u_{L_k},u_{L_{k+1}}],0)}^{([v_{U_k},v_{U_{k+1}}],1)}   \, ,
\end{equation}
this being the difference in parabolically adjusted weight of this polymer as a result of the cord move from its index $k$ to $k+1$ location.

Lemma~\ref{l.basic}(1) implies that, for given $n$, $\lim_k \weight_{n;(L_k,0)}^{(U_k,1)} = \weight_{n;(u,0)}^{(v,1)}$.
Thus,
\begin{equation}\label{e.sumwk}
 \weight_{n;(u,0)}^{(v,1)} -   \weight_{n;(x,0)}^{(y,1)} + 2^{-1/2}(v-u)^2 - 2^{-1/2}(y-x)^2 = \sum_{k=1}^\infty W_k \, ,
\end{equation}
where it is understood that the right-hand sum may have only finitely non-zero terms.

For each $k \in \N$ (at least one), we set $W^*[k]$ to be the supremum of the values of $\vert W_k \vert$ over all choices of $(u,v)$ in our construction. We now argue that there exists a pair $(u^*[k],v^*[k]) \in I \times J$ such that $W^*[k] = \big\vert W_k(u^*[k],v^*[k]) \big\vert$; in other words, we argue that the supremum specifying $W^*[k]$ is attained. First note that, if $W^*[k]=0$, any choice of $(u^*[k],v^*[k]) \in I \times J$ works. Suppose instead then that $W^*[k] > 0$. For any given $(u,v) \in I \times J$, the value $W_k(u,v)$ takes one of the forms $\paradelta\weight_{n;(x + \Theta,0)}^{(y + z_1,1)}$ or $\paradelta\weight_{n;(x + z_2,0)}^{(y + \Lambda,1)}$, where in the former case, $\Theta$ is a dyadic rational interval and $z_1 \in [0,\e]$ is a dyadic rational whose scale is by our construction at most that of~$\Theta$; and, in the latter, these statements are equally true of $\Lambda$ and $z_2$.
Now the number of such expressions for which the length of the interval $\Theta$ or $\Lambda$ in the indexing pair $(\Theta,z_1)$ or $(z_2,\Lambda)$ exceeds an arbitrary given positive value~$\delta$ is finite. On the other hand, since the map $I \times J \to \R: (u,v) \to \weight_{n;(u,0)}^{(v,1)}$ is uniformly continuous,
there exists a random value $\delta > 0$ such that none of the parabolically adjusted weight differences whose indexing pair contains an interval $\Theta$ or $\Lambda$
of length less than $\delta$
has value exceeding the positive quantity $W^*[k]/2$. The value $W^*[k]$ is thus seen to be the maximum of a certain finite set of such weight differences, and thus, it is indeed achieved.  

For given $(u,v) \in I \times J$, there are many choices of $(u^*[k],v^*[k]) \in I \times J$, and we pick the lexicographically minimal pair for definiteness.  We further specify the dyadic scale
$D^*[k] = D_k(u^*[k],v^*[k])$.

By~(\ref{e.djlb}), $D_\ell(u,v) \geq i_0  - 1 + k$ provided that $\ell \geq 2k$, whatever the value of $(u,v) \in [x,x+\e] \times [y,y+\e]$ may be.
Thus, when $\ell \in \{ 2k,2k+1 \}$,
the quantity $W^*[\ell]$ takes the form $\big\vert \paradelta\weight_{n;(x+ z,0)}^{(y+ U,1)} \big\vert$ or  $\big\vert \paradelta\weight_{n;(x+ U,0)}^{(y+ z,1)} \big\vert$
where  $z \in [0,2^{1-i_0})$ is a 
dyadic rational of scale $D^*[\ell]$ at least $i_0 - 1 + k$ and $U \subset [0,2^{1-i_0})$ is a dyadic interval whose scale is at least that of $z$. 
For this reason, it is tempting -- though, as we shortly explain, mistaken -- to use Proposition~\ref{p.dyadic}
with ${\bf k_0} = i_0 - 1$ and ${\bf i} = D^*[\ell]$ to find an upper bound on the probability 
$$
 \PP \Big( W^*[\ell] \geq K_0 \, 2^{-D^*[\ell]/2}  \, , \,  \maxmin_{n;([x-2,x+3],0)}^{([y-2,y+3],1)}\big( 2^{-2} K_0 \big) \Big) \, ,
$$
where recall from the proposition that $K_0 \geq 9$ is a parameter.
Indeed, the proposition formally implies that, 
whenever $k \in \N$ and $\ell \in \{ 2k,2k+1 \}$,
this probability is at most
$$
 2 \cdot 
 2^{2 (D^*[\ell] - i_0 +1  )} \cdot   \exp \big\{ - 2^{-3}  \Kzero^2 \big\}  \, .
$$
This application is erroneous, however, because  $D^*[\ell]$ is random and thus the choice ${\bf i} = D^*[\ell]$ is inadmissible. 
We adjust to cope with this problem by introducing a parameter $j \in \N$ that is supposed to be at least $i_0 - 1 +k$. For any given such $j$, we may now apply
 Proposition~\ref{p.dyadic}
with ${\bf k_0} = i_0 - 1$ and ${\bf i} = j$. We thus find that, for $k \in \N$, $\ell \in \{ 2k,2k+1 \}$ and $j \geq i_0 - 1 +k$, 
\begin{eqnarray*}
 & &\PP \Big( W^*[\ell] \geq K_0 \, 2^{-j/2}  \, , \, D^*[\ell] = j  \, , \, \maxmin_{n;([x-2,x+3],0)}^{([y-2,y+3],1)}\big( 2^{-2} K_0 \big) \Big) \\
 & \leq &  
 2^{1 + 2 (j - i_0 +1  )}    \exp \big\{ - 2^{-3}  \Kzero^2 \big\}  \, . 
\end{eqnarray*}

With such $j$ remaining fixed,
we now choose the parameter $K_0$ as a function of $\ell \in \N$ and $j$,
setting $K_0 = S \cdot 2^{-i_0/2 - k/4} \cdot 2^{j/2}$. 
(Since $k = \lfloor \ell/2 \rfloor$, this choice is indeed determined by $\ell$.) For reasons soon to be explained, we set the new quantity $S$ equal to $2^{9/2}$.
  Using $j \geq i_0 - 1 + k$, we see that $K_0 \geq  2^{k/4 - 1/2} S$ (and thus $K_0 \geq 9$, as is needed).
Since the event $\maxmin_{n;([x-2,x+3],0)}^{([y-2,y+3],1)}(K_0/4)$ is increasing in $K_0 > 0$, 
we~find~that
\begin{eqnarray*}
 & & \PP \Big( W^*[\ell] \geq  S \cdot 2^{- i_0/2 - k/4}   \, , \, D^*[\ell] = j  \, , \,  \maxmin_{n;([x-2,x+3],0)}^{([y-2,y+3],1)}\big( 2^{- 2 - 1/4} S \big) \Big) \\
 & \leq &  
 2^{1 + 2 (j - i_0 +1  )}  \exp \big\{ - 2^{-3 -i_0 - k/2 + j} S^2 \big\}  \, .
\end{eqnarray*}
Summing over $k \geq 1$, the two values of $\ell$ for each $k$, and the values $j \in \N$ that are at least $i_0 - 1 +k$, we learn that
\begin{eqnarray}
 & &\PP \Big( \sum_{\ell = 0}^\infty W^*[\ell] \geq  2S \sum_{k=0}^\infty 2^{-i_0/2 - k/4}  \, , \,  \maxmin_{n;([x-2,x+3],0)}^{([y-2,y+3],1)}\big(  2^{- 9/4} S    \big) \Big) \nonumber  \\
 & \leq & 2 \, \sum_{\ell = 0}^\infty \sum_{j = i_0 - 1 + \lfloor \ell/2 \rfloor}^\infty 
  2^{1 + 2(j - i_0 + 1)}  
 \exp \big\{  - 2^{-3 -i_0 - \ell/4 + j} S^2  \big\} \nonumber \\
 & \leq & 8 \, \sum_{\ell = 0}^\infty  
 2^{2 \lfloor \ell/2 \rfloor}
 \exp \big\{  - 2^{-3 -i_0 - \ell/4 + i_0 - 1 + \lfloor \ell/2 \rfloor} S^2 
     \big\} \nonumber \\
     & \leq & 8 \, \sum_{\ell = 0}^\infty  
 2^\ell
 \exp \big\{  - 2^{-5 + \ell/4} S^2 \big\} \label{e.sumexpress} \, .
 \end{eqnarray}
The final inequality is trivial, but a brief explanation is needed to verify the second.  Note that the ratio of the summands on this inequality's left-hand side  indexed by $(\ell,j+1)$ and $(\ell,j)$ is at most one-half provided that 
 $2^{-3 - i_0 - \ell/4 + j} S^2 \geq 3 \log 2$; since $j \geq i_0 - 1 + k$, while $k = \lfloor \ell/2 \rfloor$ and $\ell \geq 0$, the latter bound is seen to be valid when  $S \geq 2^{5/2} 3^{1/2} (\log 2)^{1/2}$. Since $S$  equals $2^{9/2}$, the last condition holds; and thus the bound in question emerges because its right-hand side is the sum of the concerned geometric series.

Consider then the expression $2^\ell
 \exp \big\{  - 2^{-5 + \ell/4} S^2 \big\}$.  
The ratio of each  summand, indexed by $\ell \geq 1$, to its predecessor is at most $2 \exp \big\{ - S^2 \cdot 2^{-5} ( 2^{1/4} - 1 ) \big\}$
which when, as is presently supposed, $S \geq 2^{5/2 + 2} = 2^{9/2}$ is at most $2 \exp \big\{ - 2^{2 \cdot 2} ( 2^{1/4} - 1 ) \big\}
\leq 3/4$; so that~(\ref{e.sumexpress}) is at most
$$
 4  \cdot 8 \, \exp \big\{  - 2^{-5} S^2 
     \big\} =  32 \,
 \exp \big\{ - 2^{-5} S^2 \big\} \, .
$$
By~(\ref{e.sumwk}) and the definition of the sequence $\big\{ W^*[\ell]: \ell \geq \N , \ell \geq 1 \big\}$, 
$$
  \sum_{\ell = 1}^\infty W^*[\ell] \geq  \sup_{(u,v) \in [x, x+\e] \times [y,y+\e]} \Big\vert \paradelta\weight_{n;([x,u],0)}^{([y,v],1)} \Big\vert \, .
$$
Thus,
\begin{eqnarray*}
 & &\PP  \, \bigg(  \, \sup_{(u,v) \in [x, x+\e] \times [y,y+\e]} \Big\vert \paradelta\weight_{n;([x,u],0)}^{([y,v],1)} \Big\vert \, \geq \, 2^{1-i_0/2}  (1-2^{-1/4})^{-1}S   \, , \,  \\
 & & \qquad  \qquad  \qquad  \qquad     \maxmin_{n;([x-2,x+3],0)}^{([y-2,y+3],1)}\big(  2^{-9/4} S   \big) \, \bigg) \,  \leq  \, 
 32 \,
 \exp \big\{ - 2^{-5} S^2 \big\}  \, .
  \end{eqnarray*}
We find that
\begin{eqnarray*}
 & &
  \PP \, \bigg( \, \sup 
 \Big\vert \paradelta\weight_{n;([u_1,u_2],0)}^{([v_1,v_2],1)} \Big\vert  
  \, \geq  \,  2 \cdot
  2^{1-i_0/2}  (1-2^{-1/4})^{-1}S        \, , \,  \\
 & & \qquad  \qquad  \qquad  \qquad     \maxmin_{n;([x-2,x+3],0)}^{([y-2,y+3],1)}\big(   2^{- 9/4} S    \big) \, \bigg) \,  \leq  \,    
   32 \,
 \exp \big\{ - 2^{-5} S^2 \big\}   \, ,
\end{eqnarray*}
where the supremum is over arbitrary $u_1,u_2 \in [x,x+\e]$ and $v_1,v_2 \in [y,y+\e]$.

Applying~(\ref{e.maxminbound}), and recalling that $\e \geq 2^{-i_0}$,
\begin{eqnarray*}
 & & 
  \PP \, \Big( \, \sup 
 \Big\vert \paradelta\weight_{n;([u_1,u_2],0)}^{([v_1,v_2],1)} \Big\vert  
  \, \geq \, \e^{1/2}
  2^2  (1-2^{-1/4})^{-1}S   
      \Big) \\
 & \leq &
  \PP \, \Big( \, \sup 
 \Big\vert \paradelta\weight_{n;([u_1,u_2],0)}^{([v_1,v_2],1)} \Big\vert  
  \, \geq \,  
  2^{2-i_0/2}  (1-2^{-1/4})^{-1}S   
      \Big) \\
  &  \leq  &   
   32 \,
 \exp \big\{ - 2^{-5} S^2 \big\}   
   \, + \, 
5^2 \cdot  400 \, C \exp \big\{ - c_1 2^{-10} r_0^{3/2} \big\} \, ,
\end{eqnarray*}
where $r_0 =    2^{- 9/4} S$. 
The right-hand side is at most
\begin{eqnarray*}
  & &  32 \,
 \exp \big\{ - 2^{-5} S^2 \big\}    \, + \, 
10000  \,  C \exp \big\{ -  c_1 2^{-10 - 27/8}   S^{3/2}  \big\} 
 \\
&  \leq & 10032 \, C  \exp \big\{ - c_1 2^{-14}   S^{3/2}   \big\} \, ,
\end{eqnarray*}
the displayed bound due to $S \geq 1$, $C \geq 1$ and $c_1 \leq 1$. 
Setting $R =  2^2(1-2^{-1/4})^{-1}S$ and noting that $2^{-14-3}   (1-2^{-1/4})^{3/2} \geq 2^{-21}$ 
completes the proof of Theorem~\ref{t.differenceweight}. \qed

\section{Profile regularity for general initial data: proving Theorems~\ref{t.nmodcon} and~\ref{t.wlp.one}}

For the proofs of these two theorems, it is needed  that, typically, every $f$-rewarded polymer 
that ends in the interval $[-1,1] \times \{ 1 \}$ begins in a compact subset of the $x$-axis $\R \times \{ 0 \}$.
We first introduce a suitable regular fluctuation event $\regfluc$, and then derive a result to this effect, Lemma~\ref{l.regfluc}.

\begin{definition}\label{d.regfluc}
Recall Definition~\ref{d.if} and the notational usage
of  $\rho_{n;(*:f,0)}^{(y,1)}$  from Subsection~\ref{s.polymercross}. Let  $\ovbar\coninit \in (0,\infty)^3$ and $f \in \initcond_{\ovbar\coninit}$.
For $R \geq 0$, define the event 
$\regfluc_{n;(*:f,0)}^{\big( \{-1,1 \}, 1 \big)}(R)$
that 
any $f$-rewarded line-to-point polymer $\rho_{n;(*:f,0)}^{(-1,1)}$ that ends at $(-1,1)$ begins at a location $(x,0)$ where $x \geq - (R+1)$; and any such polymer $\rho_{n;(*:f,0)}^{(1,1)}$ that ends at $(1,1)$ begins at $(x,0)$, where $x \leq R+1$.
\end{definition}


We remark that our $\regfluc$ event entails that any $f$-rewarded line-to-point polymer that ends at a location $(y,1)$ with $y \in [-1,1]$ must begin at a location $(x,0)$, where $\vert x \vert  \leq R+1$. Indeed, were such a polymer $\rho_{n;(*:f,0)}^{(y,1)}$ to begin at $(x,0)$, with $x < -(R+1)$, then, in the event $\regfluc$, it would cross  any example of $\rho_{n;(*:f,0)}^{(-1,1)}$.
The rewiring of these two polymers described in  Subsection~\ref{s.polymercross} would then furnish an example of  $\rho_{n;(x,0)}^{(-1,1)}$ with $x < -(R+1)$, 
in conflict with the occurrence of $\regfluc$.
 
\begin{lemma}\label{l.regfluc}
Let $n \in \N$, $R > 0$ and $\ovbar\coninit \in (0,\infty)^3$ satisfy 
$$
  n \geq  c^{-18} \max \Big\{  (\coninit_2 + 1)^9 \, , \,   10^{23} \coninit_1^9     \, , \, 3^{9}  \Big\} \, ,  \, R  \geq \max \Big\{ \,    39 \coninit_1  \,  , \, 5  \, , \,   3 c^{-3}  \, , \, 2 \big( (\coninit_2 + 1 )^2 +  \coninit_3 \big)^{1/2} \, \Big\} \, , 
$$
  and
 $R \leq 6^{-1}c n^{1/9}$.
Then, for any $f \in \initcond_{\ovbar\coninit}$,
$$
 \PP \Big( \neg \, 
\regfluc_{n;(*:f,0)}^{\big( \{-1,1 \}, 1 \big)}(R)  \Big) \leq 38 R \rsC       \exp \big\{ -  2^{-6} \rsc   R^3 \big( 2^{ - 1/2} - 2^{-1} \big)^{3/2}  \big\}  \, .
$$
\end{lemma}
\noindent{\bf Proof.} The event  $\neg \, 
\regfluc_{n;(*:f,0)}^{\big( \{-1,1 \}, 1 \big)}(R)$ equals $A_1 \cup A_2$, where 
 $A_1$ is the event that $y \to \mc{L}^{\downarrow;(-1,1)}_{n;0}(1,y) + f(y)$ achieves its maximum for a value of $y$ that is less than  $-1 - R$,
and  $A_2$ is the event that $y \to \mc{L}^{\downarrow;(1,1)}_{n;0}(1,y) + f(y)$ achieves its maximum for a value of $y$ that is at greater than $1 + R$.

Note the inclusion
\begin{equation}\label{e.twoleft}
\Big\{ \sup_{x \in [-\coninit_2,\coninit_2]}  \big( \mc{L}^{\downarrow;(-1,1)}_{n;0}(1,x) + f(x)  \big) > - R^2/2 \, , \, \sup_{x \leq - 1 -R} \big( \mc{L}^{\downarrow;(-1,1)}_{n;0}(1,x) + f(x) \big)  \leq - R^2/2  \Big\} \, \subseteq \, A_1^c \, .
\end{equation}
Upper bounds on the failure probability of the left-hand events will now be found.
The first event will be shown to be probable because, in view of Definition~\ref{d.if}, the function $f$ is known to assume a not highly negative value somewhere in a compact interval about the origin.
The second event is probable due to the at most linear growth of $f$ far from the origin, combined with decay estimates on the curve $\mc{L}^{\downarrow;(-1,1)}_{n;0}(1,x)$ for large $x$.
These estimates take two forms: when $x$ is large, but less than order $n^{1/9}$, the curve hews to a parabola, in accordance with the No Big Max Proposition~\ref{p.mega}(2), applied to the normalized cousin of the ensemble in question;
when $x$ becomes even larger, the curve may escape the reaches of this parabola, but it continues to decay rapidly, in accordance with collapse near infinity Proposition~\ref{p.mega}(3). 


Since $f \in \initcond_{\ovbar\coninit}$, there exists $x_0 \in [-\coninit_2,\coninit_2]$ such that $f(x_0) \geq - \coninit_3$.
As such, the first left-hand event in~(\ref{e.twoleft}) fails with a probability that satisfies
\begin{eqnarray*}
 & & \PP \Big( \sup_{x \in [-\coninit_2,\coninit_2]}  \big( \mc{L}^{\downarrow;(-1,1)}_{n;0}(1,x) + f(x) \big) \leq - R^2/2   \Big) \leq 
\PP \Big(  \mc{L}^{\downarrow;(-1,1)}_{n;0}(1,x_0)  \leq - R^2/2 + \coninit_3   \Big) \\
 & = & \PP \Big(  \scaledle^{\downarrow;(-1,1)}_{n;0}(1,x_0 + 1)  \leq - R^2/2 + \coninit_3   \Big) \\
& \leq & 
\PP \Big(   \scaledle^{\downarrow;(-1,1)}_{n;0}(1,x_0 + 1)  + 2^{-1/2} (x_0 + 1 )^2  \leq - R^2/2 + 2^{-1/2} (\coninit_2 + 1)^2 +  \coninit_3   \Big) \\
 & \leq & 
\PP \Big(    \scaledle^{\downarrow;(-1,1)}_{n;0}(1,x_0 + 1)  + 2^{-1/2} (x_0 + 1 )^2    \leq - R^2/4   \Big) \leq   \rsC  \exp \big\{ - 2^{-3} \rsc R^3 \big\} \, , 
\end{eqnarray*}
where the penultimate inequality depends on $R^2/4 \geq  2^{-1/2} (\coninit_2 + 1 )^2 +  \coninit_3$.
The final inequality  was  obtained
by applying the one-point lower tail $\rmreg(2)$ with parameter choices ${\bf z} = x_0$
and ${\bf s} = R^2/4$
to the $(c,C)$-regular ensemble  $\scaledle_{n;0}^{\downarrow;(-1,1)}$.
Since the ensemble has $n+1$ curves, this application of $\rmreg(2)$  may be made provided that $\vert x_0 \vert + 1 \leq \rsc n^{1/9}$ and $R^2/4 \in [1,n^{1/3}]$.
The first of these conditions due to $n \geq c^{-9} (\coninit_2 + 1)^9$ alongside 
$\vert x_0 \vert \leq \coninit_2$; the second we assume.

The failure probability of the second left-hand event in~(\ref{e.twoleft}) may be gauged as follows: since $f(x)$ is at most $\coninit_1 \big( 1 + \vert x \vert \big)$ for any $x \in \R$, we may note that
\begin{eqnarray}
 & &  \PP \Big( \sup_{x \leq - 1 -R} \big( \mc{L}^{\downarrow;(-1,1)}_{n;0}(1,x) + f(x) \big) > - R^2/2 \Big)  \nonumber \\
 & = & \PP \Big( \sup_{x \leq - R} \big( \scaledle^{\downarrow;(-1,1)}_{n;0}(1,x) + f(x-1) \big) > - R^2/2 \Big)  \nonumber \\
  & \leq & \PP \Big( \sup_{x \leq - R} \big(  \scaledle^{\downarrow;(-1,1)}_{n;0}(1,x) + \coninit_1 ( 2 + \vert x \vert ) \big) > - R^2/2 \Big) \, ; \label{e.aboutsplit}
\end{eqnarray}
the latter term may then be bounded above by
$$
 \sum  \PP \Big( \sup_{x \in -  R [2^j,2^{j+1}] }  \scaledle^{\downarrow;(-1,1)}_{n;0}(1,x) > - R^2/2 -   \coninit_1 \big( 2 +    2^{j+1} R  \big) \Big) \, \, + \, E_1 + E_2 \, ,
$$
 where the first sum is indexed by a parameter $j$ that varies over the integer interval $\llbracket 0,k \rrbracket$ where $k \in \N$ to chosen to be  maximal subject to 
 $2^{k+1} R \leq  3^{-1} c  n^{1/9}$. (Such a $k$ exists because we suppose that $2R \leq  3^{-1} c  n^{1/9}$.)
 The term $E_1$ corresponds to part of a dyadic scale that has been sliced in two by the value  $-3^{-1} c  n^{1/9}$: this term is specified by the expression in~(\ref{e.aboutsplit}) when the supremum in the variable $x$ is chosen to be  over  the interval $\big[ - 3^{-1} \rsc  n^{1/9} , - 2^{k+1}R \big]$. The remaining term $E_2$ is a long-range error term corresponding to the interval $\big[ - \xnmac ,  - 3^{-1} \rsc  n^{1/9} \big]$.   Since $R \leq 3^{-1} \rsc  n^{1/9}$, 
 $$
 E_2  \leq  \PP \Big( \sup \Big\{   \scaledle^{\downarrow;(-1,1)}_{n;0}(1,x) + \coninit_1 ( 2 + \vert x \vert ) 
  : x \in \big[ - \xnmac ,  - 3^{-1} \rsc  n^{1/9} \big] \Big\} \geq - 2^{-1} (c/3)^2 n^{2/9}  \Big) \, .
 $$ 
This right-hand side will be bounded above by applying collapse-near-infinity Proposition~\ref{p.mega}(3)
to the $(n+1)$-curve $(c,C)$-regular ensemble  $\scaledle_{n;0}^{\downarrow;(-1,1)}$.
We apply  Proposition~\ref{p.mega}(3)
 with its parameter ${\bm \eta}$ chosen so that  ${\bm \eta} (n+1)^{1/9} = -  3^{-1} \rsc  n^{1/9}$. In order to make the application, we first claim that the affine function $x \to \ell(x)$
in 
the proposition
lies below the function 
\begin{equation}\label{e.affinefunction}
x \to - 2^{-1} (c/3)^2 n^{2/9} -   \coninit_1 \big( 2 + \vert x \vert    \big) 
\end{equation}  
whenever $x \leq -  3^{-1} \rsc  n^{1/9}$.
 To verify this, note that, when $x = - 3^{-1} \rsc  n^{1/9}$, the assertion  takes the form $\coninit_1 \big( 2 + 3^{-1} \rsc  n^{1/9} \big) \leq \big( 2^{-1/2} - 2^{-5/2} - 2^{-1} \big) (\rsc/3)^2 n^{2/9}$, which holds due to the supposed $1 \leq 3^{-1} \rsc  n^{1/9}$ and $\coninit_1 \leq \big( 2^{-1/2} - 2^{-5/2} - 2^{-1} \big) (\rsc/9)n^{1/9}$. Confirming the claim is then a matter of checking that the gradient of $\ell$ exceeds that of the function~(\ref{e.affinefunction}), which holds due to $\coninit_1 \leq 5 \cdot 2^{3/2} c/3 \cdot n^{1/9}$.
 
We may thus apply  Proposition~\ref{p.mega}(3)
   when $n + 1  \geq 2^{45/4} \rsc^{-9}$, finding that 
$$
 E_2 \leq 
6C \exp \big\{ -  2^{-15/4} 3^{-3}  \rsc^4 n^{1/3} \big\} \, .
$$ 
 Note that
 \begin{eqnarray}
  & & 
    \PP \Big( \sup_{x \in -  R [2^j,2^{j+1}] }  \scaledle^{\downarrow;(-1,1)}_{n;0}(1,x) > - R^2/2 -   \coninit_1 \big( 2 +    2^{j+1} R  \big) \Big) \label{e.nrphi} \\
  & \leq & 
    \PP \Big( \sup_{x \in -  R [2^j,2^{j+1}] } \big(  \scaledle^{\downarrow;(-1,1)}_{n;0}(1,x)  + 2^{-1/2} x^2 \big)  >  R^2 \big(   2^{2j -1/2}   - 1/2   \big)  -   \coninit_1 \big( 2 +    2^{j+1} R  \big) \Big) \nonumber \\  
  & \leq & 
    \PP \Big( \sup_{x \in -  R [2^j,2^{j+1}] } \big(  \scaledle^{\downarrow;(-1,1)}_{n;0}(1,x)  + 2^{-1/2} x^2 \big)  >  2^{-1} R^2 \big(   2^{2j -1/2}   - 1/2   \big)   \Big) \nonumber \\
       & \leq &  6 C (   2^{j-1} R +1 )     \exp \big\{ -  2^{-6} \rsc   R^3 \big( 2^{2j - 1/2} -  2^{-1} \big)^{3/2}  \big\} \, , \nonumber
\end{eqnarray}
where in the second inequality we used $R \geq 1 \vee 39 \coninit_1$ in the form 
$$
2^{-1} R^2 \big( 2^{2j - 1/2} - 1/2 \big) \geq \coninit_1 \big( 2 + 2^{j+1}R  \big) \, \, \, \,  
\textrm{for each $j \geq 0$} \, .
$$ 
The final inequality arises from an application of 
Proposition~\ref{p.mega}(2) 
to the ensemble
$\scaledle_{n;0}^{\downarrow;(-1,1)}$. 
The parameters of the application are set to be 
$$
{\bf y} = - 2^{-1} R \big( 2^j + 2^{j+1} \big) \, , \, {\bf r} = 2^{-1} R \big( 2^{j+1} - 2^j \big) \, \,  \textrm{and} \, \, {\bf t} =  2^{-1} R^2 \big( 2^{2j - 1/2} - 2^{-1} \big) \, .
$$
The application's hypotheses are implied by  
$$
3 \cdot 2^j R  \leq c n^{1/9} \, , \, 
 2^{j+1} R  \leq \rsc n^{1/9} \, , \, 2^{-1} R^2 \big( 2^{2j - 1/2} - 1/2 \big) \in \big[ 2^{7/2} , 2 n^{1/3} \big] \, \textrm{ and }  \, n \geq c^{-18} \, .
$$
  The first three of these conditions are valid when $j \in \llbracket 0,k \rrbracket$ in light of the assumed bound $2^{k+1} R \leq 3^{-1} c n^{1/9}$ (and $c \leq 1$); indeed, they are also valid when $j = k+1$, a fact that we will use momentarily.

 The term $E_1$ is bounded above by~(\ref{e.nrphi}) with $j = k+1$, so that the preceding argument shows that
$$
 E_1 \leq 
 6 C (   2^k R +1 )    \exp \big\{ -  2^{-6} \rsc   R^3 \big( 2^{2k + 3/2} -  2^{-1} \big)^{3/2}  \big\}  \, .
$$  
Thus,
 \begin{eqnarray*}
 & & \PP \Big( \sup_{x \leq - 1 -R} \big( \mc{L}^{\downarrow;(-1,1)}_{n;0}(1,x) + f(x) \big) > - R^2/2 \Big) 
  \\
 & \leq & 6  \rsC \sum_{j=0}^{k+1}  (   2^j R +1 )       \exp \big\{ -  2^{-6} \rsc   R^3 \big( 2^{2j - 1/2} - 2^{-1} \big)^{3/2}  \big\}
\, + \, 
6C \exp \big\{ -  2^{-15/4} 3^{-3}  \rsc^4 n^{1/3} \big\} \\
 & \leq & 12 R \rsC       \exp \big\{ -  2^{-6} \rsc   R^3 \big( 2^{ - 1/2} - 2^{-1} \big)^{3/2}  \big\}
\, + \, 
6C \exp \big\{ -  2^{-15/4} 3^{-3}  \rsc^4 n^{1/3} \big\} \\
 & \leq & 18 R \rsC       \exp \big\{ -  2^{-6} \rsc   R^3 \big( 2^{ - 1/2} -  2^{-1} \big)^{3/2}  \big\} \, ,
 \end{eqnarray*}  
 where the second inequality used $R \geq (\log 4)^{1/3} 2^2 c^{-3} \big( (2^{3/2} - 2^{-1})^{3/2} - (2^{-1/2} - 2^{-1})^{3/2} \big)^{-1/3}$ in order to ensure that each ratio of consecutive summands in the sum is at most one-half; the third makes use of $1 \leq R \leq  (2^{-1/2} - 2^{-1})^{-1/2} 2^{3/4} 3^{-1} \rsc n^{1/9}$.

 Thus,
\begin{eqnarray*}
 \PP \big( A_1 \big) & \leq & 
   \rsC  \exp \big\{ - 2^{-3} \rsc R^3 \big\} 
 \, + \, 18 R \rsC       \exp \big\{ -  2^{-6} \rsc   R^3 \big( 2^{ - 1/2} - 2^{-1} \big)^{3/2}  \big\} \\
  & \leq &  19 R \rsC       \exp \big\{ -  2^{-6} \rsc   R^3 \big( 2^{ - 1/2} - 2^{-1} \big)^{3/2}  \big\} \, ,
\end{eqnarray*}
the latter inequality due to $R \geq 1$.
The same argument yields that $\PP \big( A_2 \big)$ satisfies the same upper bound. Combining the estimates completes the proof of Lemma~\ref{l.regfluc}. \qed

Theorem~\ref{t.nmodcon}'s proof will also harness estimates showing that a {\em local weight regularity} event is typical. 
For $x,y \in \R$ and $\e,r > 0$, this event is defined by
$$
\lwr_{n;([x,x+\e],0)}^{([y,y+\e],1)}\big( \e , r  \big) \, = \,
 \left\{ \sup_{\begin{subarray}{c} x_1,x_2 \in [x,x+\e] \\
    y_1, y_2 \in [y,y+\e]  \end{subarray}} \Big\vert \weight_{n;(x_2,0)}^{(y_2,1)}  - \weight_{n;(x_1,0)}^{(y_1,1)}  \Big\vert \, \leq \, r \e^{1/2} \right\} \, . 
$$

The relevant control is offered by Theorem~\ref{t.differenceweight}, except that this theorem addresses parabolically adjusted weight.
The next result is the one we will apply in proving 
Theorem~\ref{t.wlp.one}: the main new hypothesis, $\vert x - y \vert \leq \e^{-1/2}$, limits parabolic curvature.  
\begin{corollary}\label{c.ordweight}
Let 
$n \in \N$ and  $x,y \in \R$ satisfy 
$n \geq 10^{32} c^{-18}$ and   $\big\vert x - y  \big\vert \leq \e^{-1/2} \wedge 2^{-2} 3^{-1} \rsc  n^{1/18}$.
Let  $\e \in (0,2^{-4}]$ and
 $R \in \big[2 \cdot 10^4 \, , \,   10^3 n^{1/18} \big]$.
Then
\begin{equation}\label{e.ordweight}
 \PP \Big( \neg \, \lwr_{n;([x,x+\e],0)}^{([y,y+\e],1)}\big( \e , R  \big)    \Big) \leq 
 10032 \, C  \exp \big\{ - c_1 2^{-22 - 1/2} R^{3/2}     \big\}
 \, .
\end{equation}
\end{corollary}
\noindent{\bf Proof.} Recall that we denote $Q:\R \to \R$, $Q(u) = 2^{-1/2} u^2$. For $u_1,u_2 \in [x,x+\e]$
and $v_1,v_2 \in [y,y+\e]$, note that, since $\vert x - y \vert \leq \e^{-1/2}$ and $\e \leq 1$, 
$$
\big\vert Q(v_2 - u_2) - Q(v_1 - u_1) \big\vert \leq 
2\e Q' \big( \vert x - y \vert + \e \big) = 2\e \cdot 2^{1/2}  \big( \vert x - y \vert + \e \big) \leq 2^{5/2} \e^{1/2} \, .
$$
Thus, when $R \geq 2 \cdot 2^{5/2}$, the left-hand side of~(\ref{e.ordweight}) is at most~(\ref{e.differenceweight}) with ${\bf R} = R/2$. The corollary thus follows from Theorem~\ref{t.differenceweight}. \qed

\subsection{Deriving Theorem~\ref{t.nmodcon}}


For  $f \in  \initcond_{\ovbar\coninit}$, $\e  \in (0,2)$ and $\rho > 0$, define   the equicontinuity event  
 $$
 \mathsf{EquiCty}_{n;(*:f,0)}^{[-1,1]}\big( \rho , \e \big)
 = \Big\{ \, \omega_{[-1,1],\e}\big( y \to  \weight_{n;(*:f,0)}^{(y,1)} \big) <  \rho \, \Big\} \, ,
$$ 
where the modulus of continuity of a function $h: [-1,1] \to \R$ is denoted by
$$
    \omega_{[-1,1],\e}(h) \,  = \, \sup \Big\{ \, \big\vert h(x) - h(y) \big\vert : x,y \in [-1,1] \, , \,  \vert x - y \vert \leq \e \, \Big\} \, .
$$
Here is an equicontinuity claim.
\begin{lemma}\label{l.equicty}
Let $\ovbar\coninit \in (0,\infty)^3$. Set
$$
\Ctbs = 2^{17} c_1^{-2/3} \, \, \, \, \textrm{and} \, \, \, \,  \Ctbd = \max \Big\{ \,    39 \coninit_1  \,  , \,  12 c^{-3}  \, , \, 2 \big( (\coninit_2 + 1 )^2 +  \coninit_3 \big)^{1/2} \, \Big\} \, ,
$$
and define the function  $g(\e) = 2 \Ctbs \, \e^{1/2} \big( \log \e^{-1} \big)^{2/3}$. Let $\e > 0$ satisfy 
$\e \leq \Ctbd^{-6} \exp \big\{ - 2^{3/2} 10^6 \big\} C^{-2}$.
When  $n \in \N$ verifies
$n \geq 10^{32} c^{-18} \Ctbd^{18} \Ctbs^{18} \big( \log \e^{-1} \big)^{12}$, we have that, for   $f \in \initcond_{\ovbar\coninit}$, 
\begin{equation}\label{e.ece}
 \PP \Big( \neg \,   \mathsf{EquiCty}_{n;(*:f,0)}^{[-1,1]}\big( g(\e) , \e \big)  \Big) \leq  \e \, .
\end{equation}
\end{lemma}
\noindent{\bf Proof.} We first argue that, whenever $R \geq 1$,  
\begin{eqnarray}
 & & \regfluc_{n;(*:f,0)}^{\big( \{-1,1 \}, 1 \big)}(R-1) \cap \bigcap_{u \in \e \Z \cap [-R,R], 
 v \in \e \Z \cap [-1,1]} \lwr_{n;([u,u+\e],0)}^{([v,v+\e],1)} \Big( \e , \Ctbs \big( \log \e^{-1} \big)^{2/3}  \Big)  \nonumber \\
 & \subseteq & \mathsf{EquiCty}_{n;(*:f,0)}^{[-1,1]}\big( g(\e) , \e \big)  \, .   \label{e.equictyinc}
\end{eqnarray}
 To verify this inclusion, suppose that $\regfluc_{n;(*:f,0)}^{\big( \{-1,1 \}, 1 \big)}(R-1)$ occurs,
and consider $y \in [-1,1]$. By the remark made 
after Definition~\ref{d.regfluc}, all $f$-rewarded line-to-point polymers that abut at time one on $[-1,1]$
must begin at time zero somewhere on $[-R,R]$, (where here of course the present parameter value $R-1$ is involved). Thus, 
the quantity
$\weight_{n;(*:f,0)}^{(y,1)}$ is seen to equal $\weight_{n;(x,0)}^{(y,1)} + f(x)$ for some $x \in [-R,R]$.
Note that the event $\lwr_{n;(x,0)}^{([y,y+\e],1)} \Big( \e , 2\Ctbs \big( \log \e^{-1} \big)^{2/3}  \Big)$ occurs when the intersection of the $\lwr$ events displayed above occurs;
in this circumstance, we thus see that, for any $\eta \in (0,\e)$, 
$\weight_{n;(x,0)}^{(y + \eta,1)} \geq  \weight_{n;(x,0)}^{(y,1)}  + f(x) \, - \,  \e^{1/2} \cdot 2\Ctbs \big( \log \e^{-1} \big)^{2/3}$
and thus 
$\weight_{n;(*:f,0)}^{(y + \eta,1)} \geq  \weight_{n;(*:f,0)}^{(y,1)} \, - \,  \e^{1/2} \cdot 2\Ctbs \big( \log \e^{-1} \big)^{2/3}$.  Provided that we further suppose that $y + \eta \leq 1 $, the inequality with the roles of $y$ and $y+\eta$ reversed is similarly obtained, so that
$$
 \Big\vert \weight_{n;(*:f,0)}^{(y + \eta,1)} -  \weight_{n;(*:f,0)}^{(y,1)}   \Big\vert \, \leq \,  \e^{1/2} \cdot 2\Ctbs \big( \log \e^{-1} \big)^{2/3} \, .
$$
Thus,~(\ref{e.equictyinc}) is obtained. Verifying the equicontinuity claim is now a matter of arguing that the $\regfluc$ and the intersection of the $\lwr$
events on the left-hand side of~(\ref{e.equictyinc}) both have probability at least $1 - \e/2$. 

Treating the intersection $\lwr$ event first, we now set the value of $R$ equal to $1 + \Ctbd \big( \log \e^{-1} \big)^{1/3}$ where $\Ctbd$ is a further positive parameter on which we will impose certain lower bounds.

Let $u,v \in \R$ satisfying $\vert u \vert \leq R$ and $\vert v \vert \leq 1$ be given.  We apply Corollary~\ref{c.ordweight} with ${\bf x} = u$, ${\bf y} = v$, ${\bm \e} = \e$ and  ${\bf R} = \Ctbs \big( \log \e^{-1} \big)^{2/3}$ to find that
$$
 \PP \Big( \neg \,  \lwr_{n;([u,u+\e],0)}^{([v,v+\e],1)} \Big( \e , \Ctbs \big( \log \e^{-1} \big)^{2/3}  \Big) \leq
    10032 \, C  \e^{c_1 2^{-22 - 1/2} \Ctbs^{3/2}    } \, .
 $$
 Since $\vert u \vert \leq R = 1 + \Ctbd \big( \log \e^{-1} \big)^{1/3}$ and $\vert v \vert \leq 1$, this application may be made provided that $\e \in (0,2^{-4}]$,  
$n \geq 10^{32} c^{-18}$,  $\Ctbd \big( \log \e^{-1} \big)^{1/3} + 2  \leq \e^{-1/2} \wedge 2^{-2} 3^{-1} \rsc  n^{1/18}$, and 
$\Ctbs \big( \log \e^{-1} \big)^{2/3} \in \big[ 2 \cdot 10^4 \, , \,   10^3 n^{1/18} \big]$. Thus, it may be made for $\e > 0$ sufficiently small, and with $n$ exceeding an $\e$-determined level whose order is $\big( \log \e^{-1} \big)^{12}$. 
 
 Allowing $u$ and $v$ to vary over $\e \Z \cap [-R,R]$ and $\e \Z \cap [-1,1]$, the probability that any of the $\lwr$ events so indexed fails is seen to be at most
 \begin{equation}\label{e.atmostebytwo}
    \big( 2 R \e^{-1} + 1 \big)    \big( 2 \e^{-1} + 1  \big)  \cdot  10032 \, C  \e^{c_1 2^{-22 - 1/2} \Ctbs^{3/2}    }
 \end{equation}
 and thus at most $\e/2$ since $\Ctbs$ satisfies $c_1 2^{-22 - 1/2} \Ctbs^{3/2} - 2 > 1$, and $\e > 0$ is small enough.  
 
Lemma~\ref{l.regfluc} shows that the failure probability of the $\regfluc$ event is governed by a similar bound. Indeed, setting ${\bf R} = R - 1$ in the lemma, we see that
\begin{equation}\label{e.rhsub}
 \PP \Big( \neg \, 
\regfluc_{n;(*:f,0)}^{\big( \{-1,1 \}, 1 \big)}(R-1)  \Big) \leq 38 (R-1) \rsC       \exp \big\{ -  2^{-6} \rsc   (R - 1)^3 \big( 2^{ - 1/2} - 2^{-1} \big)^{3/2}  \big\}  \, ,
\end{equation}
provided that 
$n \geq  c^{-18} \max \big\{  (\coninit_2 + 1)^9 \, , \,   10^{23} \coninit_1^9     \, , \, 3^{9}  \big\}$,
 $$
    R  \geq 1 + \max \Big\{ \,    39 \coninit_1  \,  , \, 5  \, , \,   3 c^{-3}  \, , \, 2 \big( (\coninit_2 + 1 )^2 +  \coninit_3 \big)^{1/2} \, \Big\} \, , 
  $$
  and
 $R - 1 \leq 6^{-1}c n^{1/9}$.
 
 Recalling that $R = 1 + \Ctbd \big( \log \e^{-1} \big)^{1/3}$, we see that, in essence since $2^{-6} \rsc   \Ctbd^3 \big( 2^{ - 1/2} - 2^{-1} \big)^{3/2}$ exceeds one, and 
 for $\e > 0$ small enough, 
\begin{equation}\label{e.regflucbound} 
 \PP \Big( \neg \, 
\regfluc_{n;(*:f,0)}^{( \{-1,1 \}, 1 )}(R)  \Big) \,  \leq  \, \e/2 \, .
\end{equation}
 We infer then from~(\ref{e.equictyinc}) that Lemma~\ref{l.equicty} holds. \qed

\noindent{\bf Proof of Theorem~\ref{t.nmodcon}.}
Set $X_n:[-1,1] \to \R$, $X_n(y) = \weight_{n;(*:f,0)}^{(y,1)}$, and let $g:(0,1) \to (0,\infty)$ be specified via Lemma~\ref{l.equicty}. Setting $c' = 10^{-29/12} c^{-3/2} \Ctbd^{3/2} \Ctbs^{3/2}$,
this lemma may be used to show that 
\begin{equation}\label{e.omegarho}
 \PP \Big(  \omega_{[-1,1],\rho}\big( y \to  X_n(y) \big) \leq g(2\rho)    \, \, \,  \forall \, \rho \in \big( e^{-c' n^{1/12}} ,  2^{-j}  \big)  \Big) \geq 1 - 2^{1-j} 
\end{equation}
whenever $j \in \N$ satisfies $2^{-j} \leq \e_0$ where $\e_0 = \Ctbd^{-6} \exp \big\{ - 2^{3/2} 10^6 \big\} C^{-2}$ (which is the upper bound on $\e$ in Lemma~\ref{l.equicty}). Indeed, this bound is obtained by noting that  $g$ is increasing on the interval $\big( 0 , e^{-4/3} \big)$ and $\e_0 \leq e^{-4/3}$, and applying 
Lemma~\ref{l.equicty} on decreasing dyadic scales ${\bm \e} = 2^{-j}, 2^{-j- 1}, \cdots$ for as long as the lemma's hypothesis $\e \geq \exp \big\{ - c' n^{1/12} \big\}$ permits.

Define the random variable $\zeta = \zeta_n \in \big[  e^{-c' n^{1/12}} , \e_0 \big]$ to be the
maximal value on this interval such that,  for all $\rho \in \big(   e^{-c' n^{1/12}} , \zeta)$, 
$\omega_{[-1,1],\rho} \big( y \to  X_n(y) \big) \leq g(2\rho)$; if no such value exists, set $\zeta = 0$. 
We see from~(\ref{e.omegarho}) that
$\PP \big( \zeta \leq s \big) \leq 4s$ for $s \in (  e^{-c' n^{1/12}} ,\e_0]$. 

We now fix $y,z \in [-1,1]$ with $y +  2 e^{-c' n^{1/12}} < z \leq y + e^{-1}$. 
Define $K \in \N$ to be the random integer $\lceil (z-y) \zeta^{-1} \rceil$.
Setting $h = (z-y)K^{-1}$, it is readily verified that $h \in (e^{-c' n^{1/12}},\zeta]$, and so we find 
\begin{eqnarray}
 \big\vert X_n(z) - X_n(y) \big\vert & \leq & \sum_{k=0}^{K-1} \big\vert X_n\big( y + (k+1)h \big) - X_n\big(y + k h \big) \big\vert \leq K g(2h) \nonumber
  \\
  & \leq & 2^{3/2} K \Ctbs h^{1/2} \big( \log h^{-1} \big)^{2/3} \, . \label{e.khbound}
\end{eqnarray} 
In the case that $K=1$, we may now apply  $h \leq \zeta \leq e^{-4/3}$ in order to learn that
$\big\vert X_n(z) - X_n(y) \big\vert$ is at most  $2^{3/2} K \Ctbs \zeta^{1/2} \big( \log \zeta^{-1} \big)^{2/3}$.

Now suppose that $K \geq 2$. Distinctive to this case is the bound  $h \geq \zeta/2$, which follows from $(z-y) \zeta^{-1} > 1$.
Before we use this bound, note that,
since $K h = (z-y)$, the quantity~(\ref{e.khbound})  equals 
$$
 2^{3/2} \Ctbs (z-y)^{1/2} h^{-1/2} \bigg( \frac{\log h^{-1}}{\log (z-y)^{-1}} \bigg)^{2/3} (z-y)^{1/2} \big( \log (z-y)^{-1} \big)^{2/3}
$$
and thus, in view of $z-y \leq e^{-1}$ and $\zeta/2 \leq h  \leq e^{-4/3}$, may be bounded above by
\begin{eqnarray*}
 & & 2^{3/2} \Ctbs h^{-1/2} \big( \log h^{-1} \big)^{2/3} (z-y)^{1/2} \big( \log (z-y)^{-1} \big)^{2/3} \\
 & \leq & 
 2^{3/2} \Ctbs (\zeta/2)^{-1/2} \big( \log (\zeta/2)^{-1} \big)^{2/3} (z-y)^{1/2} \big( \log (z-y)^{-1} \big)^{2/3} \, .
\end{eqnarray*}
Further using $\zeta \leq 1/2$, we find that 
 $2^{3/2 + 1/2 + 2/3} \Ctbs \zeta^{-1/2} \big( \log   \zeta^{-1} \big)^{2/3} (z-y)^{1/2} \big( \log (z-y)^{-1} \big)^{2/3}$
serves as an upper bound on  
 $\big\vert X_n(z) - X_n(y) \big\vert$ in the case that $K \geq 2$.
 Whether this case applies, or rather $K=1$, we see that  
 the random variable
$$
 S_n : = \sup \bigg\{ \big\vert X_n(z) - X_n(y) \big\vert (z-y)^{-1/2} \big( \log (z-y)^{-1} \big)^{-2/3} : -1 \leq y,z \leq 1  \, , \,  2 e^{-c' n^{1/12}} < z - y \leq e^{-1} \bigg\}
$$
is bounded above by  $2^{8/3} \Ctbs \zeta^{-1/2} \big( \log \zeta^{-1} \big)^{2/3}$. 
Recalling that  $\PP \big( \zeta \leq s \big) \leq 4s$ for $s \in ( e^{-c' n^{1/12}} , \e_0]$,
we see that, for such~$s$, 
$$
\PP \Big(  S_n \geq 2^{3/2} \Ctbs  s^{-1/2} \big( \log s^{-1} \big)^{2/3} \Big) \leq 4s \, .
$$
Set $r = 2^{8/3} \Ctbs  s^{-1/2} \big( \log s^{-1} \big)^{2/3}$.
We {\em claim} that  $r \geq 2^{10/3} \Ctbs e$ implies that $s \leq d^{-2} r^{-2} \big( \log r \big)^{4/3}$. Indeed, setting
 $h = s^{-1/2}$ and $K = d r$, with $d = 2^{-10/3} \Ctbs^{-1}$, we have $K =  h ( \log h )^{2/3}$. We also have $K \geq e$, and this implies $h \leq K$, whence $h \geq K (\log K)^{-2/3}$ and also $s \leq d^{-2} r^{-2} \big( \log (dr) \big)^{4/3}$. Since $d \leq 1$, we have our claim. Furthermore, the condition that $s \leq \e_0$ is ensured when $d^{-2} r^2 (\log r)^{4/3} \leq \e_0$, and two omitted lines of working imply that $r \geq 2^{1/2} d^{-3/2} \e_0^{-3/4}$ is enough to ensure the latter.

 That is, setting $r_0 = r_0(\ovbar\coninit) = 2^{10/3} \Ctbs e \vee 2^{11/2} \Ctbs^{3/2} \Ctbd^{9/2} C^{3/2} \exp \big\{ 3 \cdot 2^{-1/2} \cdot 10^6 \big\}$, we have found that $r \geq r_0$ implies that
$$
\PP \big(  S_n \geq r  \big) \leq 2^{26/3} \Ctbs^2 r^{-2} \big( \log r \big)^{4/3} \vee 4  e^{-c' n^{1/12}}  \, .
$$
Noting that $2^{26/3} \Ctbs^2 \leq 2^{43} c_1^{-4/3}$, and $c_1 = 2^{-5/2} c \wedge 8^{-1} \geq 2^{-3}c$ in view of $c \geq 1$,
completes the proof of Theorem~\ref{t.nmodcon}. \qed


\subsection{Deriving Theorem~\ref{t.wlp.one}}

$\empty$

\noindent{\bf Proof of Theorem~\ref{t.wlp.one}(1).} The result is Lemma~\ref{l.basic}(3).

\noindent{\bf Proof of Theorem~\ref{t.wlp.one}(2).} 
Let $f \in \initcond_{\ovbar\coninit}$ be given. From~\cite[Theorem~$8.2$]{Billingsley}, the sequence of probability measures $\big\{ \nu_{n;(*:f,0)}^{([-1,1],1)} : n \in \N \big\}$  is tight if, first, the one-point distribution is tight, in the sense that for all $\eta > 0$, there exists $K > 0$ such that, 
 for all $n \in \N$, 
\begin{equation}\label{e.onepoint}
\PP \Big( \, \big\vert \weight_{n;(*:f,0)}^{(0,1)} \big\vert \leq K \, \Big) \geq 1 - \eta \, ;
\end{equation}
and, second, if, for each $\e > 0$
and
$\eta > 0$,
there exist 
$\rho > 0$ and $n_0 \in \N$
such that, for $n \geq n_0$,
\begin{equation}\label{e.equictyclaim}
\PP \Big( \mathsf{EquiCty}_{n;(*:f,0)}^{[-1,1]}\big( \rho , \e \big) \Big) \geq 1 - \eta \, .
\end{equation}
Moreover, 
if a choice of $n_0 = n_0(\e,\eta)$ such that~(\ref{e.onepoint}) and~(\ref{e.equictyclaim}) hold whenever $n \geq n_0$ may be made independently of  $f \in \initcond_{\ovbar\coninit}$, then the collection of measures
$\big\{ \nu_{n;(*:f,0)}^{([-1,1],1)}: n \in \N \big\}$ is $\initcond_{\ovbar\coninit}$-uniformly tight, where  the indexing variable is $f \in  \initcond_{\ovbar\coninit}$.
A little work is needed to use the proof of~\cite[Theorem~$8.2$]{Billingsley} to establish this last assertion.
We need to understand that, if the two bounds~(\ref{e.onepoint}) and~(\ref{e.equictyclaim}) hold whenever $n \geq n_0(\e,\eta)$, we are able to assert that the same bounds also hold whenever $n \geq n_0$ where the new selection of $n_0$ is made merely as a function of $\ovbar\coninit$. For this, what is needed is that, for a {\em given} value of $n$ that exceeds an $\ovbar\coninit$-determined level, these two bounds may be asserted with the parameters $K$ and $\eta$ being selected independently of $f \in \initcond_{\ovbar\coninit}$. We omit this fact's proof, but mention that the essence of the derivation lies in the argument for Lemma~\ref{l.basic}(1) and~(3).

Theorem~\ref{t.wlp.one}(2)
thus follows from equicontinuity Lemma~\ref{l.equicty} and the next {\em uniform boundedness}~lemma. \qed

\vspace{-2mm}
 Let  $f \in  \initcond_{\ovbar\coninit}$. For $K > 0$, define the 
 event  
 $$
  \mathsf{UnifBd}_{n;(*:f,0)}^{[-1,1]}(K) \, = \, \Big\{ \sup_{y \in [-1,1]} \big\vert  \weight_{n;(*:f,0)}^{(y,1)} \big\vert \leq K \Big\} \, .
 $$ 
\begin{lemma}\label{l.unifbd}
For any $\e > 0$ small enough, there exists $K = K(\e,\ovbar\coninit) > 0$ such that, for all $f \in \initcond_{\ovbar\coninit}$,
$$
 \PP \Big( \neg \,    \mathsf{UnifBd}_{n;(*:f,0)}^{[-1,1]}(K) \Big) \leq \e 
$$
whenever $n$ exceeds a certain constant multiple of $\big( \log \e^{-1} \big)^{12}$.
\end{lemma}
\vspace{-2mm}
\noindent{\bf Proof.}
We begin by arguing that, 
for $R \geq \coninit_2$,
\begin{eqnarray}
 & & \regfluc_{n;(*:f,0)}^{\big( \{-1,1 \}, 1 \big)}(R-1) \cap  \maxmin_{n;([-R,R],0)}^{([-1,1],1)}( R^2 ) \label{e.unifbdinc} \\
 & \subseteq &  \mathsf{UnifBd}_{n;(*:f,0)}^{[-1,1]}\Big( R^2 + 2^{-1/2} (R+1)^2 + \max \{ \coninit_3 , \coninit_1(1 + R) \} \Big) \, . \nonumber
\end{eqnarray}
Indeed, it was noted after~(\ref{e.equictyinc})  that $\weight_{n;(*:f,0)}^{(y,1)} = \sup \big\{ \weight_{n;(x,0)}^{(y,1)} + f(x) : x \in [-R,R] \big\}$
when the event $\regfluc_{n;(*:f,0)}^{\big( \{-1,1 \}, 1 \big)}(R-1)$ occurs;  since $R \geq \coninit_2$, $- \coninit_3 \leq \sup_{\vert x \vert \leq R} f(x)  \leq  \coninit_1(1 + R)$.
On the event $\maxmin_{n;([-R,R],0)}^{([-1,1],1)}( R^2 )$,  $\big\vert \weight_{n;(x,0)}^{(y,1)} \big\vert \leq R^2 + 2^{-1/2} (R+1)^2$ whenever $\vert x \vert \leq R$ and $\vert y \vert \leq 1$; this proves~(\ref{e.unifbdinc}).

Set $R = 1 + \Ctbd \big( \log \e^{-1} \big)^{1/3}$ as in the proof of Lemma~\ref{l.equicty}.
 We now apply 
Corollary~\ref{c.maxminweight}   with parameter settings ${\bf t_1} = 0$, ${\bf t_2} = 1$, ${\bf x} = -R$, ${\bf y} = -1$, ${\bf a} = \lceil 2R \rceil$, ${\bf b} = 2$ and ${\bf r} = R^2$ to find that 
$$
\PP \Big( \neg \, \maxmin_{n;([-R,R],0)}^{([-1,1],1)}( R^2 ) \Big) \leq  (2R+1) \cdot 400 C \exp \big\{ - \rsc_1 2^{-10} R^3 \big\}
$$
provided that $n$ exceeds an $\e$-determined level (which is of the order $\big( \log \e^{-1} \big)^{12}$, in order that the hypothesis ${\bf r} \leq 4 n^{1/{18}}$ be satisfied).   This upper bound is at most $\e/2$ for small enough $\e > 0$,
since $\rsc_1 2^{-10} \Ctbd^3 > 1$ holds in view of $\Ctbd \geq 12c^{-3}$, $c_1 \geq 2^{-5/2}c$ and $c \leq 1/2$.

Set $K = R^2 + 2^{-1/2} (R+1)^2 + \max \{ \coninit_3 , \coninit_1(1 + R) \}$. From~(\ref{e.unifbdinc}), we combine the last inference with~(\ref{e.regflucbound}) to obtain Lemma~\ref{l.unifbd}. \qed

\noindent{\bf Proof of Theorem~\ref{t.wlp.one}(3).}
 Let $\nu \in \wlp_{\ovbar{\coninit}}$. For a sequence $\big\{ f_n : n\in \N \big\}$ of elements of  $\initcond_{\ovbar\coninit}$, set $\nu_n = \nu_{n;(*:f_n,0)}^{([-1,1],1)}$. Then, for some such sequence of functions, and along a certain subsequence of~$n$, the sequence $\nu_n$ converges to $\nu$ weakly.


Recall that $r$ is a positive parameter that is at least $r_0$. For $\e > 0$ and $m \in \N$, the set 
$$
 \left\{ \, h \in \mc{C}: \sup_{\begin{subarray}{c} y,z \in [-1,1], \\
    2\exp \{-c' m^{1/12} \} < z - y < e^{-1}  \end{subarray}} \frac{\big\vert \, h(z) - h(y) \, \big\vert}{(z-y)^{1/2} \big( \log (z-y)^{-1} \big)^{2/3}} \, > \, r - \e \, \right\}  
$$
is open in $\mc{C}$. Recall the notation that $X$  is $\nu$-distributed. Applying the Portmanteau theorem along the convergent subsequence of $\nu_n$, we learn that 
$$
 \nu 
 \left( \,  \sup_{\begin{subarray}{c} y,z \in [-1,1], \\
    2\exp \{-c' m^{1/12} \} < z - y < e^{-1}  \end{subarray}} \frac{\big\vert \, X(z) - X(y) \, \big\vert}{(z-y)^{1/2} \big( \log (z-y)^{-1} \big)^{2/3}} \, > \, r - \e \, \right)  
$$
is at most the limit infimum along the concerned subsequence of $n \in \N$  of the left-hand side of~(\ref{e.nmodcon}) when $m$ replaces $n$ in the subscripted lower bound on $z - y$; $f_n$ replaces $f$; and $> r - \e$ replaces $\geq r$. 
We then apply Theorem~\ref{t.nmodcon}, consider $m \to \infty$ and then $\e \searrow 0$ to obtain Theorem~\ref{t.wlp.one}(3). \qed


\bibliographystyle{plain}

\bibliography{airy}

\end{document}